\documentclass[a4paper,twoside]{memoir}

\usepackage{myPackages}
\usepackage{myAssignments3}
\usepackage{myOperators}

\usepackage{fancyhdr,floatpag}
\floatpagestyle{fancy}

 \addto\captionsenglish{} 

\counterwithin{thm}{chapter}
\counterwithin{defn}{chapter}
\counterwithin{lem}{chapter}
\counterwithin{prop}{chapter}
\counterwithin{cor}{chapter}

\begin{document}

\includepdf[fitpaper=true]{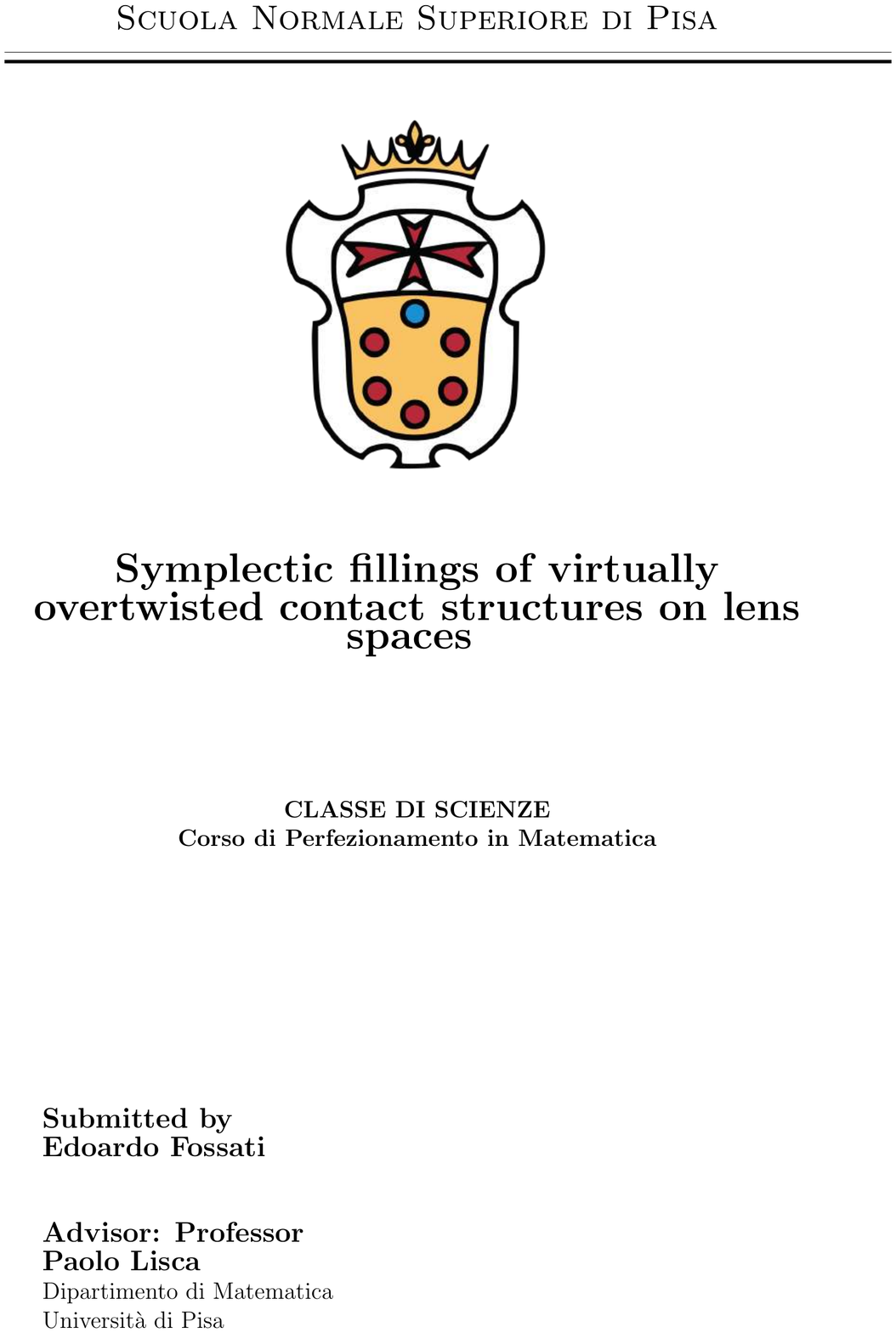}  

\clearpage{\pagestyle{empty}\cleardoublepage}
\thispagestyle{empty}

\begin{center}
\huge{\textbf{Abstract}}
\end{center}

\vspace{1cm}

Symplectic fillings of standard tight contact structures on lens spaces are understood and classified. The situation is different if one considers non-standard tight structures (i.e. those that are \emph{virtually overtwisted}), for which a classification scheme is still missing. In this work we use different approaches and employ various techniques to improve our knowledge of symplectic fillings of virtually overtwisted contact structures. 

\vspace{0.35cm}

We study curves configurations on surfaces to solve the problem in the case of a specific family of lens spaces. Then we give general constraints on the topology of Stein fillings of any lens space by looking at algebraic properties of integer lattices and at geometric slicing of solid tori. Furthermore, we try to place these manifolds in the context of algebraic geometry, in order to determine whether Stein fillings can be realized as Milnor fibers of hypersurfce singularities, finding a series of necessary conditions for this to happen. In the concluding part of the thesis, we focus on the connections between planar contact 3-manifolds and the theory of Artin presentations.

\clearpage{\pagestyle{empty}\cleardoublepage}
\thispagestyle{empty}

\tableofcontents*

\clearpage{\pagestyle{empty}\cleardoublepage}
\thispagestyle{empty}

\begin{center}
\huge{\textbf{Introduction}}
\end{center}
\addcontentsline{toc}{chapter}{Introduction}
 
\vspace{1cm}

In this thesis we study some topological properties of symplectic 4-manifolds with non-empty (and prescribed) boundaries. If the boundary of a symplectic manifold $(X,\omega)$ is $\omega$-convex (i.e. $\mathcal{L}_{\nu}\omega=\omega$ for an outer vector field $\nu$ defined in a neighborhood of $\partial X$), then the restriction of $\omega$ to the 3-dimensional boundary is a contact 1-form $\alpha=\omega(\nu,-)$: this is the basic fact that ties symplectic and contact topology together. 

A contact structure $\xi$ is a nowhere integrable planes distribution. The first distinction within contact geometry is between overtwisted and tight structures: $\xi$ is an \emph{overtwisted} contact structure on $Y$ if there exists an embedded disk $D\hookrightarrow Y$ such that $\xi$ agrees with $TD$ along the boundary $\partial D$. If such disk does not exist, the structure is called \emph{tight}. The classification of overtwisted contact structures on a 3-manifold is essentially reduced to a problem in homotopy theory: there is a unique overtwisted contact structure in every homotopy class of oriented plane fields, see \cite{eliashover}. However, all the contact structures that arise from a symplectic filling (see Definition \ref{weakfilling}) as $\xi=\ker (\omega(\nu,-))$ are tight \cite{eligrom}, and their properties are more geometric in nature. As a drawback, this kind of geometric approach that uses differential forms is not ideal for topologists. Thanks to the works of Donaldson \cite{donaldsonLefschetz} and Gompf \cite{gompfLefschetz}, the way of looking at symplectic 4-manifolds has changed in favor of a more topological point of view: up to certain conditions, symplectic 4-manifolds are the same as Lefschetz fibrations, in the sense that such a fibration supports a symplectic form and, vice versa, any such manifold can be given the structure of a Lefschetz fibration.

Of particular importance in the context of contact 3-manifolds are those Lefschetz fibrations over the disk, with bounded fibers and positive monodromy: this is one way of characterizing Stein domains \cite{loipiergallini}, \cite{akbulut}. A \emph{Stein domain} can be defined as a smooth submersion $X\to D^2$ whose general fiber is a surface with boundary, away from a finite set of critical points where the local behavior is modeled by the holomorphic map $(z,w)\mapsto z^2+w^2$. To each singular value corresponds a singular fiber, obtained by collapsing to a point a simple closed curve on a neighboring smooth surface $\Sigma$, see \cite{gs}. The advantage of this correspondence is that all the information regarding the underlying diffeomorphism type of the Stein domain is encoded by a chosen positive factorization of the monodromy: the (ordered) collection of those simple curves on the general fiber is the supporting set of as many positive (i.e. right) Dehn twists, whose product $\t_1 \t_2\cdots\t_n$ is called the \emph{monodromy}, usually indicated with the letter $\varphi$. 

The element $\varphi$ of the mapping class group $\Gamma(\S)$ is not enough, alone, to describe the symplectic (in fact Stein) 4-manifold $X$ we started from, but a positive factorization has to be specified. This is the key distinction between the 4-dimensional description of a symplectic manifold and the 3-dimensional description of a contact manifold: the pair $(\Sigma,\varphi)$ constitutes a 3-dimensional open book decomposition, and identifies a 3-manifold $Y$ together with a contact structure $\xi$ on it, which is said to be supported by the open book itself. The resulting contact manifold $(Y,\xi)$ in unique up to contactomorphism \cite{giroux}. The pair $(Y,\xi)$ is determined by the element $\varphi\in\Gamma(\Sigma)$ and is independent of the chosen factorization: a different positive factorization of $\varphi$ would describe a different symplectic 4-manifold who has, nevertheless, the same contact boundary $(Y,\xi)$. 

\noindent In synthesis, open book decompositions are the boundaries of Lefschetz fibrations, and by running through the possible positive factorizations of the monodromy it is possible to find all the symplectic 4-manifold with a prescribed contact boundary (this holds when the genus of $\S$ is zero thanks to the work of Wendl \cite{wendl}, otherwise it is more complicated). This is what is usually referred to as the problem of studying symplectic or Stein fillings of a given contact 3-manifold.

\vspace{0.5cm}

This dissertation is about Stein fillings of certain contact structures on lens spaces. These structures are understood and classified in the work of Honda \cite{honda}, which is therefore the starting point of the theory.

We review some of the basic facts about contact and symplectic topology in Chapter \ref{generalities}, which serves as an introduction for all of the remaining chapters. These are independent one from the other, and can be read separately. In Chapter \ref{classificationfillings} we focus our attention on a specific family of lens spaces, i.e. those which arise via surgery on the Hopf link, and we classify their Stein fillings. We will prove the following:

\begin{thm*} Let $L$ be the lens space resulting from Dehn surgery on the Hopf link with framing $-a_1$ and $-a_2$, with $a_1,a_2\geq 2$. Let $\xi_{vo}$ be a virtually overtwisted contact structure on $L$. 
Then $(L,\xi_{vo})$ has:
\begin{itemize}
\item a unique (up to diffeomorpism) Stein filling if $a_1\neq 4\neq a_2$;
\item two homeomorphism classes of Stein fillings, distinguished by the second Betti number $b_2$, if at least one of $a_1$ and $a_2$ is equal to 4 and the corresponding rotation number is $\pm 2$. Moreover, the diffeomorphism type of the Stein filling with bigger $b_2$ is unique. If the rotation number is not $\pm 2$, then we have again a unique filling.
\end{itemize}
\end{thm*}

Chapter \ref{topologicalconstraints} is dedicated to general restrictions on the topology of minimal fillings of lens spaces, such as the Euler characteristic and the fundamental group. Among various results, we will show:

\begin{thm*} 
Let $\xi$ be any tight contact structure on $L(p,q)$. Let $W$ be a minimal symplectic filling of $L(p,q)$ and let $l=\length(p/q)$. Then $\chi(W)\leq 1+l$.
\end{thm*}

\begin{thm*}
Let $W$ be a symplectic filling of $(L(p,q),\xi_{vo})$, with $p=m^2$ and $q=mk-1$, for some $m>k>0$ and $(m,k)=1$. Then $\chi(W)\geq 2$.
\end{thm*}

Then, inspired by an open question in the book \cite{nemethi}, we look for necessary conditions for realizing lens spaces as boundary of the Milnor fiber of a complex hypersurface singularity: Chapter \ref{boundarymilnor} deals with this problem, after recalling some algebraic geometry terminology. The main result from this chapter is the following theorem: 

\begin{thm*}
Let $\xi_{vo}$ be a virtually overtwisted structure on $L(p,q)$. If we are in one of the cases below, then $(L(p,q),\xi_{vo})$ is not the boundary of the Milnor fiber of any complex hypersurface singularity:
\begin{itemize}
	\item[a)] $p/q=[a_1,a_2,\ldots,a_n]$ and $a_i$ is odd for some $i$;
	\item[b)] $p/q=[2x_1,2x_2]$;
	\item[c)] $p/q=[2x_1,2x_2,\ldots,2x_n]$, with $x_i>1$ for every $i$ ($n\geq 3$) and either:
\begin{itemize}
		\item[i)] $q^2\not\equiv 1 \pmod{p}$ or
		\item[ii)] $q^2\equiv 1 \pmod{p}$ and $n$ is even.
\end{itemize}
\end{itemize}
\end{thm*}

Finally, the concluding part of the thesis, Chapter \ref{artinandcontact}, looks for connections between contact geometry and the theory of Artin presentations.

\paragraph{Acknowledgments} 

First and foremost, I wish to thank my advisor, professor Paolo Lisca, for introducing me to the topic, for all the time he dedicated to me and for all the support during these three years in Pisa. He was able to suggest me several interesting problems, and then guide me towards the solutions of a good portion of them.
Many thanks are due to professor Andr\'as Stipsicz, who taught me a lot of beautiful mathematics and was always enthusiast to share his knowledge with me. I appreciated every conversation I had with him, without which my whole PhD experience would have been less amusing and pleasant. Their mentoring skills made this thesis work possible.

I am also grateful to those, in particular to András Némethi, who have raised some of the questions that I tried to answer in this dissertation, and to those who have dedicated time fixing small (and less small) details, especially Marco Golla.
My gratitude goes as well to all the fantastic colleagues that I met in Pisa and Budapest, who helped me carrying out my research work daily, among whom Giulio, Andrea, Marco, Carlo, Antonio, Fabio, Viktória and Kyle. Lastly, I thank the Rényi Institute of Mathematics for the hospitality, where part of this work has been done.


\clearpage{\pagestyle{empty}\cleardoublepage}
\thispagestyle{empty}

\chapter{Background notions}\label{generalities}

\section{Generalities on contact 3-manifolds}

In this opening chapter, we recall some of the basic notions of contact and symplectic topology, following mainly the book of Özbağcı
 and Stipsicz \cite{ozbagci}.

\begin{defn}
A \emph{contact structure} on a 3-manifold $M$ is a nowhere integrable planes distribution $\xi$. If there exists an embedded disk $D\hookrightarrow M$ such that $\xi$ agrees with $TD$ along the boundary $\partial D$, then the contact structure $\xi$ is said to be \emph{overtwisted}, otherwise it is called \emph{tight}. 
\end{defn}

\begin{defn}
Two contact 3-manifolds $(Y_1,\xi_1)$ and $(Y_2,\xi_2)$ are \emph{contactomorphic} if there exists a diffeomorphism $f:Y_1\to Y_2$ such that $df(\xi_1)=\xi_2$. Two contact structures $\xi$ and $\xi'$ on a 3-manifold Y are \emph{isotopic} if there is a contactomorphism $f:(Y,\xi)\to (Y,\xi')$ which is isotopic to the identity.
\end{defn}

A contact structure $\xi$ on an oriented 3-manifold $Y$ can also be described as the kernel of a 1-form $\alpha$ such that $\alpha\wedge d\alpha \neq 0$. If the orientation of $Y$ coincides with the one given by $\alpha\wedge d\alpha$, then we say that $\xi$ is positive. This condition is independent of the choice of $\alpha$ with $\xi=\ker\alpha$. We will always assume that the contact structure we deal with are positive.

\begin{defn} 
The \emph{standard contact structure} on $\R^3$ is the kernel of the 1-form:
\[\alpha_{st}=dz+x\,dy. \]
This extends to a well defined contact structure $\xi_{st}=\ker \alpha_{st}$ on $S^3$, which we refer to as the \emph{standard contact structure} on $S^3$.
\end{defn}

\begin{defn} 
We say that a knot $K\subseteq Y$ is in \emph{Legendrian position} with respect to a contact structure $\xi$ on $Y$ if
\[T_x K\subseteq\xi_x,\]
for every $x\in K$. We simply refer to \emph{Legendrian knots and links} when $(Y,\xi)=(S^3,\xi_{st})$.
\end{defn}

For a Legendrian knot $K\subseteq (\R^3,\xi_{st})$ we define the two \emph{classical invariants} (invariants of the contact isotopy type, see \cite[Sections 4.1, 4.2]{ozbagci}):
\begin{itemize}
\item the \emph{Thurston-Bennequin number} $\tb(K)\in\Z$ is the measure of the contact framing (i.e. the vector field orthogonal to $\xi$ along $K$) with respect to the Seifert framing of $K$.
\item The \emph{rotation number} $\rot(K)\in\Z$ is the winding number of $TK$ calculated in any trivialization of $\xi$ along $K$. For this definition to make sense, we need to choose an orientation of $K$ (for a nullhomologous knot $K$ in a general 3-manifold we need to fix a Seifert surface for it, and in general $\rot(K)$ will depend on this choice).
\end{itemize}

By looking at the \emph{front projection} $p(K)$ of an oriented Legendrian knot $K\subseteq (R^3,\xi_{st})$ to the $yz$-plane, we can compute its Thurston-Bennequin number and rotation number in the following way:
\begin{align*}
\tb(K)= &w(p(K))-\frac{1}{2}(c_D+c_U),\\
\rot(K)= &\frac{1}{2}(c_D-c_U),
\end{align*}
where $w(p(K))$ is the writhe of the knot diagram, and $c_D$ and $c_U$ are respectively the number of down and up cusps, see \cite[Lemmas 4.2.3, 4.2.4]{ozbagci}. Following \cite[Section 11]{ozbagci}, we now define an operation which is an extension to Legendrian knots of the usual Dehn surgery, and has the effect of producing a new contact manifold out of a Legendrian knot $K$ inside $(Y,\xi)$. 

\begin{defn}
\emph{Contact $(\pm 1)$-surgery} on $K\subseteq (Y,\xi)$ is an integral Dehn surgery along $K$ with framing given by its contact framing $\pm 1$. The contact structure on $Y\smallsetminus \nu(K)$ extends to the solid torus $S^1\times D^2$ glued in while performing the surgery, which is endowed with the unique tight contact structure that makes the boundary convex (see Definition \ref{convexdivid}) and whose dividing curves on the boundary agree with the ones coming from the exterior of the knot, see \cite[Section 11.2]{ozbagci}. We refer to contact $(-1)$-surgery as \emph{Legendrian surgery}.
\end{defn}


\begin{thm*}[\cite{dinggeiges}]
Every closed contact 3-manifold is the result of contact $(\pm 1)$-surgery on a Legendrian link in the standard $S^3$.
\end{thm*}

Contact geometry has become popular among topologists thanks to the work of Giroux \cite{giroux}. Open book decompositions can be used to study contact 3-manifolds in the way that we now describe. 
Given a fibered link $L$ inside a (compact, connected, oriented) 3-manifold $Y$, we look at the fibration structure of the complement
\[f:Y\smallsetminus\nu(L)\to S^1,\]
whose fiber $\Sigma$ is a surface with boundary. The link $L$ is called the \emph{binding} of the open book decomposition, and $\Sigma$ is the \emph{page}. This structure can be also described abstractly as a pair $(\Sigma,\varphi)$, where $\varphi$ is the element of the mapping class group $\Gamma(\Sigma)$ which corresponds to the monodromy of the previous fibration. The 3-manifold $Y$ can be recovered by capping off the mapping torus $M_{\varphi}$ with solid tori, one for every boundary component (glued identifying the boundary of the meridian disk with the $S^1$-factor of $M_{\varphi}$, so that the core curves of the $D^2\times S^1$'s will recover the binding $L$). 

\begin{defn} 
A contact structure $\xi$ on $Y$ is \emph{supported by the open book} with page $\Sigma$ and binding $L$ if $\xi$ can be expressed as the kernel of a contact 1-form $\alpha$ such that $L$ is transverse to the contact planes, $d\alpha$ is a volume form on each page and $\alpha(L)>0$ with the orientation of $L$ induced by the pages. In particular, $(Y,\xi)$ is said to be \emph{planar} if $Y$ admits an open book decomposition with planar pages supporting $\xi$.
\end{defn}

In the work \cite{thurwink}, the authors start from such a decomposition to show that any 3-manifold has a contact forms. Later, open book decompositions appear in Giroux correspondence up to the notion of \emph{stabilization}:
start with an abstract open book decomposition $(\Sigma,\varphi)$ and modify the page by attaching a 1-handle to $\Sigma$. Call $\Sigma'$ the new surface and consider a simple closed curve $\gamma$ intersecting the cocore of the new 1-handle once. We get another open book decomposition described by the pair $(\Sigma',\varphi\circ \t^{\pm}_{\gamma})$, which is called a positive/negative stabilization of $(\Sigma,\varphi)$. The inverse procedure is called positive/negative destabilization.
As explained in \cite[Chapter 9]{os}, a positive stabilization has no effect on the associated contact 3-manifold: if $(\Sigma,	\varphi)$ describes $(Y,\xi)$, then $(\Sigma',\varphi\circ \t_{\gamma})$ describes $(Y,\xi)\#(S^3,\xi_{st})$ which is contactomorphic to $(Y,\xi)$. On the other hand, a negative stabilization corresponds to a connected sum with an overtwisted $(S^3,\xi_{ot})$, hence the resulting contact structure will be overtwisted as well. We are now ready to state the following:

\begin{thm*}[Giroux]
On any closed oriented 3-manifold there is a one-to-one correspondence between the set of isotopy classes of contact structures and the open book decompositions up to positive stabilization/destabilization.
\end{thm*}

\begin{figure}[ht!]
\centering
\includegraphics[scale=0.4]{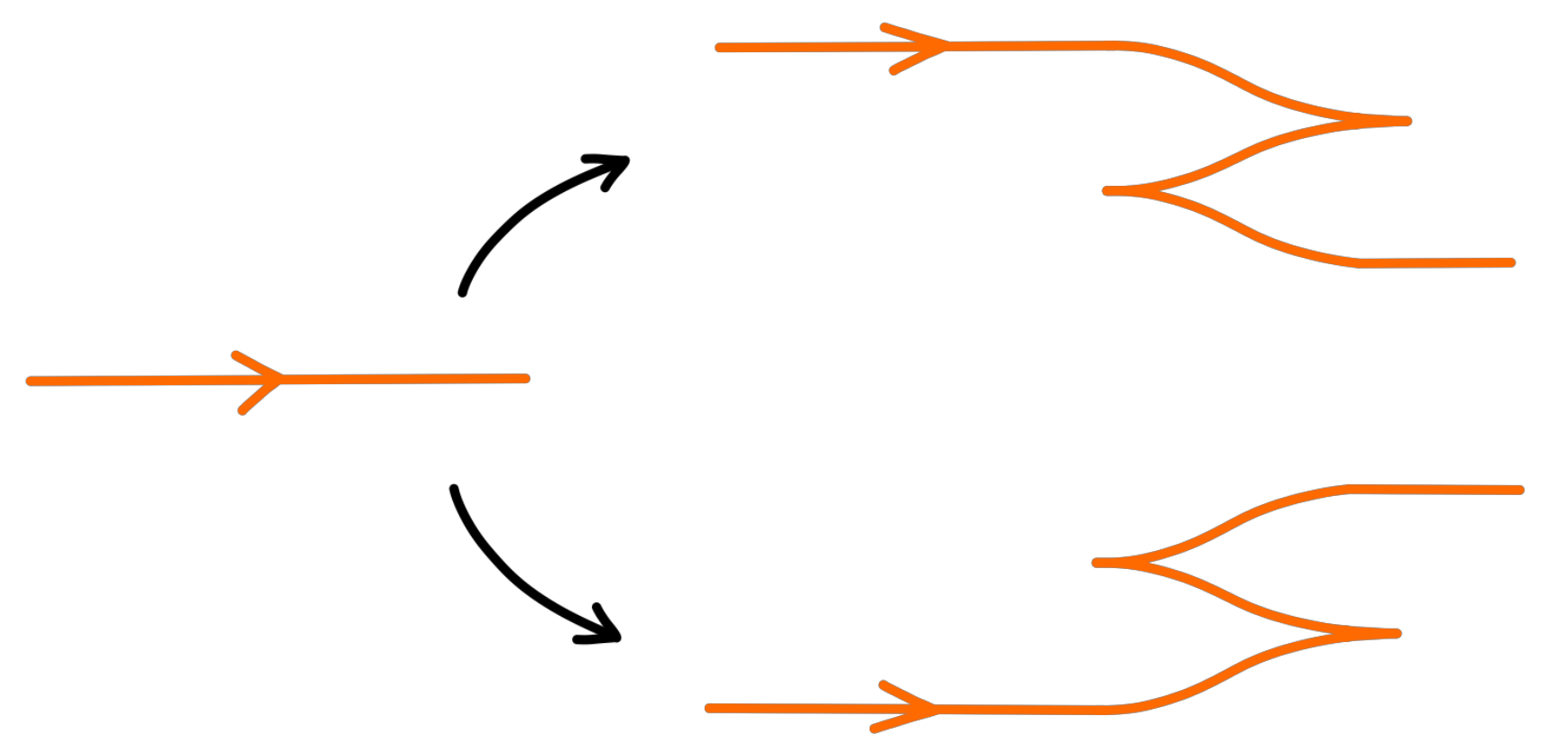}
\caption{Positive and negative stabilizations.}
\label{stabilization}
\end{figure}

Given a Legendrian knot, there is a procedure, called \emph{stabilization} (not to be confused with the stabilization of an open book decomposition), which reduces its Thurston-Bennequin number. After choosing an orientation of the knot, the stabilization can be either positive or negative, according to the effect it has on the rotation number (see Figure \ref{stabilization}). 

\noindent Recall that the rotation number of an oriented Legendrian knot in $(S^3,\xi_{st})$ can be computed in the front projection by the formula
\[\rot(K)=\frac{1}{2}(c_D-c_U),\]
In particular, a positive stabilization increases the rotation number by 1 and a negative stabilization decreases it by 1.

\begin{figure}[ht!]
\centering
\includegraphics[scale=0.5]{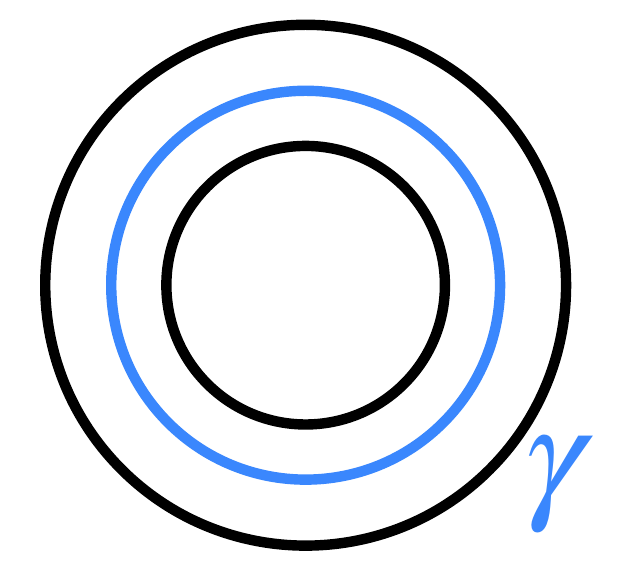}
\caption{Open book decomposition for $(S^3,\xi_{st})$.}
\label{S3page}
\end{figure}

As discussed above, we can present contact 3-manifolds in two different ways: via a Legendrian link in $(S^3,\xi_{st})$ on which to perform contact $(\pm 1)$-surgery, or via an open book decomposition $(\Sigma,\varphi)$. In what follows, we give an idea of how to move from a representation to the other, at least in the case when the link is as simple as possible, i.e. a 1-component unknot. We start with the open book decomposition of $(S^3,\xi_{st})$, which has an annular page and monodromy given by a positive Dehn twist along the core curve $\gamma$, see Figure \ref{S3page}. 

\begin{figure}[h!]
\centering
\begin{subfigure}[t]{.5\textwidth}
  \centering
  \includegraphics[scale=0.5]{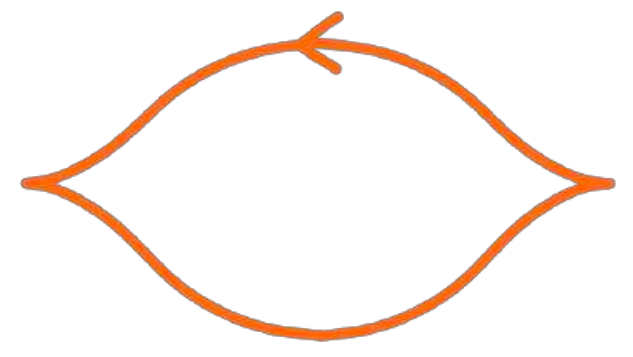}
  \caption{Legendrian unknot with $\tb=-1.$}
  \label{legunknot}
\end{subfigure}%
\begin{subfigure}[t]{.5\textwidth}
  \centering
 \includegraphics[scale=0.5]{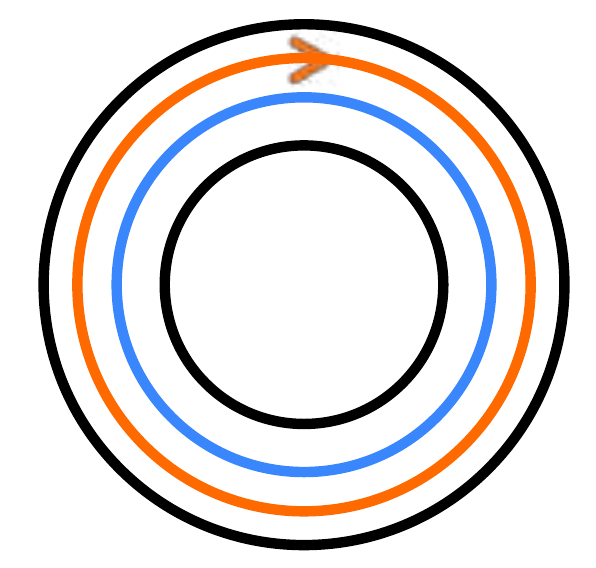}
  \caption{Resulting page.}
  \label{legunknotpage}
\end{subfigure}
\caption{Placing a knot on a page.}
\label{knotannulus}
\end{figure}

Now a Legendrian unknot $K$ with $\tb=-1$ in $(S^3,\xi_{st})$ can be placed on a page of the previous decomposition, simply by drawing a parallel copy of $\gamma$, see \cite{etnyre}. Performing Legendrian surgery on $K$ gives a contact manifold ($\R P^3$ with its unique tight structure) with supporting open book whose page is an annulus and whose monodromy is $\tau_{\gamma}\tau_{\gamma}$ (Figure \ref{knotannulus}). 

\begin{figure}[h!]
\centering
\begin{subfigure}[t]{.3\textwidth}
  \centering
  \includegraphics[scale=0.5]{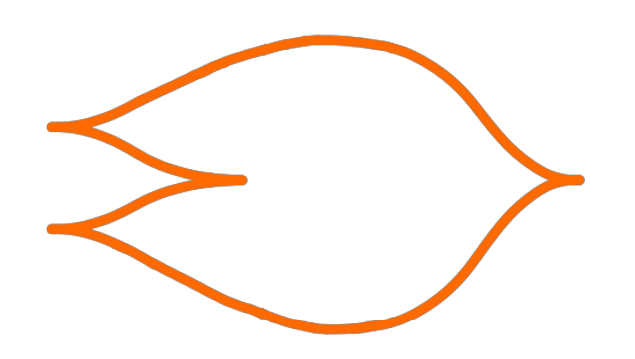}
  \caption{Positive stabilization.}
  \label{posstab1}
\end{subfigure}%
\begin{subfigure}[t]{.3\textwidth}
  \centering
 \includegraphics[scale=0.07]{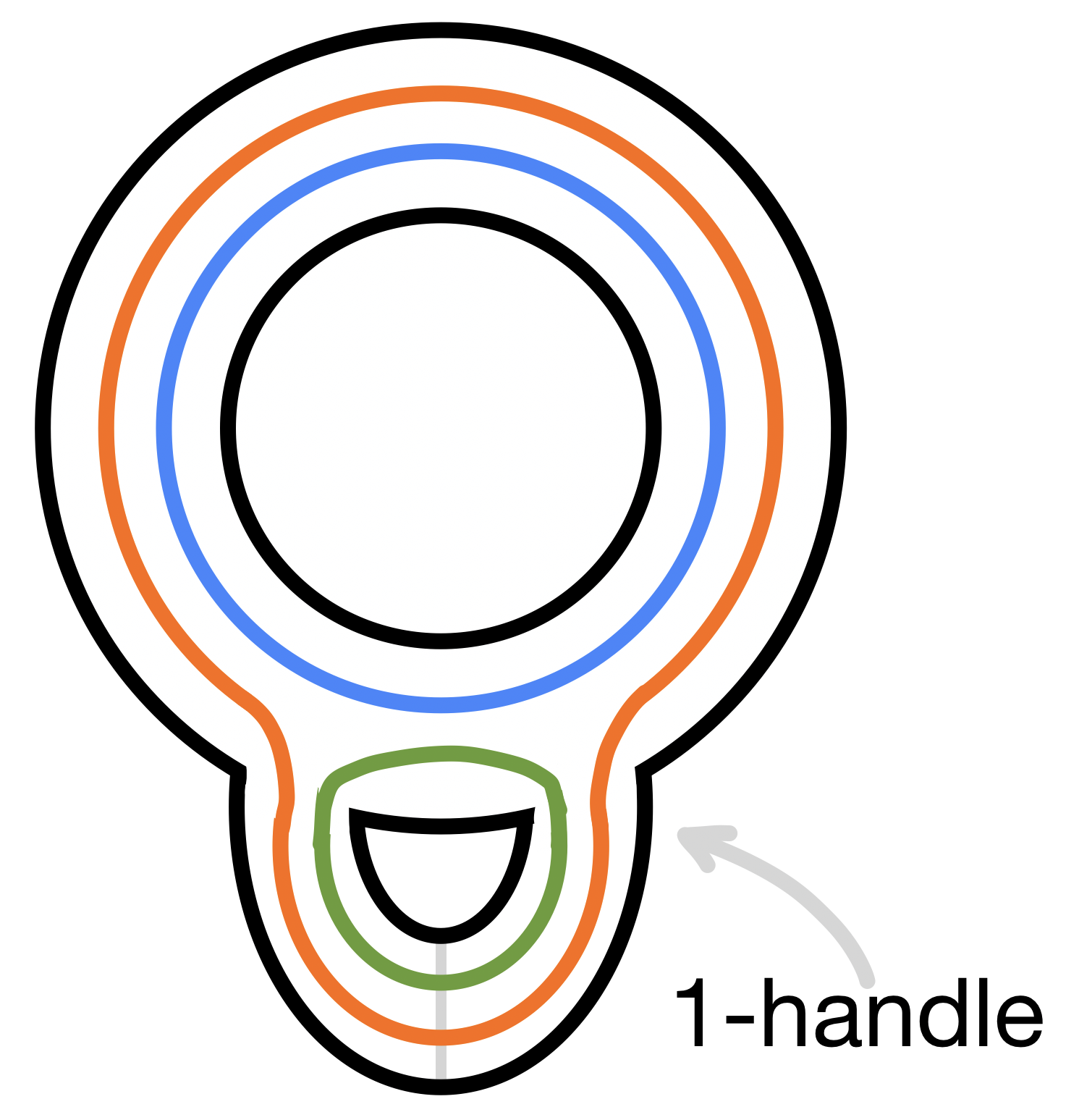}
  \caption{Effect on the open book.}
  \label{posstab2}
\end{subfigure}
\begin{subfigure}[t]{.3\textwidth}
  \centering
 \includegraphics[scale=0.5]{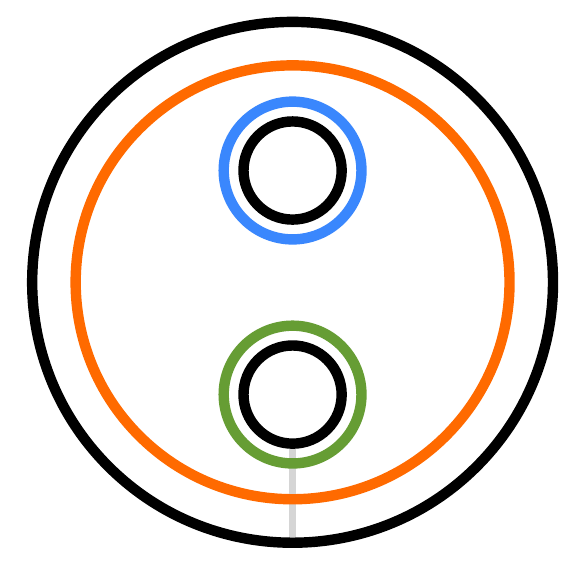}
  \caption{Resulting page.}
  \label{posstab3}
\end{subfigure}
\begin{subfigure}[t]{.3\textwidth}
  \centering
  \includegraphics[scale=0.5]{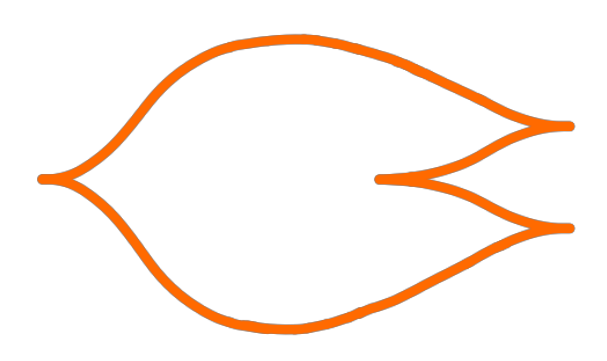}
  \caption{Negative stabilization.}
  \label{negstab1}
\end{subfigure}%
\begin{subfigure}[t]{.3\textwidth}
  \centering
 \includegraphics[scale=0.06]{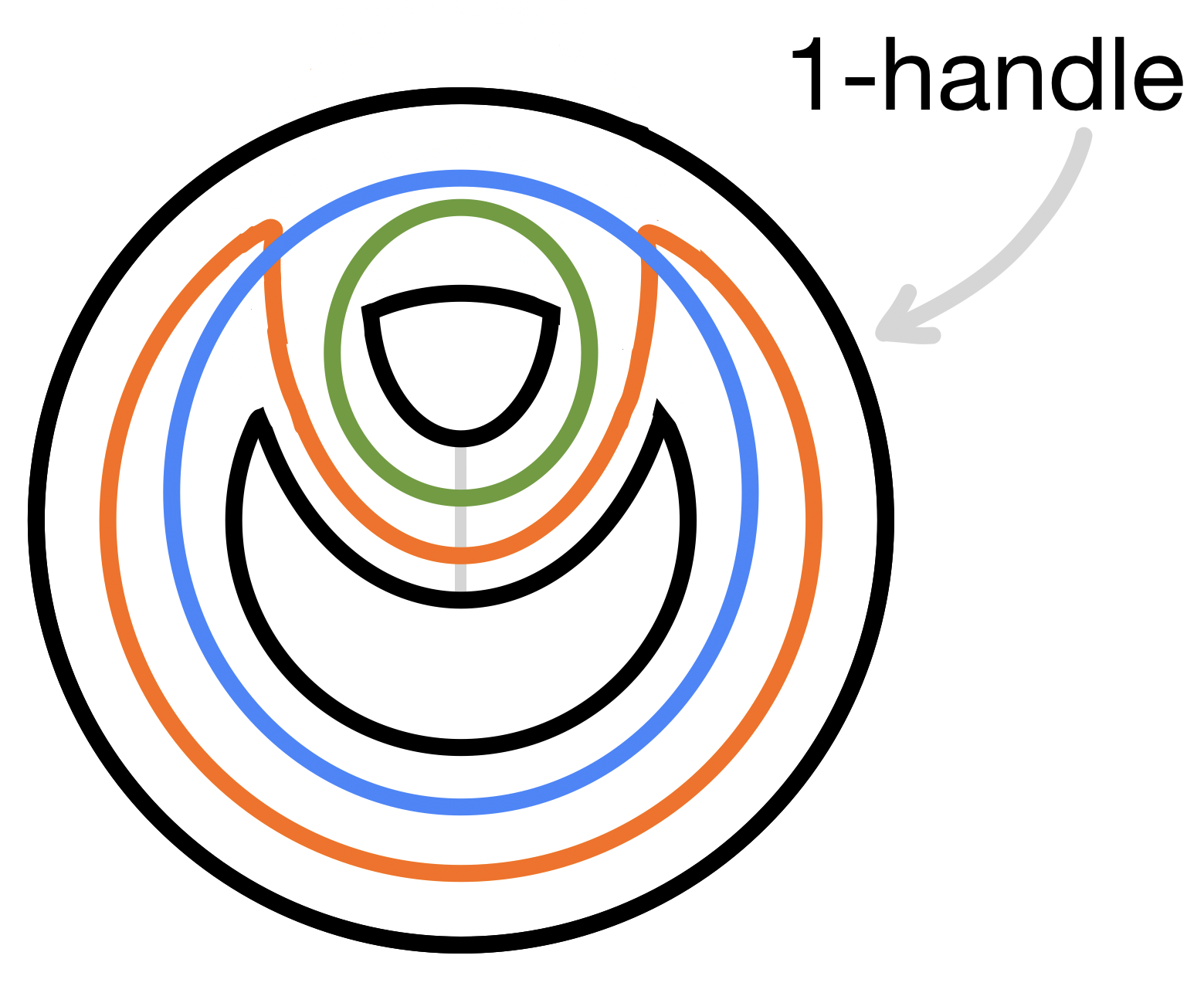}
  \caption{Effect on the open book.}
  \label{negstab2}
\end{subfigure}
\begin{subfigure}[t]{.3\textwidth}
  \centering
 \includegraphics[scale=0.5]{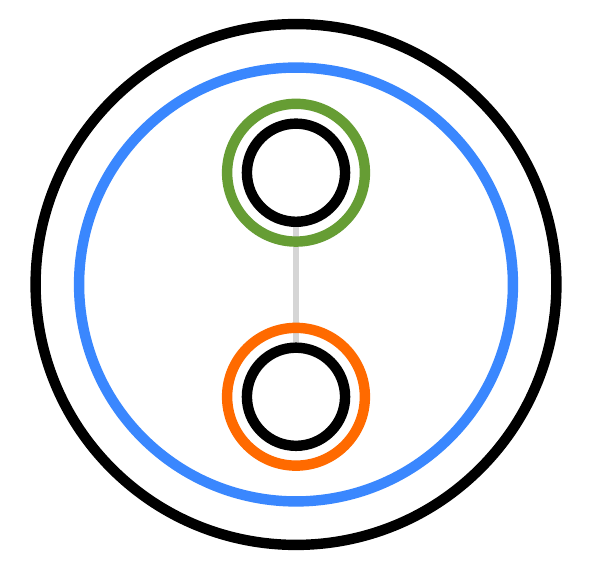}
  \caption{Resulting page.}
  \label{negstab3}
\end{subfigure}
\caption{Placing a stabilized knot on a page.}
\label{posnegstab}
\end{figure}

To reduce the Thurston-Bennequin number of the unknot, we add a positive or negative stabilization. Before drawing this Legendrian unknot $K'$ on a page of an open book decomposition of $(S^3,\xi_{st})$ we need to attach a 1-handle to the annulus and modify the monodromy by adding a positive Dehn twist along a curve intersecting the cocore of this new 1-handle once. We always attach the 1-handles on a connected component of the boundary, so that the total number of boundary components increases by one every time. Then we can draw $K'$ on the page by sliding the core curve of the annulus over the 1-handle. According to where we attach the 1-handle, we get a positive or negative stabilization (compare with Figure \ref{posnegstab} and with \cite{etnyre}): in order to distinguish between positive and negative, we need to pick an orientation of $K'$ on the page and we orient it in the clockwise direction.

For the rest of the work we will always use this \textbf{orientation convention}, specified for diagrams of Legendrian unknots: these knot diagrams in the front projection are oriented in the counter-clockwise direction,  dwhile curves on the planar page of an open book are oriented in the clockwise direction, as in Figure \ref{knotannulus}.

 
\section{The problem of understanding Stein fillings}\label{steinfillings}

The problem of classifying Stein (or more generally, symplectic) fillings of contact 3-manifolds has come to the attention of topologists since the pioneering work of Eliashberg \cite{eliashberg}. Showing that $S^3$ with its standard tight contact structure has a unique Stein filling, which is $(D^4,J_{st})$, is the starting point of this active research area. In the last years several other works have appeared on this subject. McDuff showed in \cite{mcduff} that $L(p,1)$, endowed with the standard tight contact structure, has a unique Stein filling when $p\neq 4$, and two different Stein fillings when $p=4$. Later, Lisca \cite{lisca} extended McDuff's results and gave a complete list of the Stein fillings of $(L(p,q),\xi_{st})$. Lens spaces surely represent a class of 3-manifolds for which many results are known: it comes from the fact that, in general, even trying to classify all the tight contact structures (up to isotopy) on a 3-manifold is hard, but at least on lens spaces this list is available thanks to the work of Honda \cite{honda} (see Section \ref{contactlens}). Partial results about fillings are available when one considers non-standard tight contact structures: Plamenevskaya and Van Horn-Morris \cite{plamenVHM} showed that the virtually overtwisted structures on $L(p,1)$ have a unique Stein filling. Another classification result about fillings of virtually overtwisted structures on certain families of lens spaces is due to Kaloti \cite{kaloti}. 

Following the notes of Özbağcı \cite{ozbagcilectures}, we recall some definitions and results in the theory of symplectic fillings.

\begin{defn}\label{weakfilling}
A contact 3-manifold $(Y,\xi)$ is said to be \emph{weakly symplectically fillable} if there is a compact symplectic 4-manifold $(W,\omega)$ such that $\p W = Y$ as oriented manifolds, and $\left.\omega\right|_{\xi}>0$. In this case we say that $(W,\omega)$ is a \emph{weak symplectic filling} of $(Y, \xi)$.
\end{defn}

\begin{defn} 
A contact 3-manifold $(Y,\xi)$ is said to be \emph{strongly symplectically fillable} if there is a compact symplectic 4-manifold $(W,\omega)$ such that $\p W = Y$ as oriented manifolds, $\omega$ is exact near the boundary and a primitive $\beta$ can be chosen in such a way that  $\ker(\left.\beta\right|_Y)=\xi$. In this case we say that $(W,\omega)$ is a \emph{strong symplectic filling} of $(Y, \xi)$.
\end{defn}

Before giving the definition of a Stein fillable contact 3-manifold, recall that a \emph{Stein manifold} is a complex manifold $(X,J)$ which admits a proper holomorphic embedding into some $\C^N$. Let $\phi:X\to\R$ be a function which is proper, bounded below and such that the associated Hermitian form $H_{\phi}$ is positive definite (see \cite[page 4]{ozbagcilectures}). Then, a \emph{Stein domain} $(W,J)$ is the preimage $\phi^{-1}(-\infty,t]$, for some regular value $t\in \R$, with the restricted complex structure $J=\left.J\right|_W$.

\begin{defn} 
A contact 3-manifold $(Y,\xi)$ is said to be \emph{Stein fillable} if there is a Stein domain $(W,J)$ such that $\p W = Y$ as oriented manifolds,
and $\xi$ is isotopic to $TY\cap JTY$. In this case we say that $(W,J)$ is a \emph{Stein filling} of $(Y, \xi)$.
\end{defn}

We often talk about \emph{minimal} fillings in the sense that the underlying smooth 4-manifold is minimal, i.e. it does not contain smoothly embedded spheres of self-intersection equal to $-1$. Stein domains are an example of minimal fillings of their contact boundary, see \cite[Theorem 10.3.1]{ozbagci}.

\begin{nrem}
By a result of Eliashberg and Gromov \cite{eligrom}, we know that weakly fillable contact structures are always tight. A Stein filling is in particular a strong filling where the symplectic form is exact (i.e. an \emph{exact symplectic filling}), and a strong filling is a weak filling. In the literature there are examples of contact 3-manifolds which are:
\begin{itemize}
\item tight but non fillable \cite{tightnonfill}, \cite{tightnonfill2};
\item weakly but not strongly fillable \cite{weaknonstrong}, \cite{weaknonstrong2};
\item strongly but not Stein fillable \cite{ghiggstrong}.
\end{itemize}
The situation is different when we deal with a contact structure on a rational homology sphere: in this case, a weak symplectic filling can be modified into a strong symplectic filling, see \cite{ohtaono} and \cite{eliashstrong}.
\end{nrem}

In this work we will study exclusively planar contact structures, for which the problem of understanding symplectic fillings is simplifyed by the following:

\begin{thm*}[\cite{niederkruger}]
If $(Y, \xi)$ is a planar contact 3-manifold, then every weak symplectic filling $(W,\omega)$ of $(Y,\xi)$ is symplectically deformation equivalent to a blow up of a Stein filling of $(Y,\xi)$.
\end{thm*}

\begin{defn} Let $X$ be a 4-manifold with boundary. A \emph{Lefschetz fibration} over a disk is a smooth map $f:X\to D^2$ subject to the following conditions:
\begin{itemize}
\item[1)] the map $f$ is a submersion away from a finite set of critical points $x_1,\ldots,x_n$ in the interior of $X$;
\item[2)] these critical points are mapped via $f$ to pairwise different points in $D^2$;
\item[3)] around each of these critical points there are complex coordinates $(z,w)$ such that $f$, in this local model, looks like $(z,w)\to z^2+w^2$.
\end{itemize}
The general fiber of $f$ is a smooth surface $\S_{g,b}$ of genus $g$ and with $b$ boundary components. Let $\{\lambda_1,\ldots, \lambda_n\}$ be an ordered set of oriented loops based at a fixed point in the interior of $D^2\smallsetminus \{f(x_1),\ldots,f(x_n)\}$, such that $\lambda_i$ encircles only the singular value $f(x_i)$, as in Figure \ref{lefdisk}.

\begin{figure}[ht!]
\centering
\includegraphics[scale=0.6]{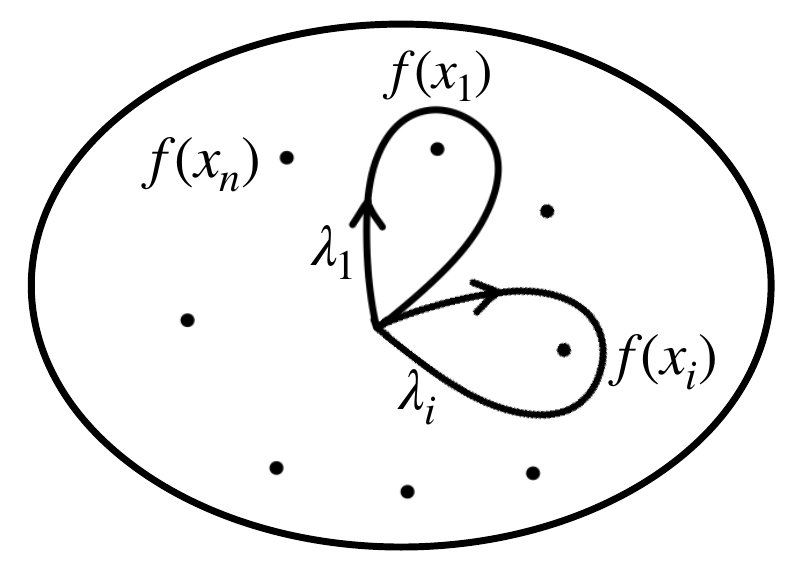}
\caption{Loops around the singular values.}
\label{lefdisk}
\end{figure}

These oriented loops freely generate the fundamental group $\pi_1(D^2\smallsetminus \{f(x_1),\ldots,f(x_n)\})$. We look at the monodromy of the fibration: consider the $\S_{g,b}$-bundle over the loop $\lambda_i$. This bundle has a monodromy, which is a positive or negative Dehn twist $\t_{\gamma_i}^{\pm}$ in the mapping class group $\Gamma_{g,b}=\Gamma(\S_{g,b})$ along a simple embedded curve $\gamma_i\subseteq \S_{g,b}$ (see \cite{gs}). If all the Dehn twists corresponding to the loops $\lambda_i$'s are positive, then the Lefschetz fibration $f$ is called \emph{positive}. In this thesis we will deal only with positive Lefschetz fibrations. The product of these ordered Dehn twists is the \emph{monodromy of $f$}:
\[\t_{\gamma_1}\t_{\gamma_2}\cdots \t_{\gamma_n}\in \Gamma_{g,b}.\]
If the curves $\gamma_i$ for $i=1,\ldots,n$ are all homologically non-trivial in the fiber, then the Lefschetz fibration is said to be \emph{allowable}.
\end{defn}

As already mentioned in the introduction, Loi-Piergallini \cite{loipiergallini} and independently Akbulut-Özbağcı
 \cite{akbulut} showed that Stein domains can be understood in terms of topological data: a positive allowable Lefschetz fibration over a disk has a Stein structure on its total space, and, vice versa, any Stein domain can be given such a fibration structure. 

Giroux correspondence allows one to represent a contact 3-manifold via a compatible open book decomposition. As a consequence, factorizing the monodromies of \emph{all} the compatible open book decompositions into products of positive Dehn twists is a (theoretical) solution to produce a complete list of Stein fillings for the corresponding contact 3-manifold. In the case of \emph{planar} open book decomposition, the factorization problem is easier thanks to the following theorem of Wendl:

\begin{thm*}[\cite{wendl}]
If a contact structure $\xi$ on a 3-manifold $Y$ is supported by an open book decomposition with planar page, then every strong symplectic filling of $(Y,\xi)$ is symplectic deformation equivalent to a blow-up of a positive allowable Lefschetx fibration compatible with the given open book.
\end{thm*}

\noindent Another important theorem regards the representation of Stein domain by means of surgery diagram:

\begin{thm*}[\cite{eliashbergtopological}, \cite{gompf}]
A smooth handlebody consisting of a 0-handle, some 1-handles and some 2-handles admits a Stein structure if the 2-handles are attached to the Stein domain $\hbox{\large$\natural$} (S^1 \times D^3)$ along Legendrian knots in $\hbox{\large$\#$} (S^1 \times S^2,\xi_{st})$ such that the attaching framing of each Legendrian knot is $-1$ relative to the contact framing, and $\xi_{st}$ refers to the standard tight contact structure on $\hbox{\large$\#$} (S^1 \times S^2)$. Conversely, any Stein domain admits such a handle decomposition.
\end{thm*}

\section{Contact structures on lens spaces}\label{contactlens}

\begin{figure}[ht!]
\centering
\includegraphics[scale=0.4]{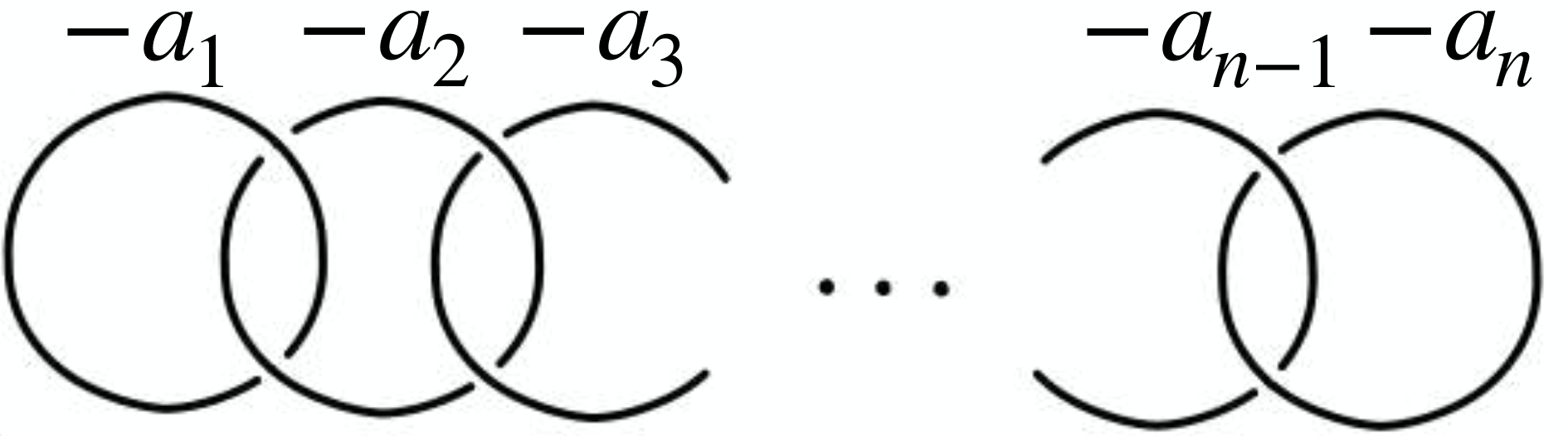}
\caption{Surgery link producing $L(p,q)$.}
\label{chain}
\end{figure}

\noindent Given a pair of coprime integers $p>q>0$, we consider the continued fraction expansion
\[\frac{p}{q}=[a_1,a_2,\ldots a_n]=a_1-\frac{1}{a_2-\frac{1}{\ddots -\frac{1}{a_n}}},\]
with $a_i\geq 2$ for every $i$. 
As a smooth oriented 3-manifold, $L(p,q)$ is the integral surgery on a chain of unknots with framings $-a_1,-a_2,\ldots -a_n$, see Figure \ref{chain}. 

\noindent To equip $L(p,q)$ with a tight contact structure, we put the link of Figure \ref{chain} into Legendrian position with respect to the standard tight contact structure of $S^3$, in order to form a linear chain of Legendrian unknots. We do this in such a way that the Thurston-Bennequin number of the $i^{th}$ component is $-a_i+1$. Remember that the convention that we use here is that  each knot in the front projection is oriented in the counter-clockwise direction. The information about the tight contact structure that we get by performing Legendrian surgery on this link is encoded by the position of the zig-zags. 

\begin{defn}\label{utvot}
A tight structure $\xi$ on $Y$ is called \emph{universally tight} if its pullback to the universal cover $\widetilde{Y}$ is tight. The tight structure $\xi$ is called \emph{virtually overtwisted} if its pullback to some finite cover $\widehat{Y}$ is overtwisted.
\end{defn}

\begin{nrem}
A consequence of the geometrization conjecture is that the fundamental group of any 3-manifold is \emph{residually finite} (i.e. any non trivial element is in the complement of a normal subgroup of finite index), and this implies that any tight contact structure is either universally tight or virtually overtwisted, see \cite{honda}.
\end{nrem}

\begin{figure}[h!]
\centering
\begin{subfigure}[t]{.5\textwidth}
  \centering
  \includegraphics[scale=0.5]{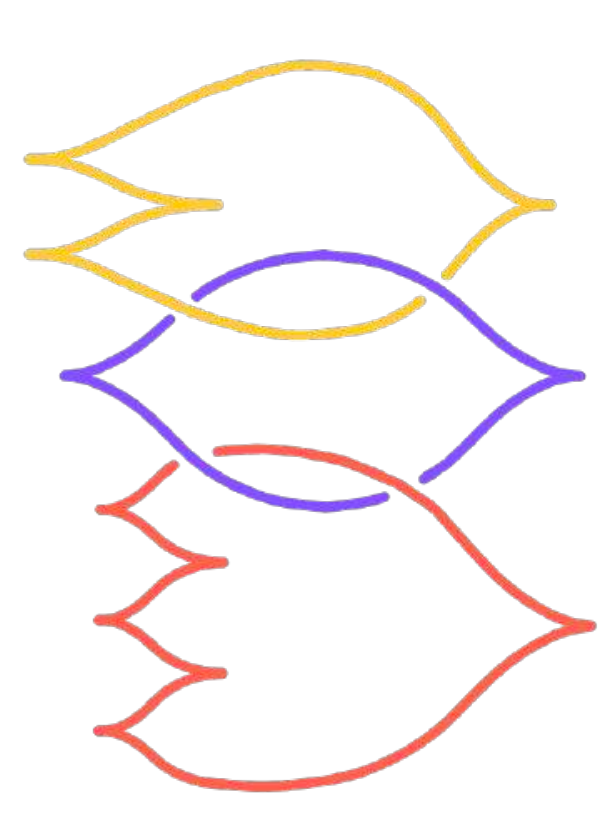}
  \caption{A universally tight structure.}
  \label{egut}
\end{subfigure}%
\begin{subfigure}[t]{.5\textwidth}
  \centering
 \includegraphics[scale=0.5]{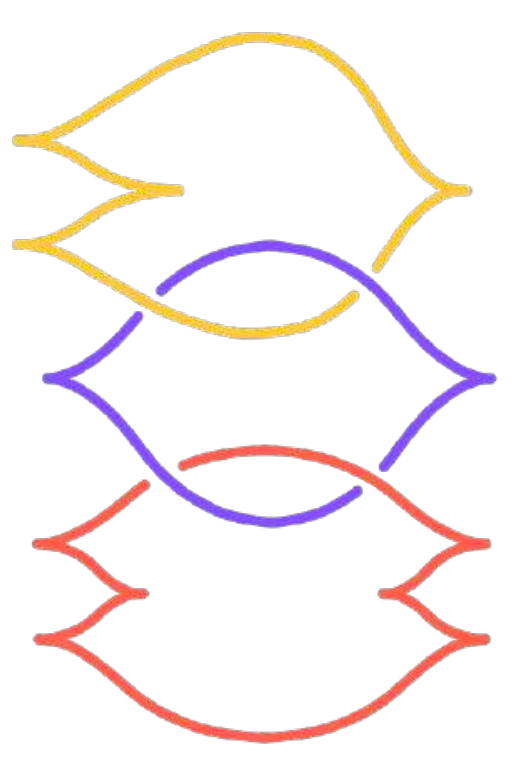}
  \caption{A virtually overtwisted structure.}
  \label{egvot}
\end{subfigure}
\caption{Comparing universally tight and virtually overtwisted contact structures.}
\label{egutvot}
\end{figure}

From the classification of tight contact structures on lens spaces (see \cite{honda}), if we look at a Legendrian realization of the link of Figure \ref{chain} given by a chain of Legendrian unknots, we can tell if the resulting contact structure will be universally tight or virtually overtwisted: if we have only stabilizations of the same type, either all positive or all negative, (i.e. zig-zags on the same side) then the contact structure will be universally tight, otherwise, if we have both positive and negative stabilizations, it will be virtually overtwisted (see Figure \ref{egutvot} for an example with $n=3$).

\begin{nrem}
\hypertarget{initialremark}{In} both cases, by attaching $n$ 4-dimensional 2-handles to $B^4$ with framing specified by the Thurston-Bennequin number of each component (decreased by 1), we get a Stein domain whose boundary has the contact structure specified by the Legendrian link. This shows that every tight structure on any lens space is Stein fillable. We will often write just "filling" to mean "Stein filling".
\end{nrem}

\chapter{Contact surgery on the Hopf link: classification of fillings}\label{classificationfillings}

Let $H\subseteq S^3$ be the two-components Hopf link. After choosing a Legendrian representative of $H$ with respect to the standard tight contact structure on $S^3$, we perform contact $(-1)$-surgery on the link itself. We get a lens space together with a tight contact structure on it, which depends on the chosen Legendrian representative. In this chapter, we classify its minimal symplectic fillings up to homeomorphism (and often up to diffeomoprhism), extending the results of \cite{lisca} which covers the case of universally tight structures, and the article of \cite{plamenVHM} which describes the fillings of $(L(p,1),\xi_{vo})$.

\begin{thm}\label{theorem1} Let $L$ be the lens space resulting from Dehn surgery on the Hopf link with framing $-a_1$ and $-a_2$, with $a_1,a_2\geq 2$. Let $\xi_{vo}$ be a virtually overtwisted contact structure on $L$. 
Then $(L,\xi_{vo})$ has:
\begin{itemize}
\item a unique (up to diffeomorpism) Stein filling if $a_1\neq 4\neq a_2$;
\item two homeomorphism classes of Stein fillings, distinguished by the second Betti number $b_2$, if at least one of $a_1$ and $a_2$ is equal to 4 and the corresponding rotation number is $\pm 2$. Moreover, the diffeomorphism type of the Stein filling with bigger $b_2$ is unique. If the rotation number is not $\pm 2$, then we have again a unique filling.
\end{itemize}
\end{thm}

\noindent We want to classify the fillings of the virtually overtwisted structures on $L(p,q)$ when 
\[\frac{p}{q}=[a_1,a_2].\]
Since all the tight contact structures on lens spaces are planar (see \cite[Theorem 3.3]{schonenberger}), we can apply Wendl's result on planar contact structures \cite{wendl} to our case. We will prove Theorem \ref{theorem1} by combining techniques coming from mapping class group theory with results by Schönenberger \cite{schonenberger}, Plamenevskaya-Van Horn-Morris \cite{plamenVHM}, Kaloti \cite{kaloti} and Menke \cite{menke}.

\section{Proof of the classification theorem} \label{classificationthm}

Recall from Section \ref{steinfillings} that a symplectic filling is called \emph{exact} if the symplectic form is exact. Stein domains are examples of exact fillings of their boundary. In \cite{menke} it is proved the following:

\begin{thm*}[\cite{menke}] Let $K$ be an oriented Legendrian knot in a contact 3-manifold $(M,\xi)$ and let $(M',\xi')$ be obtained from $(M,\xi)$ by Legendrian surgery on $S_+S_-(K)$, where $S_+$ and $S_-$ are positive and negative stabilizations, respectively. Then every exact filling of $(M',\xi')$ is obtained from an exact filling of $(M,\xi)$ by attaching a symplectic 2-handle along $S_+S_-(K)$.
\end{thm*}

\noindent Then we can derive an immediate corollary (everything is meant up to diffeomorphism).

\begin{cor*} Let $(L,\xi)$ be obtained by Legendrian surgery on the Hopf link.
\begin{itemize} 
\item[a)] Suppose that both components have been stabilized positively and negatively. Then $(L,\xi)$ has a unique Stein filling.
\item[b)] When just one component is positively and negatively stabilized, and the other one has topological framing different from $-4$, then again $(L,\xi)$ has a unique Stein filling.
\item[c)] If only one component is positively and negatively stabilized, and the other one has framing $-4$, then $(L,\xi)$ has two distinct fillings, coming from the two Stein fillings of $(L(4,1),\xi_{st})$.
\end{itemize}
\end{cor*}

\begin{figure}[h!]
  \centering
  \includegraphics[scale=0.5]{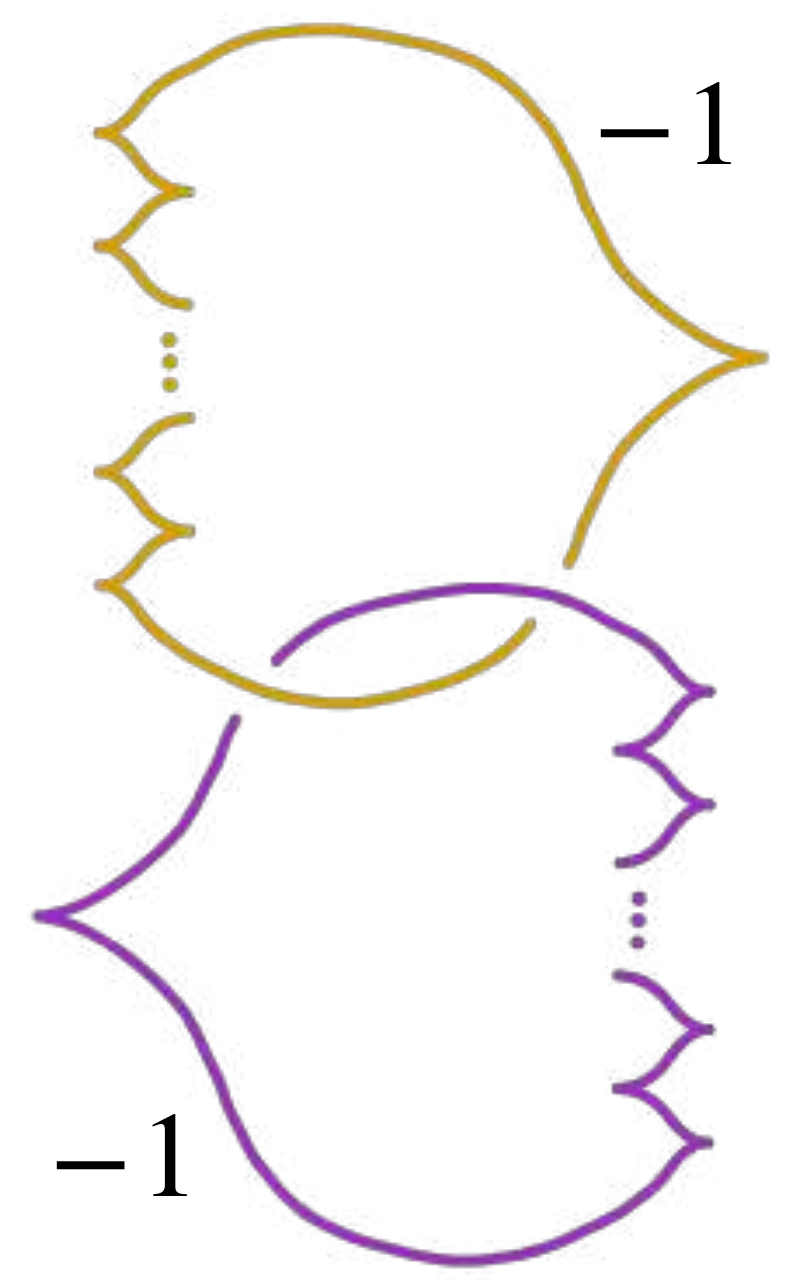}
  \caption{ }
  \label{hopf}
\end{figure}

\noindent The case that does not follow from the theorem of Menke, among the virtually overtwisted structures, is when one component of the link has all the stabilizations on one side and the other component on the opposite side. The rest of this chapter is devoted to cover this missing case. We will derive the classification of the fillings of the contact structures on $L(p,q)$ as in Figure \ref{hopf} in various steps. 

\subsection*{Step 1: place the link on a planar page}

We will focus on those lens spaces $L=L(p,q)$ with $\frac{p}{q}=[a_1,a_2]$. 
Starting from the link of Figure \ref{hopf}, we construct a planar open book decomposition of $(S^3,\xi_{st})$ so that the link itself can be placed on a page inducing the same framing as the contact framing: in order to do so, we first slide one component over the other. This changes the isotopy class of the link, but does not change the contact type of the 3-manifold obtained by Legendrian $(-1)$-surgery, see \cite{ding}. 

\begin{figure}[h!]  
  \centering
 \includegraphics[scale=0.5]{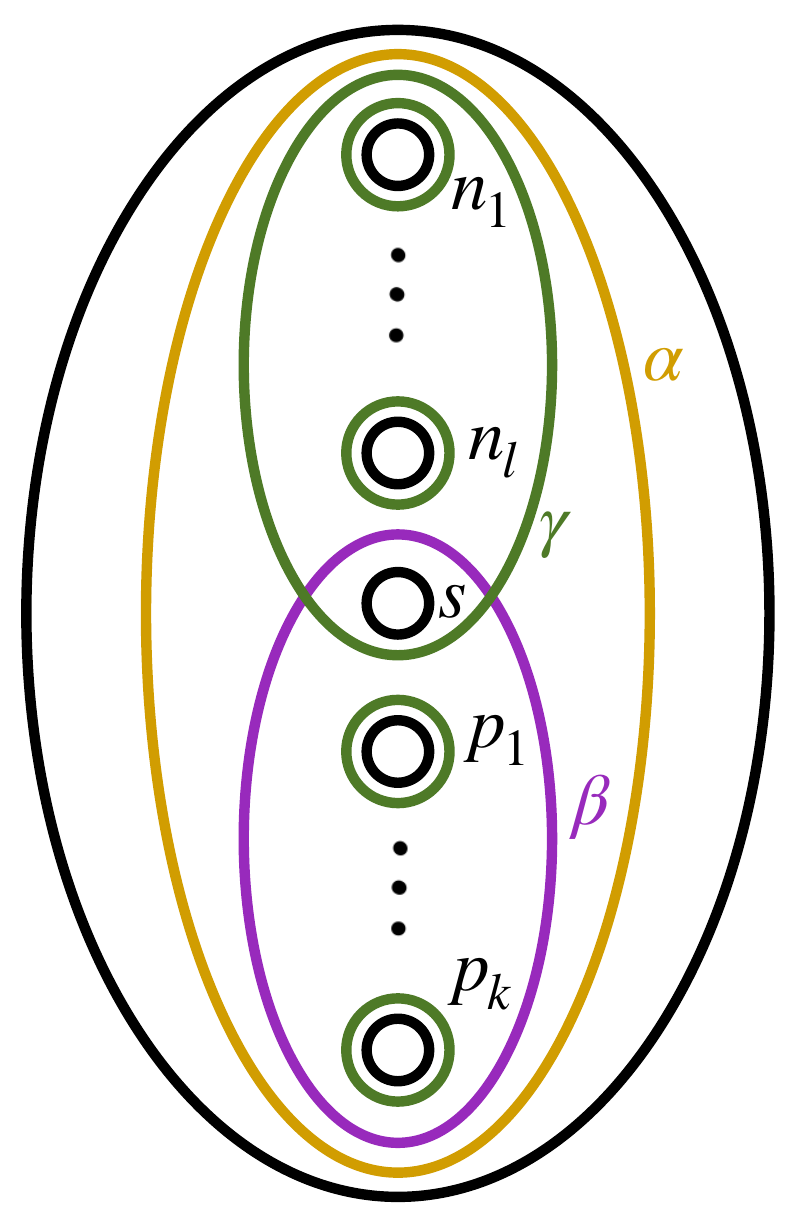}
\caption{Placing the (rolled-up) link on a page.}
\label{page}
\end{figure}

\begin{prop} Figure \ref{page} represents the page (and the monodromy) of an open book decomposition compatible with the contact structure obtained by performing Legendrian surgery on the link of Figure \ref{hopf}.
\end{prop}

\begin{prf}
Proceeding as described in Chapter \ref{generalities}, we start with the green curve $\gamma$ around the single inner hole $s$. Then the page is stabilized with the holes (and corresponding stabilizing green curves) $p_1,\ldots,p_k$ by attaching Hopf bands to the outer boundary component of the page. So nothing happens to the original curve $\gamma$, and we can place the yellow curve $\alpha$ going around every hole, corresponding to the yellow component of the link of Figure \ref{hopf}. The next step is to slide the purple curve around the yellow one in Figure \ref{hopf}, so that in the page we can place a parallel copy of $\alpha$, which is the curve $\beta$. Now $\beta$ has to be stabilized negatively $l$ times. To do this, we attach $l$ Hopf bands (with the corresponding stabilizing green curves) to the interior of the hole $s$, and we slide the purple curve on these 1-handels. The result is therefore what appears in Figure \ref{page}.
\end{prf}

\noindent Call $p_i$ the stabilizing curves corresponding to the positive stabilization of the knot, and $n_i$ the stabilizing curves corresponding to the negative ones.

This is how we can place the link (after sliding) on a planar page of an open book for $(S^3,\xi_{st})$, as described in \cite{schonenberger}. Here is the advantage of such a construction: performing Legendrian surgery on the link $\alpha\cup\beta$ gives the same contact 3-manifold as the one described by the abstract open book decomposition with that page and monodromy given by post-composing the original monodromy with $\t_{\alpha}\t_{\beta}$, that we eventually call $\varphi$. Hence, by Giroux's correspondence, we have that the contact type of $(L(p,q),\xi_{vo})$ is encoded in the pair $(\S,\varphi)$ with 
\begin{equation}\label{monod}
\varphi=\t_{\alpha}\t_{\beta}\t_{\gamma}\t_{p_1}\cdots \t_{p_k}\t_{n_1}\cdots \t_{n_l},
\end{equation}
which in turn already describes a Stein filling of $(L(p,q),\xi_{vo})$.

The theorem of Wendl (that we recalled in Section \ref{steinfillings}) implies that it is enough to find all the possible factorizations into positive Dehn twists of a \emph{given} planar monodromy in order to get all the Stein fillings of the contact manifold it represents.

\subsection*{Step 2: compute the possible homological configurations} \label{secmatrixmult}

In their work, Plamenevskaya and Van Horn-Morris \cite{plamenVHM} introduce the multiplicity and joint multiplicity of one or of a pair of holes for a given element in the mapping class group of a planar surface which is written as a product of positive Dehn twists. 

To define these numbers we need the "cap map", which is induced by capping off all but one (respectively two) interior component, while the outer boundary component is never capped. In the first case, we get the mapping class group of the annulus, which is isomorphic to $\Z$, generated by a positive Dehn twist along the core curve. In the second case, we get the mapping class group of a pair of pants, isomorphic to a free abelian group of rank 3, generated by the three positive Dehn twists around each boundary component. By projecting onto the third summand (that comes from the outer boundary component) we get the joint multiplicity around the other two components. We denote by $m(-)$ the multiplicity of a single hole, and by $m(-,-)$ the joint multiplicity of a pair of holes.

\begin{minipage}[c]{.40\textwidth}
\centering
\[\xymatrix{
\Gamma(\Sigma) \ar[r]^-{\mbox{cap}} \ar[rd]_-{m(-)} & \Gamma(\Sigma_{0,2}) \ar[d]^-{\simeq} \\
& \Z
}\]
\end{minipage}%
\hspace{10mm}%
\begin{minipage}[c]{.40\textwidth}
\centering
\[\xymatrix{
\Gamma(\Sigma) \ar[r]^-{\mbox{cap}} \ar@/_/[rdd] _-{m(-,-)} & \Gamma(\Sigma_{0,3}) \ar[d]^-{\simeq} \\
& \Z\oplus\Z\oplus\Z \ar[d]^-{\mbox{pr}_3}\\
& \Z
}\]
\end{minipage}

\begin{figure}[h!]
\centering
\begin{subfigure}[t]{.5\textwidth}
  \centering
  \includegraphics[scale=0.5]{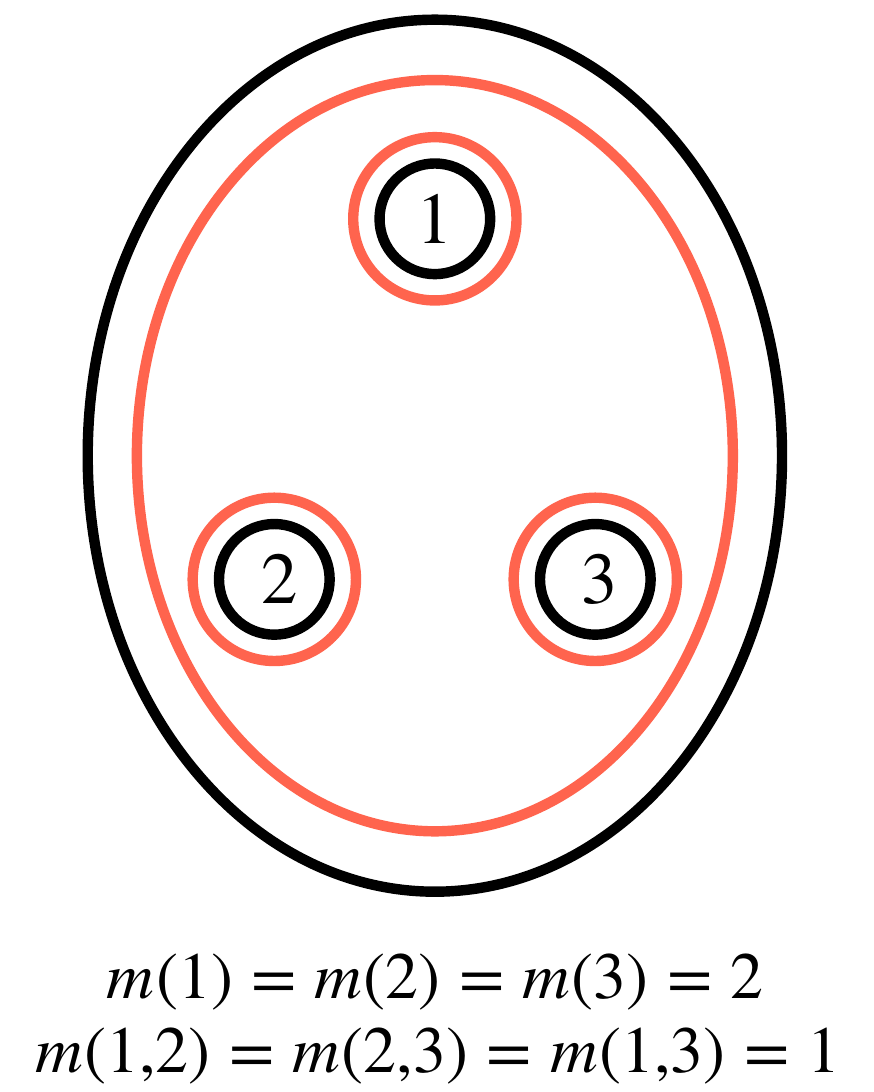}
\end{subfigure}%
\begin{subfigure}[t]{.5\textwidth}
  \centering
 \includegraphics[scale=0.5]{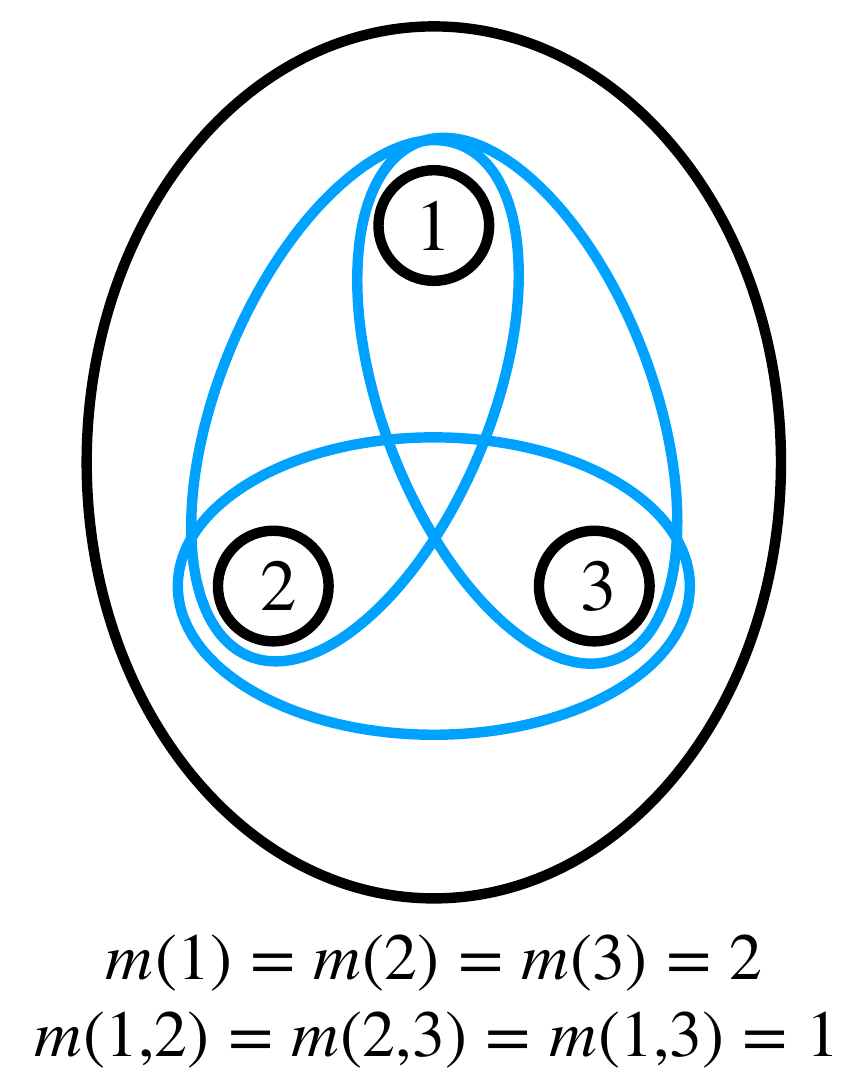}
\end{subfigure}
\caption{Lantern relation and multiplicities.}
\label{lanternmult}
\end{figure}

These multiplicities are independent of the positive factorization we chose for the monodromy.
The reason is that the lantern substitution preserves these numbers (compare with Figure \ref{lanternmult}) and the commutator relation does it too. Since these relations generate all relations in the planar mapping class group (see \cite{margalit}), we can, by introducing negative Dehn twists as well, apply them repeatedly until every curve encloses at most 2 holes.

We compute these numbers for the monodromy that we got from Step 1 by looking at Figure \ref{page}. In that figure, the holes called $n_i$'s are the ones corresponding to the negative stabilizations, while the $p_i$'s come from the positive ones and $s$ is the starting hole of the annulus, i.e. the one without a boundary-parallel curve around it. By applying the definition of the cap maps, we directly compute:

\begin{itemize}  
\item $m(n_i)=m(p_i)=3$
\item $m(s)=3$ 
\item $m(n_i,n_j)=m(p_i,p_j)=2$
\item $m(n_i,s)=m(p_i,s)=2$
\item $m(n_i,p_j)=1$
\end{itemize}

\noindent Starting from this collection of numbers, we try to reconstruct the homology classes of the curves appearing as support for the positive Dehn twists in the monodromy. We have to distinguish two cases: the first case is when $a_1\neq 4$ and $ a_2\neq 4$ (remember that $[a_1,a_2]$ is the continuous fraction expansion of $p/q$), the second case is when at least one of them is equal to 4. We postpone this second case to \hyperlink{thesentence}{Step 5}. 

\begin{prop}\label{conf}
Assume $a_1\neq 4$ and $ a_2\neq 4$. Then there is a unique homology configuration of curves whose associated multiplicities and double multiplicities are as above, and that gives the same monodromy we started from (Equation \eqref{monod}). In particular, the second Betti number of any Stein filling of the corresponding contact manifold is 2.
\end{prop}

\begin{prf}
We say that a simple closed curve is a \emph{multi-loop} if it encloses at least two holes, in order to distinguish it from a boundary-parallel curve.

Suppose we have another positive factorization of $\varphi$ (see Equation \ref{monod}), and call $\nu_{p_i}$ the number of multi-loops around the hole $p_i$ in this new factorization. 

From the fact that $m(p_i)=m(n_j)=m(s)=3$, we have the upper bound $\nu_{p_i},\nu_{n_j},\nu_s\leq 3$. Moreover, $m(p_i,s)=m(n_j,s)=2$ implies that $\nu_{p_i},\nu_{n_j},\nu_s\geq 2$. Hence we know
\[2\leq\nu_{p_i},\nu_{n_j},\nu_s\leq 3\]
for every $i$ and $j$. We present in details the longest combinatorial part:

\paragraph{Case $a_1,a_2>4$.}
We claim that, in the case $a_1>4$ (which translates into the fact that there are at least three positive holes), there is a multi-loop encircling all the positive holes and $s$. Since $m(p_1,s)=m(p_2,s)=2$ and $m(s)=3$, there must be a curve $a$ encircling $s,p_1,p_2$. If, by contradiction, there exists a hole $p_{\overline{\imath}}$ which is not encircled by $a$, then by the fact that $m(s,p_{\overline{\imath}})=2$, one finds other two curves $b$ and $c$ such that: $b$ encircles $s,p_{\overline{\imath}},p_2$ but not $p_1$, and $c$ encircles $s,p_{\overline{\imath}},p_1$ but not $p_2$. Since $m(s)=3$ and $m(s,n_1)=2$ we must have one of the three curve $a,b,c$ encircling $n_1$ as well, say it is $a$; but then it is impossible to obtain $m(n_1,p_{\overline{\imath}})=1$ without contradicting either $m(s)=3$ or $m(s,n_1)=2$.
This shows that there is a multi-loop $\beta'$ encircling (at least) all the positive holes and $s$. Similarly, there is a multi-loop $\gamma'$ encircling (at least) all the negative holes and $s$.

Now there are two cases:
\begin{itemize}
\item[1)] $\beta'$ and $\gamma'$ coincide, hence the multi-loop $\beta'=\gamma'$ encircles all the hole (i.e. it is parallel to the outer boundary component), and from now on it will be referred to as $\alpha'$;
\item[2)] $\beta'$ and $\gamma'$ are distinct in homology, hence one sees that $\beta'$ cannot encircle negative holes and $\gamma'$ cannot encircle positive holes (this uses the fact that $a_2>4$).
\end{itemize}
We are left to see how we can place the other curves in these two cases in order to get the multiplicities as in previous factorization (Equation \eqref{monod}):
\begin{itemize}
\item[1)] we just forget about $\alpha'$ by lowering all the multiplicities by 1 and then we do again the computation as above. We end up with a curve around the positive holes and $s$ (homologous to $\beta'$), and a curve around the negative holes and $s$ (homologous to $\gamma'$).
\item[2)] We do the same computation as above, by starting with a curve around $s$ and assuming that there is a hole not encircled by it. We get a contradiction with $m(s)=3$. This shows that there must be a curve encircling all the holes, hence parallel to the outer boundary component (i.e. homologous to $\alpha'$).
\end{itemize}
In both cases we get $\nu_s=3$ and $\nu_{p_i}=\nu_{n_j}=2$ for all $i,j$, and we end up with three multi-loops $\alpha'$, $\beta'$ and $\gamma'$, with $\alpha'$ going around all the holes, $\beta'$ around $\{s,n_1,\ldots,n_l\}$, $\gamma'$ around $\{s,p_1,\ldots, p_k\}$. In this way, all the conditions on the joint multiplicities are met, and we just need to add boundary-parallel loops around all the $p_i$'s and $n_j$'s in order to get $m(n_j)=m(p_i)=3$ as required.

\paragraph{Case $a_1=3$ or $a_2=3$.} This case is easier from the combinatorial point of view, and gives the same result. On the other hand, the case $a_1=4$ or $a_2=4$ gives rise to an extra configuration, as discussed later in Step 5.
\end{prf}

This tells us that the homology of any Stein filling for each one of these lens space is fixed. In particular, another factorization of $\varphi$ must be of the form
\[\varphi=\t_{\alpha'}\t_{\beta'}\t_{\gamma'}\t_{p_1}\cdots \t_{p_k}\t_{n_1}\cdots \t_{n_l},\]
where $\alpha',\beta'$ and $\gamma'$ are simple closed curves on $\S$ such that $[\alpha]=[\alpha'], [\beta]=[\beta'], [\gamma]=[\gamma']$ in $H_1(\S;\Z)$. Notice that we do not need to worry about the $p_i$'s and $n_j$'s because the fact that they homologically enclose just one hole implies that they are boundary-parallel, and so their homotopy (and therefore isotopy) class is already determined. Also the homotopy class of $\alpha'$ is determined (since it is boundary-parallel to the outer component) and so $\t_{\alpha}=\t_{\alpha'}$.
Therefore, if the configuration of curves we started from is like the one of Figure \ref{arcs1}, then we already know how to place the boundary-parallel curves appearing in any other factorization (compare with Figure \ref{arcs2}).

In light of this, using the previous factorization of $\varphi$, we see that all the $\t_{p_j}$'s and $\t_{n_i}$'s, together with $\t_{\alpha}$ and $\t_{\alpha'}$, cancel out, leaving us with:
\[\t_{\beta}\t_{\gamma}=\t_{\beta'}\t_{\gamma'}.\]
This relation holds in $\Gamma_{0,k+l+2}=\Gamma(\S_{0,k+l+2})$, and we are asking ourselves if there can be a pair of curves $\{\beta',\gamma'\}$ on $\S$ with the homological condition that $[\beta']=[\beta]$ and $[\gamma']=[\gamma]$ inside $H_1(\S_{0,k+l+2};\Z)$, and such that the product of the corresponding Dehn twists is isotopic to $\t_{\beta}\t_{\gamma}$. We will reduce the problem from $\S_{0,k+l+2}$ to $\S_{0,4}$.

\subsection*{Step 3: reduce the number of boundary components} 

\begin{defn} A diffeomorphism $f:\Sigma\to \Sigma$ which restricts to the identity on $\partial\Sigma$ is said to be \emph{right-veering} if $f(\eta)$ is to the right of $\eta$ at its starting point (or $f(\eta)=\eta$), for every properly embedded oriented arc $\eta\subseteq\Sigma$, isotoped (relative to the boundary) to minimize the number of intersection points of $\eta\cap f(\eta)$.
\end{defn}

\begin{figure}[ht!]
\centering
\begin{subfigure}[t]{.5\textwidth}
  \centering
  \includegraphics[scale=0.4]{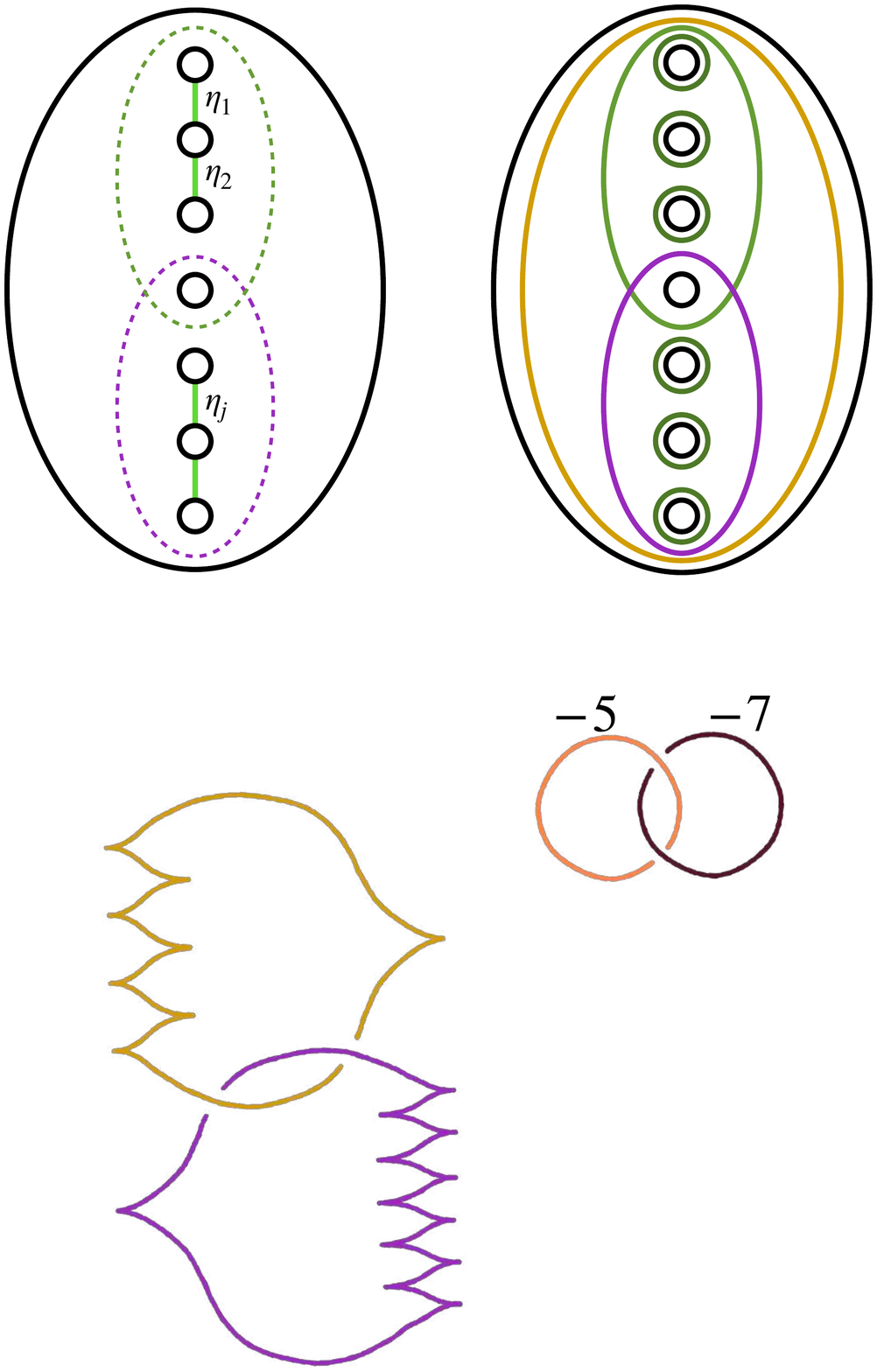}
\caption{Original configuration of curves.}
\label{arcs1}
\end{subfigure}%
\begin{subfigure}[t]{.5\textwidth}
\centering
 \includegraphics[scale=0.4]{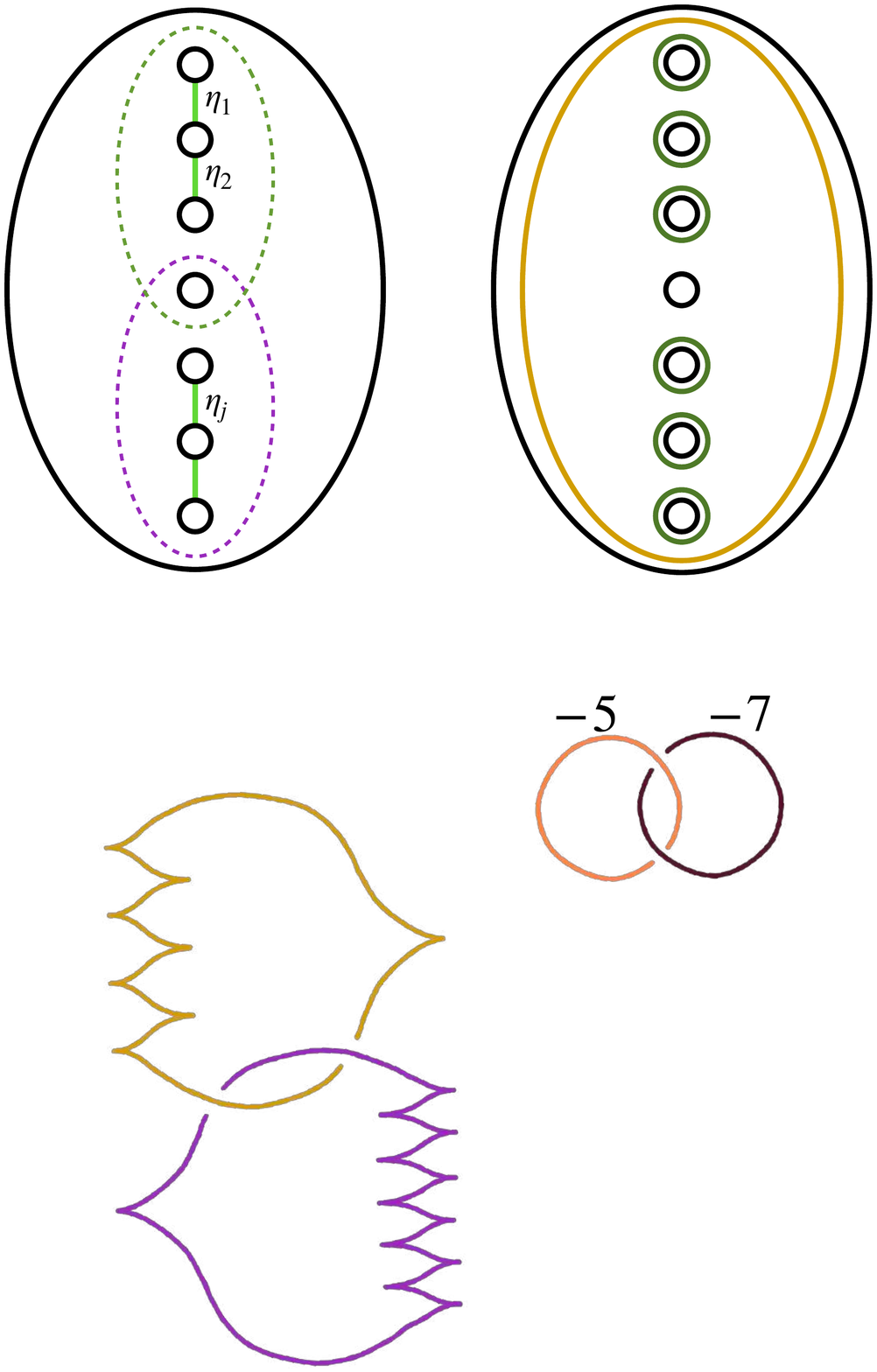}
\caption{Boundary parallel curves.}
\label{arcs2}
\end{subfigure}
\begin{subfigure}[t]{.5\textwidth}
  \centering
  \includegraphics[scale=0.4]{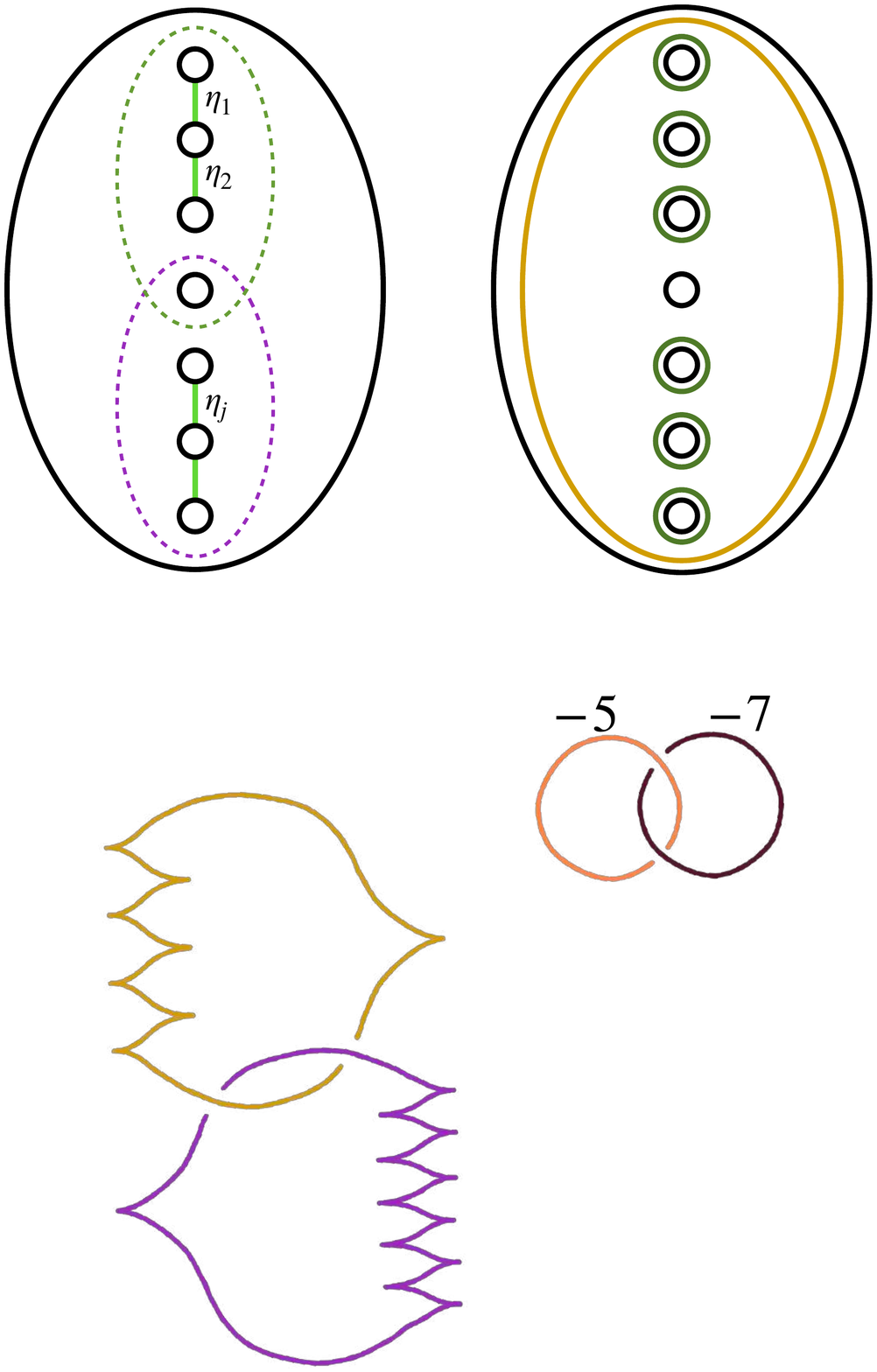}
\caption{Cutting the surface along arcs.}
\label{arcs3}
\end{subfigure}%
\begin{subfigure}[t]{.5\textwidth}
  \centering
 \includegraphics[scale=0.5]{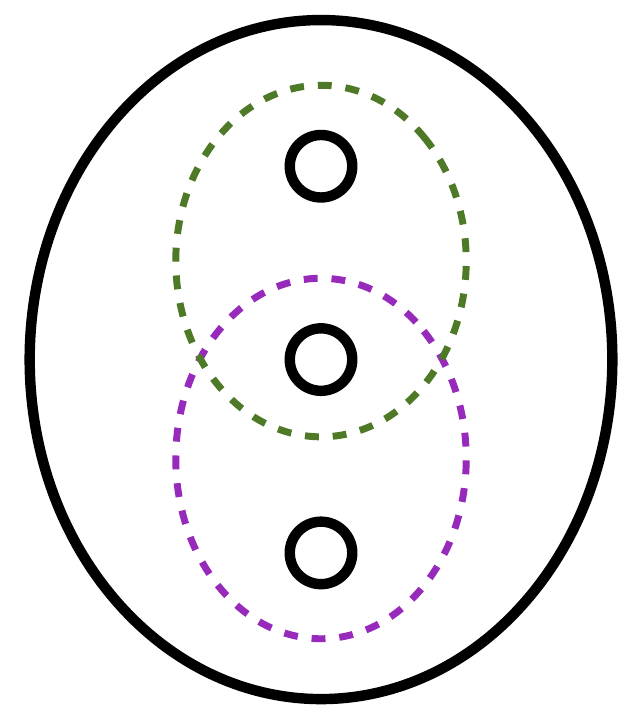}
\caption{Resulting surface and curves.}
\label{arcs4}
\end{subfigure}

\caption{Configuration of curves.}
\end{figure}

\noindent In \cite{rightveering} it is proved that any positive Dehn twist is right-veering, and that the composition of right-veering homeomorphisms is still right-veering. 
Consider the green arcs $\eta_j$ drawn in Figure \ref{arcs3}. The dashed curves $\gamma'$ and $\beta'$ are drawn like that just to indicate their homology class, while the isotopy classes are still unknown.

The arcs are disjoint from $\beta\cup\gamma$, therefore $\t_{\beta}\t_{\gamma}(\eta_j)=\eta_j$. If any of the curves $\{\beta',\gamma'\}$ (respectively purple and green) crosses one of the $\eta_j$'s, say $\gamma'$, then we would need $\t_{\beta'}$ to move the arc back to the initial position, since 
\[\eta_j=\t_{\beta}\t_{\gamma}(\eta_j)=\t_{\beta'}\t_{\gamma'}(\eta_j).\]

But this is impossible because $\t_{\beta'}$ is right-veering as well:  the tangent vector at the starting points of $\tau_{\gamma'}(\eta_j)$ is to the right of the tangent vector at $\eta_j(0)$, and when we apply $\t_{\beta'}$ we move the first vector further to the right (or we leave it where it is, depending on whether $\beta'$ intersects $\tau_{\gamma'}(\eta_j)$). This implies that $\t_{\beta'}\t_{\gamma'}(\eta_j)$ cannot be isotopic to $\eta_j$.

Therefore all the arcs $\eta_j$'s are disjoint from $\beta'\cup\gamma'$ as well. So we cut along these arcs and look for a configuration of $\{\beta',\gamma'\}$ on the resulting surface, which has still genus zero but now has just 4 boundary components (compare with Figure \ref{arcs4}). We are left with:
\[\t_{\beta}\t_{\gamma}=\t_{\beta'}\t_{\gamma'} \in \Gamma_{0,4}.\]

\subsection*{Step 4: identify the possible geometric configurations}\label{kalotilemma}

To conclude the argument for the uniqueness of the filling we need a result which is proved by Kaloti:

\begin{lem*}[\cite{kaloti}] Suppose there are two simple closed curves $\beta',\gamma'$ on $\S_{0,4}$ with $[\beta']=[\beta]\in H_1(\Sigma_{0,4};\Z)$, $[\gamma']=[\gamma]\in H_1(\Sigma_{0,4};\Z)$ and such that $\t_{\beta}\t_{\gamma}=\t_{\beta'}\t_{\gamma'} \in \Gamma_{0,4}$. 
Then there exists a diffeomorphism $\Gamma_{0,4}$ taking $\beta\mapsto \beta'$ and $\gamma\mapsto\gamma'$. 
\end{lem*}

\noindent The consequence of this lemma is that the two curves $\beta'$ and $\gamma'$ are, up to diffeomorphism, the same as $\beta$ and $\gamma$, and therefore the filling that the pair $\{\beta',\gamma'\}$ describes, together with the boundary-parallel curves, is the one described by Figure \ref{hopf}.

\subsection*{Step 5: the extra filling when $a_i=4$ for $i\in\{1,2\}$}\label{extra}

\hypertarget{thesentence}{When} $a_2$ (or $a_1$) is equal to 4, then, as expected, there is another homology configuration which is coherent with the single and joint multiplicities computed above. It is given by applying the lantern relation to the original configuration. This has the effect of reducing by one the total number of curves appearing in the factorization (hence $b_2$ of the corresponding filling is 1 and not 2). 

\begin{prop} 
If $a_2= 4$, then there are two possible homology configurations of curves respecting the single and double multiplicities computed above. 
\end{prop}

\begin{prf} 

If we go through the computation of homological configurations of curves as we did in Proposition \ref{conf}, we find this time two of them, due to the possibility of performing a lantern substitution once (notice in fact that the case when $a_1=a_2=4$ allows again just two configurations, and not three, because applying the first substitution changes the configuration preventing the second-one from being possible).

One configuration has been already described element-wise in Proposition \ref{conf} and corresponds to a Stein filling with $b_2=2$ (unique up to diffeomorphism). 

The other one is homologous (curve by curve) to the configuration we get after applying the lantern substitution to Figure \ref{page} with $l=2$ (i.e. $a_2=4$); the proof of its homology uniqueness is derived as in the proof of Proposition \ref{conf}, and it is omitted here.
\end{prf}

Then we proceed as above: all the boundary-parallel curves are placed and ignored, since, again, their homology classes determine their isotopy classes. Therefore, we can cut along appropriate arcs and reduce the number of holes appearing in the factorization, which, in turn, is the same thing as starting with a Legendrian knot whose Thurston-Bennequin number is smaller. So we can focus on the minimal possible example (after having cut along the maximal system of arcs), which has $[a_1,a_2]=[3,4]$, producing $L(11,4)$, by the computation $3-\frac{1}{4}=\frac{11}{4}$, with the contact structure specified by Figure \ref{hopf(3,4)} whose compatible open book decomposition has page as in Figure \ref{pagehopf(3,4)}. Here we immediately see that we can (uniquely) apply the lantern relation on the set of four curves given by the yellow one, the red one and the two boundary-parallel green curves (compare with Figure \ref{pagelantern}). After the substitution we get a Stein filling with $\chi=2$.

\begin{figure}[ht!]
\centering
\begin{subfigure}[t]{.45\textwidth}
  \centering
 \includegraphics[scale=0.4]{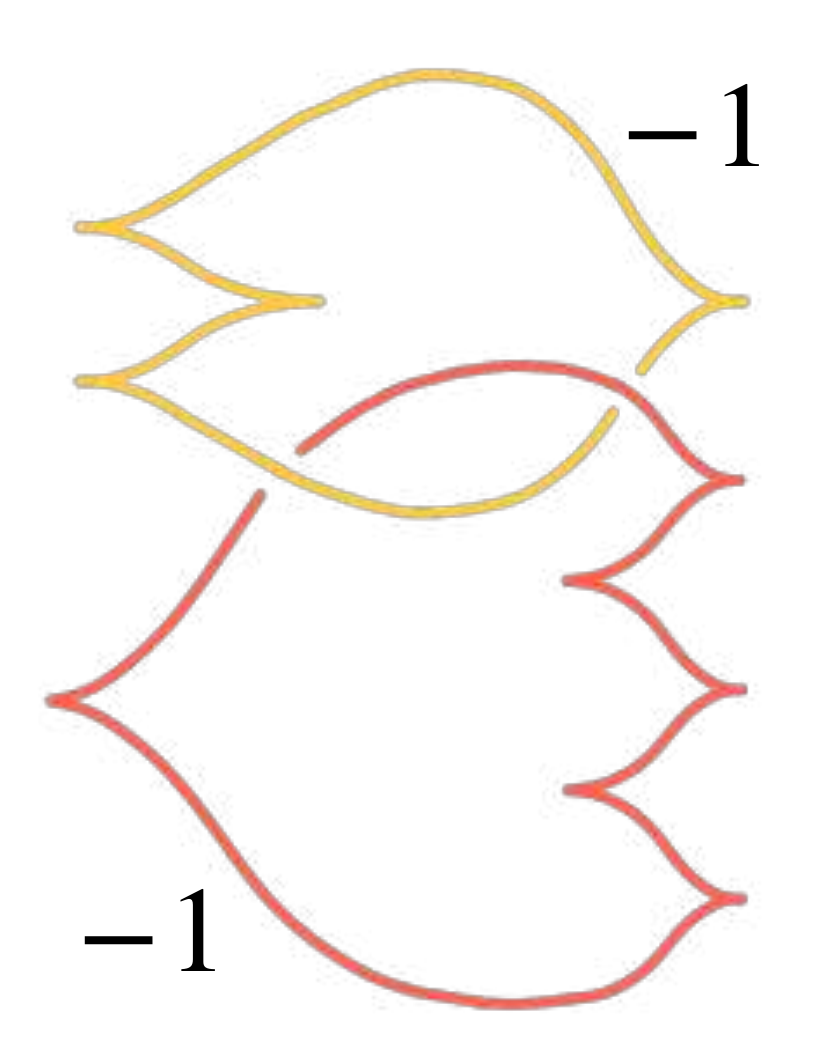}
\caption{ }
  \label{hopf(3,4)}
\end{subfigure}
\begin{subfigure}[t]{.45\textwidth}
  \centering
 \includegraphics[scale=0.3]{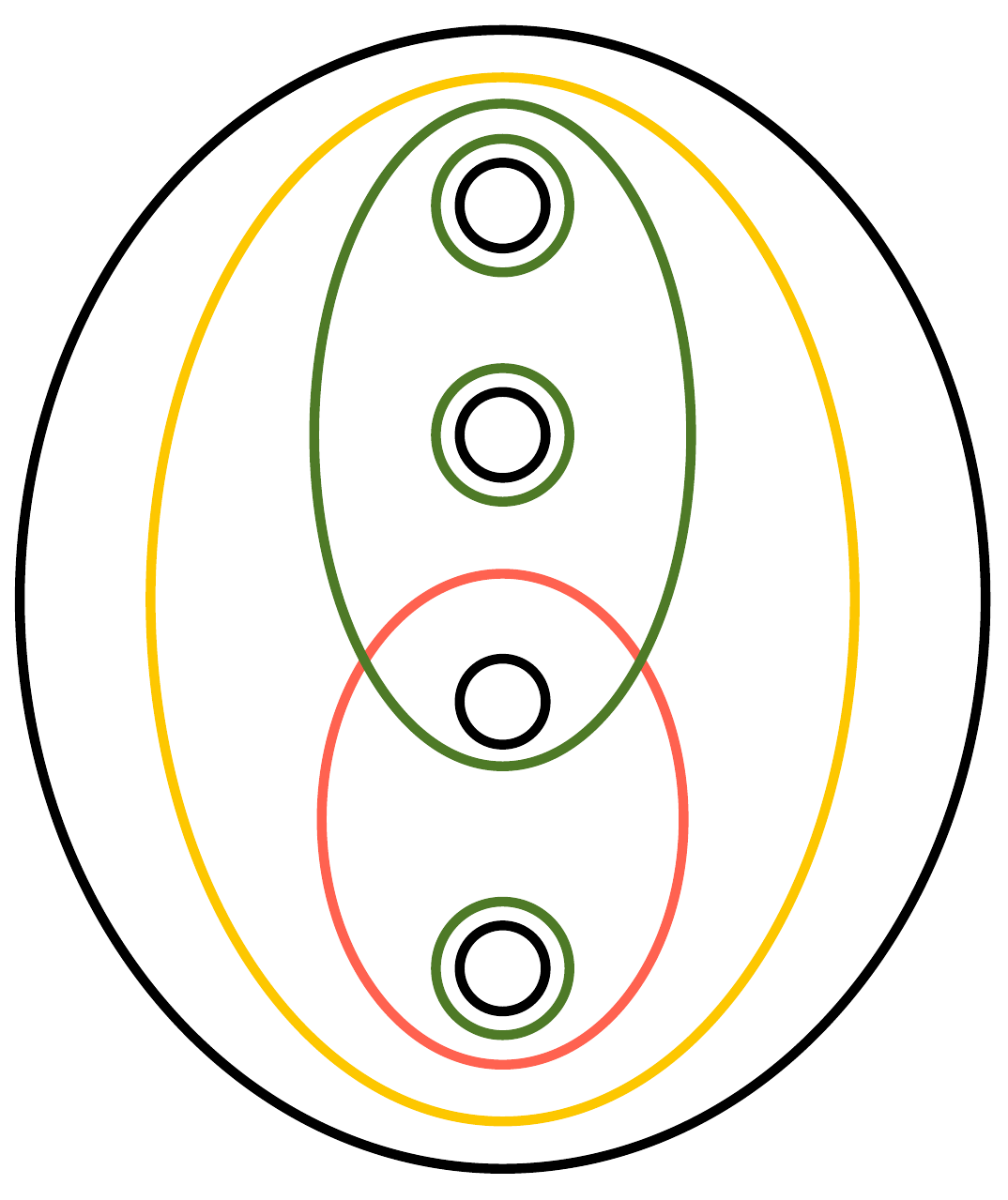}
  \caption{ }
  \label{pagehopf(3,4)}
\end{subfigure}
\caption{Tight structure on $L(11,4)$ with two fillings.}
\label{contactxi}
\end{figure}
 
\begin{prop} Let $X$ be a Stein filling of the contact 3-manifold described by Figure \ref{hopf(3,4)}, with $\chi(X)=2$. Then the homeomorphism type of $X$ is unique.
\end{prop}

\begin{prf} In order to derive our statement we need three facts:

\begin{itemize}
\item[1)] since we have determined the homology configuration of curves appearing in the factorization of the monodromy, we can use Figure \ref{pagelantern} to compute $H_1(X)$: the holes of the surface correspond to the 1-handles of $X$, and the curves themselves are the attaching circles of the 2-handles. It is immediate to see that $H_1(X)=0$, hence $\pi_1(X)$ is perfect. 
But $\pi_1(X)$ is a quotient of $\Z/11\Z\simeq \pi_1(L(11,4))$, see \cite[page 216]{ozbagci}, hence it is abelian. We conclude that $X$ is simply connected.

\item[2)] The intersection form of $X$ is characterized by the homological configuration of curves in the open book decomposition, hence it is uniquely determined (and it is isomorphic to $[-11]$).

\item[3)] The fundamental group of the boundary of $X$ is $\pi_1(L(11,4))\simeq \Z/11\Z$.
\end{itemize}

\noindent Then \cite[Proposition 0.6]{boyer} applies and tells that $X$ is unique up to homeomorphism, as claimed. The reason why we cannot describe all the (potential) diffeomorphism types is because, after applying the lantern substitution, we cannot solve explicitly the geometric configuration problem on $\Sigma_{0,5}$ with the curves we got (see Figure \ref{pagelantern}), and, moreover, there is no arc which is disjoint from these curves and on which to cut open in order to reduce the number of boundary components.
\end{prf}

\section{Final remark}

\begin{figure}[ht!]
\centering
\begin{subfigure}[t]{.45\textwidth}
  \centering
 \includegraphics[scale=0.35]{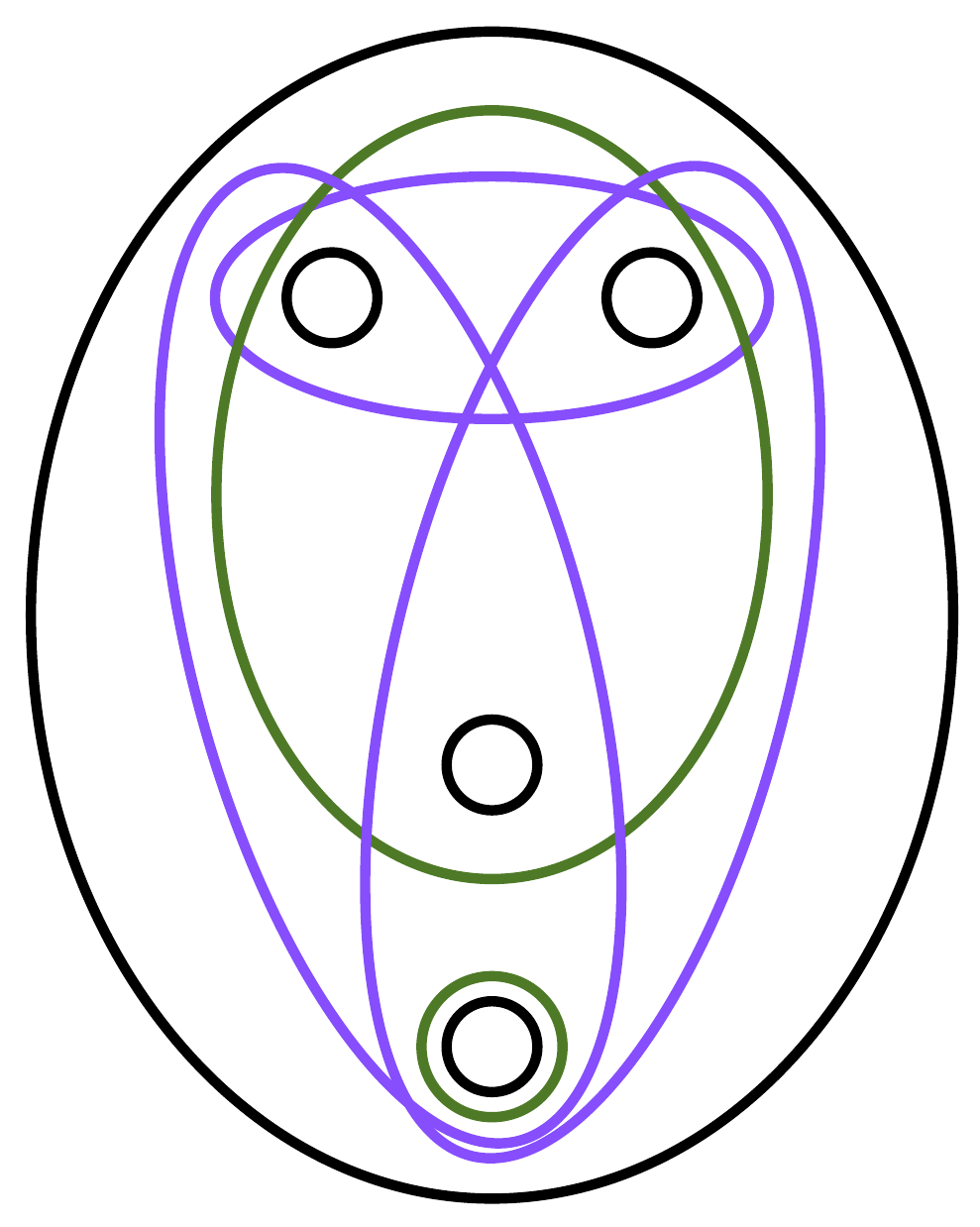}
\caption{The configuration of curves after lantern substitution.}
  \label{pagelantern}
\end{subfigure}
\begin{subfigure}[t]{.45\textwidth}
  \centering
 \includegraphics[scale=0.40]{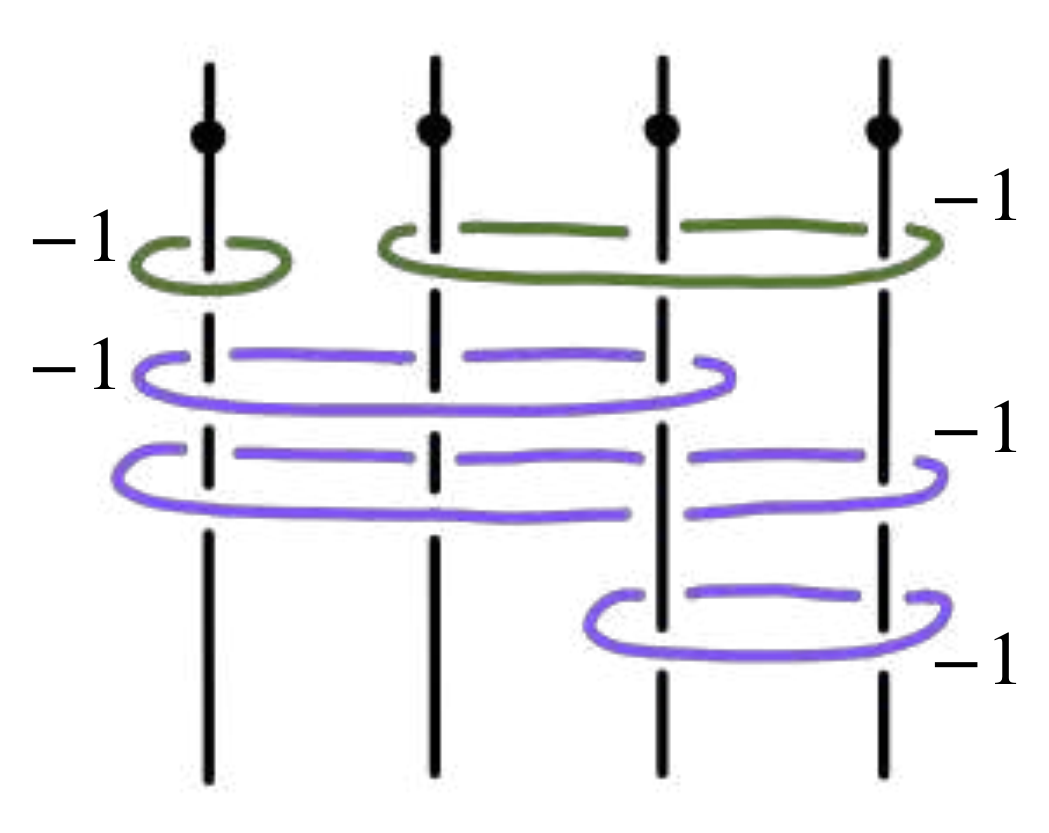}
  \caption{The corresponding 4-manifold X.}
  \label{steinlantern}
\end{subfigure}
\caption{Description of the filling with $\chi=2$.}
\end{figure}

Starting from the explicit configuration of curves corresponding to $L(11,4)$ with the contact structure of Figure \ref{hopf(3,4)}, we can apply the lanter substitution (see Figure \ref{pagelantern}) and then draw a Kirby diagram of the corresponding Stein domain $(X,J)$, see Figure \ref{steinlantern}. 

By performing handle calculus we get a new diagram of $X$ which is simpler in the following sense: this smooth handle decomposition of $X$ consists of a 0-handle and a single 2-handle, attached along the torus knot of type $(-5,2)$, pictured in Figure \ref{t(5,2)}, with framing $-11$.

\begin{figure}[ht!]
\centering
\includegraphics[scale=0.5]{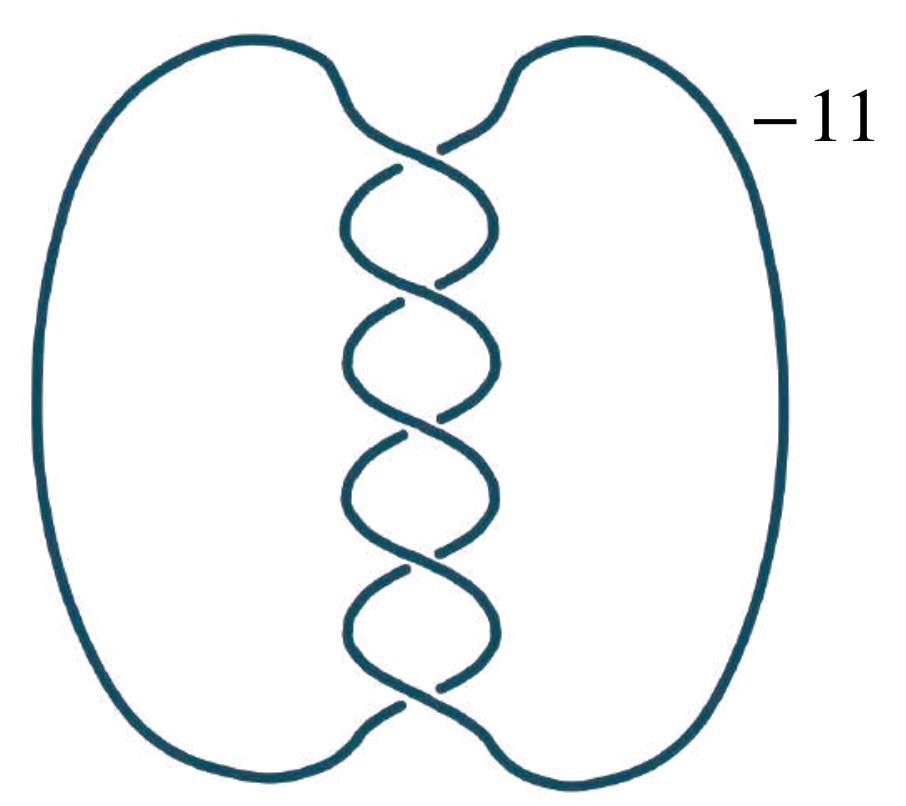}
\caption{Torus knot $T(-5,2)$.}
\label{t(5,2)}
\end{figure}

\begin{figure}[ht!]
\centering
\begin{subfigure}[t]{.45\textwidth}
  \centering
 \includegraphics[scale=0.45]{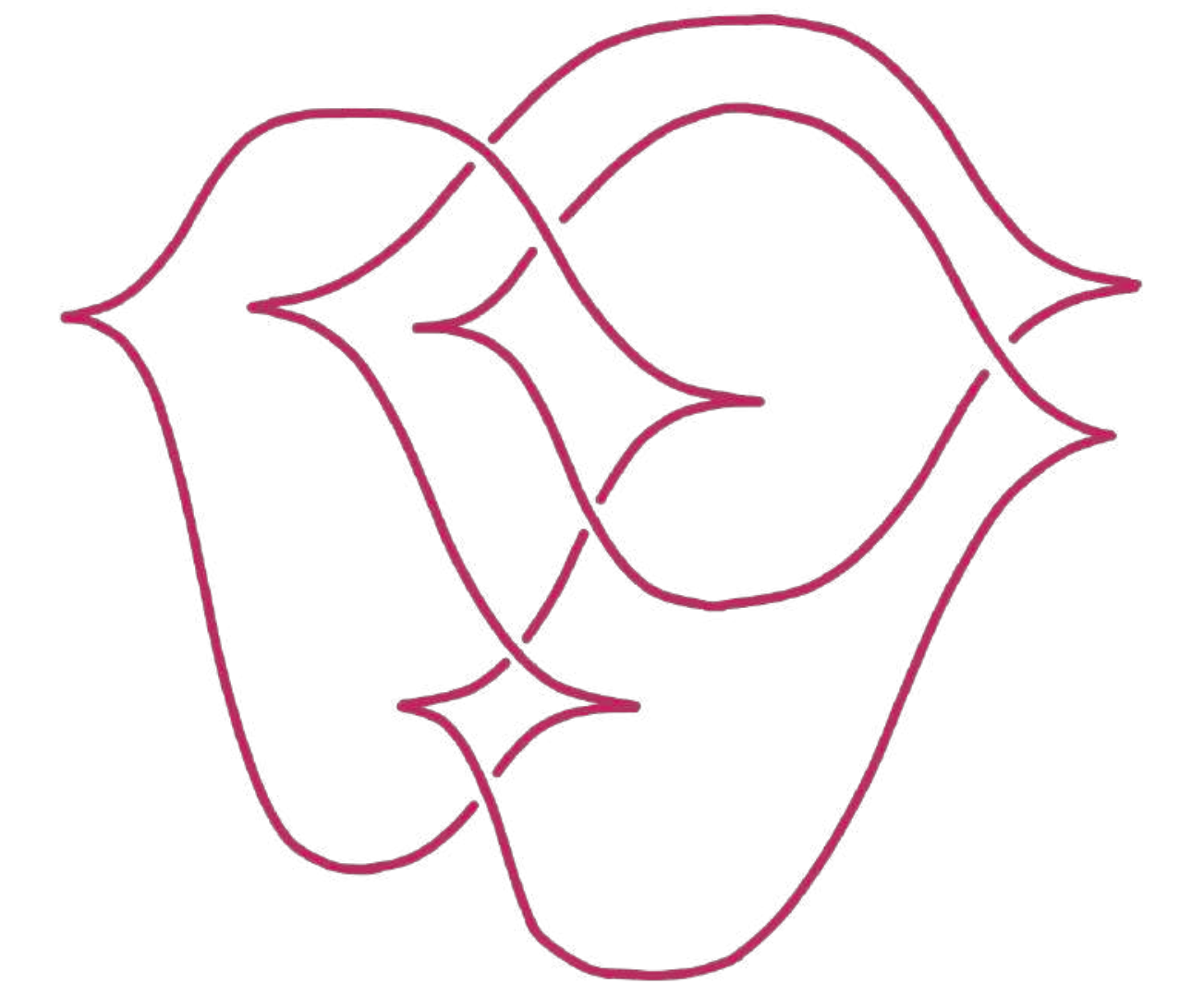}
\caption{Virtually overtwisted structure.}
  \label{t(5,2)1}
\end{subfigure}
\begin{subfigure}[t]{.45\textwidth}
  \centering
 \includegraphics[scale=0.5]{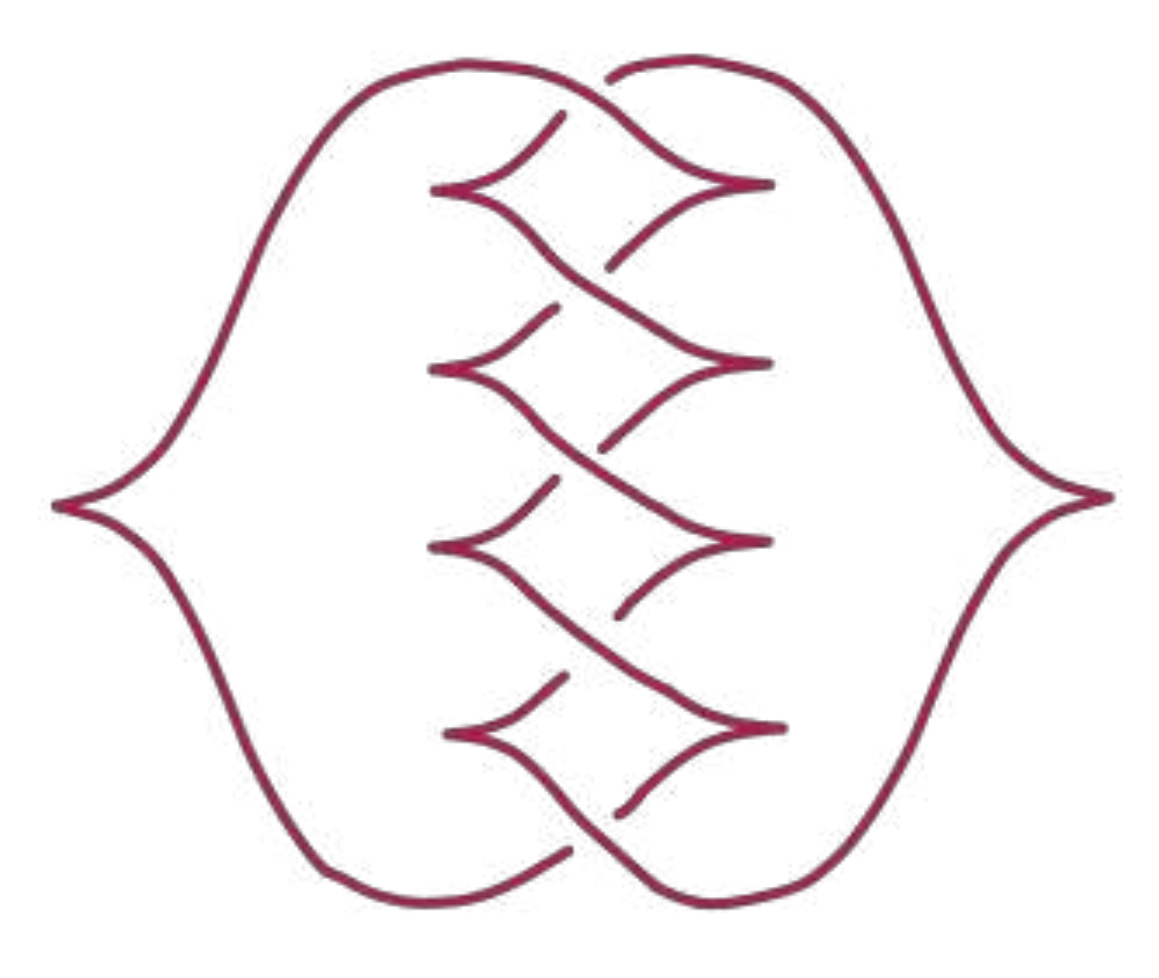}
  \caption{Universally tight structure.}
  \label{t(5,2)2}
\end{subfigure}
\caption{Different Legendrian representatives of $T(-5,2)$.}
\label{legrepr}
\end{figure}

In order to encode the Stein structure of $(X,J)$ in this handle decomposition we need a Legendrian representative of $T(-5,2)$ with Thurston-Bennequin number equal to $-10$: combining \cite[Theorem 4.3]{torusknotsimple} and \cite[Theorem 4.4]{torusknotsimple} we see that there are just two such Legendrian isotopy classes which maximize the Thurston-Bennequin number (equal to $-10$), distinguished by the rotation numbers, respectively $\pm 1$ and $\pm 3$ depending on the orientations (see Figure \ref{legrepr}). We want to understand which one suits to our case. Let $J_1$ and $J_2$ be the two Stein structures on $X$ described respectively by Figures \ref{t(5,2)1} and \ref{t(5,2)2}. 

\begin{prop} 
The Stein domain $(X,J)$ with a handle decomposition consisting of a single 2-handle is the one described by Figure \ref{t(5,2)1}, i.e. $J=J_1$. 
\end{prop}

\begin{prf} 

Call $\xi$ the contact structure described by Figure \ref{contactxi}. Remember that the two open book decompositions of Figures \ref{pagehopf(3,4)} and \ref{pagelantern} represent the same contact structure.
To prove the proposition, it is enough to check that the induced contact structure on $\partial (X,J_1)$ is isotopic to $\xi$. This is achieved by computing the 3-dimensional invariant:
\[d_3(\partial (X,J_1))=\frac{1}{4}(c_1(X,J_1)^2-3\sigma(X)-2\chi(X)).\]
If we call $K$ the Legendrian knot of Figure \ref{t(5,2)1}, then the first summand is given by 
\[\rot(K)\cdot [-11]^{-1}\cdot\rot(K),\]
while $\sigma(X)=-1$ and $\chi(X)=2$. By putting everything together we obtain
\[d_3(\partial (X,J_1))=-\frac{3}{11}.\]
On the other hand, the link of Figure \ref{hopf(3,4)} gives a Stein filling $(W,J_0)$ of $(L(11,4),\xi)$ with two 2-handles such that 
\[c_1(W,J_0)^2=[1,-2]\cdot Q_W^{-1}\cdot [1,-2]^T,  \qquad \sigma(W)=-2,	\qquad \chi(W)=3,\]
where $Q_W$ is the matrix of the intersection form, which is just the linking matrix
\[  \begin{bmatrix}
   -3 & 1 \\
   1 & -4 \\
   \end{bmatrix}.\] 
The computation shows again that
\[d_3(\partial (W,J_0))=-\frac{3}{11}.\]
Moreover, in the case when the rotation number of the second component of the Legendrian link is 0, the $d_3$ invariant of the resulting contact structure is $-1/11$. According to Honda's classification of tight contact structures on $L(11,4)$, this computation covers all the three (up to contactomorphism) possible cases, see next paragraph for computation in the universally tight case.

Therefore, we conclude that
\[(L(11,4),\xi)=\partial (W,J_0)=\partial(X,J_1).\]
Hence $J=J_1$, as wanted.
\end{prf}

\begin{figure}[ht!]
\centering
\includegraphics[scale=0.4]{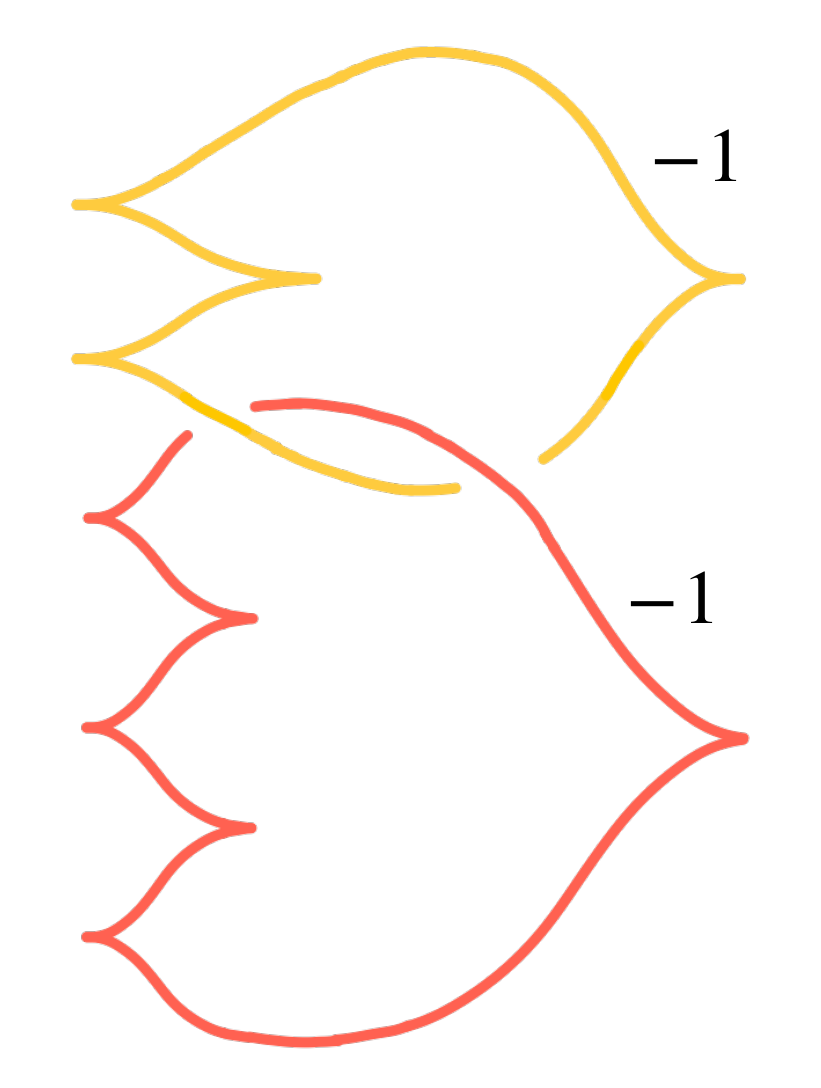}
\caption{Universally tight contact structure on $L(11,4)$.}
\label{hopf34ut}
\end{figure}

To conclude, we check that the other Legendrian representative of the torus knot $T(-5,2)$ with Thurston-Bennequin number $-10$ gives a different contact structure: by performing contact $(-1)$-surgery on the Legendrian knot of Figure \ref{t(5,2)2}, we get $L(11,4)$ with a universally tight structure. This is proved by comparing its $d_3$ invariant with the one computed from Figure \ref{hopf34ut}.

\noindent In both cases we get 
\[d_3=-\frac{5}{11}. \]
Therefore the two different Stein fillings of $(L(11,4),\xi_{ut})$ are described by the handle diagrams of Figures \ref{hopf34ut} and \ref{t(5,2)2}.

\subsection{Further generalization problems}
One can imagine of following the steps of Section \ref{classificationthm} to produce and classify the Stein fillings of those virtually overtwisted contact structure on lens spaces obtained from Legendrian surgery on a 3-components chain of unknots, or even on a longer one. Extending Proposition \ref{conf} is just a matter of carefully studying the combinatorics of the multiplicities numbers, but no substantial difficulty should arise here, at least in the case of $\length(p/q)=3$. The critical point of the proof that we presented is Kaloti's lemma of Section \ref{kalotilemma}, which has no known analogous for surfaces with more boundary components. If a result which identifies a unique configuration of curves in a base case were available, one might try to reproduce the steps in the proof of the classification theorem and extend Theorem \ref{theorem1} to lens spaces $L(p,q)$ with $\length(p/q)=3$.

\chapter{Topological constraints for Stein fillings}\label{topologicalconstraints}

As discussed in Section \ref{steinfillings}, classifying symplectic fillings (up to homeomorphism, diffeomorphism or symplectic deformation equivalence) of a given contact 3-manifold can be a very hard task, even though some progress has been made in the last years. A more modest approach is trying to give some constraints on the topological invariants of the Stein fillings, even if a complete classification is missing. If we restrict to planar contact structures, then studying Stein fillings is enough if we want to understand weak symplectic fillings, since these are symplectically deformation equivalent to blow ups of Stein fillings, see \cite[Theorem 2]{niederkruger}.
Some topological constraints for Stein fillings of planar contact structures have already been found (see for example \cite{etnyre}, \cite{c1planar}, \cite{wendl}, \cite{wand}), and here we specifically focus on lens spaces $(L(p,q),\xi)$.

Throughout this chapter we often refer to the length $l=\length(p/q)$ of the expansion $\frac{p}{q}=[a_1,a_2,\ldots, a_l]$. To this expansion we can associate a negative linear graph $\Lambda(p,q)$ and a corresponding negative definite 4-manifold $X_{\Lambda(p,q)}$ realized as a plumbing. We give a sharp upper bound on the possible values of the Euler characteristic for a minimal symplectic filling of a tight contact structure on a lens space:

\begin{thm} \label{maxchi}
Let $\xi$ be any tight contact structure on $L(p,q)$. Let $W$ be a minimal symplectic filling of $L(p,q)$ and let $l=\length(p/q)$. Then
\[\chi(W)\leq 1+l.\]
\end{thm}

This estimate is obtained by looking at the topology of the spaces involved, extending this way what we already knew from the universally tight case to the virtually overtwisted one. Moreover, the upper bound is always realized by a minimal symplectic filling ($X_{\Lambda(p,q)}$ itself supports a Stein structure inducing the prescribed contact structure on its boundary) whose intersection form and fundamental group are uniquely determined:

\begin{thm}\label{uniqueintform}
Let $\xi$ be any tight contact structure on $L(p,q)$ and let $l=\length(p/q)$. Let $X$ be a minimal symplectic filling of $(L(p,q),\xi)$ with $b_2(X)=l$. Then the intersection form $Q_X$ is isomorphic to the intersection form of $X_{\Lambda(p,q)}$. Moreover, $X$ is simply connected.
\end{thm}

We also prove the following corollary, regarding the uniqueness (in certain cases) of the filling with maximal Euler characteristic:

\begin{cor}\label{homeo}
Let $\xi$ be a tight contact structure on $L(p,q)$ and let $l=\length(p/q)$. Let $X$ be a minimal symplectic filling of $(L(p,q),\xi)$ with $b_2(X)=l$. Assume that $p\in\{2,4,s^n,2s^n\}$, for some odd prime $s$ and positive integer $n$. Then $X$ is homeomorphic to $X_{\Lambda(p,q)}$.
\end{cor}

On the other hand, the first and third Betti numbers of a Stein filling $W$ of a lens space $L(p,q)$ are always zero \cite[page 216]{ozbagci}, hence we have an obvious lower bound on the value $\chi(W)$, which is $\chi(W)=b_0(W)=1$: this is realized precisely when $(L(p,q),\xi)$ bounds a Stein rational homology ball. Lisca proved in \cite{liscaribbon} that, in order to guarantee the existence of rational balls with boundary a lens space, the numbers $p$ and $q$ must fall into one of three families with specific numerical conditions.

Among those, we restrict to the case when $p$ and $q$ are of the form $p=m^2$ and $q=mk-1$, for some $m>k>0$ with $(m,k)=1$. It is known \cite[Corollay 1.2c]{lisca} that $L(m^2,mk-1)$ endowed with a universally tight contact structure bounds a Stein rational ball and we use this fact to prove that in the virtually overtwisted case this never happens, concluding that:

\begin{thm} \label{minchi}
Let $W$ be a symplectic filling of $(L(p,q),\xi_{vo})$, with $p=m^2$ and $q=mk-1$, for some $m>k>0$ and $(m,k)=1$. Then $\chi(W)\geq 2$.
\end{thm}

Theorem \ref{minchi} can be generalized to the other families of lens spaces which are known to bound a smooth rational homology ball: these balls do not support any symplectic structure, i.e. none of the virtually overtwisted contact structures can be filled by a Stein rational ball, see \cite[Proposition A.1]{gollastar}.

Then we turn our attention to covering maps: since an overtwisted disk lifts to an overtwisted disk, all the coverings of a universally tight structure are themselves tight. The situation is less clear when we consider virtually overtwisted structures. By starting with such a structure on a lens space, we know that it lifts to an overtwisted structure on $S^3$, but what happens to all the other intermediate coverings? 

One of the problems we faced when studying contact structures along covering maps is the mysterious behavior of numbers: for example, there is no understanding on how the lengths of $p/q$ and $p'/q$ are related, if $p'$ is a divisor of $p$. This makes the problem hard even to organize, since we could not glimpse any clear scheme or pattern for stating reasonable guesses. Theorem \ref{thmcovering} is the only stance of a general result which does not depend on specific examples.

\begin{thm}\label{thmcovering}
Let $p,\,q$ and $d$ be such that $q<p<dq$. Then every virtually overtwisted contact structure on $L(p,q)$ lifts along a degree $d$ covering to a structure which is overtwisted.
\end{thm}

In Section \ref{coveringsection} we give a series of examples using the description of tight structures given by Honda in \cite{honda} to study explicit cases of covering maps between contact lens spaces. The last part of this chapter is dedicated to the study of the fundamental group of Stein fillings of virtually overtwisted structures on lens spaces, combining the results above about Euler characteristic with what we developed on the behavior of coverings. Recall that the fundamental group of any Stein filling is a quotient of the fundamental group of its boundary, see Remark \ref{b1zero}. 
As a consequence of Theorem \ref{maxchi}, we will prove Theorem \ref{pi1b2} and provide some specific examples and applications.

\begin{thm}\label{pi1b2} 
Let $X$ be a Stein filling of $(L(p,q),\xi)$ with $\pi_1(X)=\Z/d\Z$, for $p=dp'$. Then
\[\chi(X)\leq \frac{1+l'}{d},\]
where $l'=\length(p'/q')$, with $q'\equiv q \pmod{p'}$.
\end{thm}

\section{Upper bound for the Euler characteristic}

The goal of this section is to prove Theorems \ref{maxchi} and \ref{uniqueintform}. First, recall that a vertex $v$ of a weighted graph is a \emph{bad vertex} if 
\[w(v)+d(v)>0,\]
where $w(v)$ and $d(v)$ are respectively the weight and the degree (i.e. the number of edges containing $v$) of the vertex.

Graphs with no (or at most one) bad vertex are studied in the context of Heegaard Floer homology, links of singularities and planar contact structures, in several works including for example \cite{badvertex}, \cite{nemethifive}, \cite{nemethilinks} and \cite{ghigginigolla}.

We show that:

\begin{thm}\label{maxtree}
Let $\Gamma$ be a negative definite plumbing tree with $k$ vertices, none of which is a bad vertex. Call $\overline{Y}$ the plumbed 3-manifold associated to $\Gamma$ and assume that $\overline{Y}$ is a rational homology sphere. Denote by $Y$ the manifold with the opposite orientation. Let $X$ be a negative definite smooth 4-manifold with no $(-1)$-class in $H_2(X;\Z)$ such that $\partial X=Y$. Then
\[b_2(X)\leq 1+\sum_{i=1}^k (|w(v_i)|-2).\]
\end{thm}

\begin{prf} 
Let $P=P_{\Gamma}$ the plumbed 4-manifold associated to $\Gamma$, whose oriented boundary is $\overline{Y}$, and whose intersection form is $Q_P=Q_{\Gamma}$. Form the closed manifold $W=X\cup_{\partial}P$ by gluing the two manifolds along the boundary, see Figure \ref{glue1}. 
\begin{figure}[ht!]
	\center
	\includegraphics[scale=0.5]{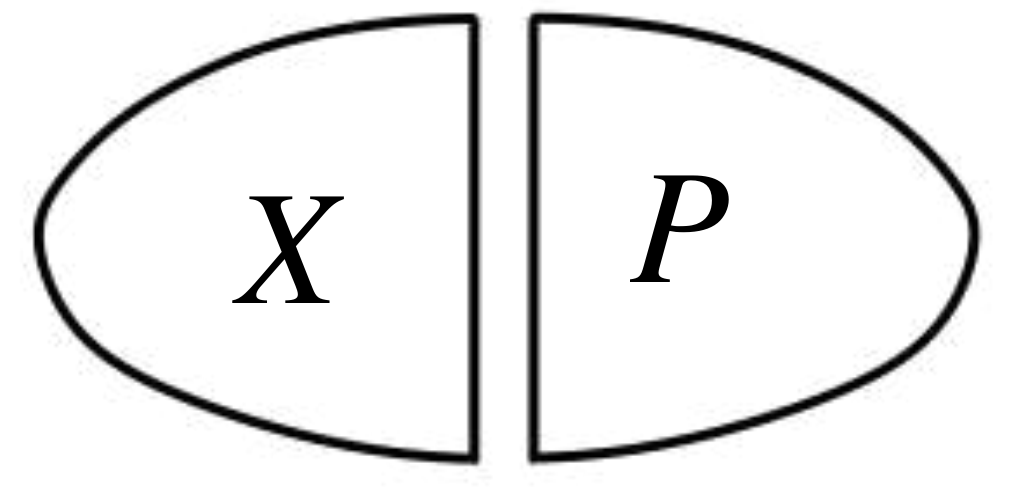}
		\caption{The closed manifold $W=X\cup_{\partial}P$.}
	\label{glue1}
\end{figure}

\noindent We get a closed smooth 4-manifold $W$ whose intersection form $Q_W$ is negative definite, and hence, by Donaldson's theorem \cite{donaldson}, isomorphic to $\langle -1\rangle^r$, for some $r$. Since $Y$ is a rational homology sphere, we have that
\[\rank(Q_X)+\rank(Q_{\Gamma})=b_2(X)+b_2(P)=b_2(X)+k=r.\]
A priori, $Q_X$ is a sub-lattice of finite index $n$, for some $n>0$:
\[Q_X\leq (Q_{\Gamma})^{\perp},\]
but the following lemma shows that they coincide. 

\begin{lem} In the setting of above we have an isomorphism $Q_X\simeq (Q_{\Gamma})^{\perp}$.
\end{lem}

\begin{prf}
Look at the exact sequence of the pair $(W,X)$:
\[ H_3(W,X)\to H_2(X) \to H_2(W) \to H_2(W,X) \]
and notice that, by excision and Poincaré-Lefschetz duality
\[H_3(W,X)\simeq H_3(P,\partial P)\simeq H^1(P)=0\] 
and similarly
\[H_2(W,X)\simeq H_2(P,\partial P)\simeq H^2(P).\]
The latter is free, because its torsion comes from $H_1(P)$, which is $0$. Therefore, the inclusion $H_2(X)\hookrightarrow H_2(W)$ has a free quotient, being this a subgroup of the free group $H_2(W,X)$. Hence, if we take a class $\alpha\in H_2(W)$ with the property that $n\alpha$ is inside $H_2(X)$, we automatically get $\alpha\in H_2(X)$. In particular, if $\alpha\in (Q_{\Gamma})^{\perp}$, then $n\alpha\in (H_2(X),Q_X)$, with $n$ equal to the index $Q_X\leq (Q_{\Gamma})^{\perp}$, and hence $\alpha\in (H_2(X),Q_X)$. So $Q_X\simeq (Q_{\Gamma})^{\perp}$
\end{prf}

\noindent The isomorphism 
\[Q_X\simeq (Q_{\Gamma})^{\perp}\]
implies that we have an embedding $\Gamma\hookrightarrow \langle -1\rangle^r$ with the property that there is no $(-1)$-class in the orthogonal, otherwise this would come from $X$, which, by assumption, does not have any. We call such an embedding \emph{irreducible}. 
Notice that, up to isomorphism, there is a unique maximal irreducible embedding
\[\Gamma\hookrightarrow \langle -1\rangle^t,\]
where by \emph{maximal} we mean that the dimension $t$ (which, by the argument below, is finite) of the ambient lattice cannot be bigger. First of all notice that at least one irreducible embedding exists: we will explicitly describe the construction of one of them, which turns out to be the maximal one. Since the sum of the weights of the graph is finite, $t$ is itself finite.
To embed the graph $\Gamma$ in such a way that there is no $(-1)$-class in the orthogonal complement implies that all the elements in the canonical basis $\{e_1,\ldots,e_t\}$ of $\langle -1\rangle^t$ appear in the image of some vertex of the graph. Therefore, to obtain the maximal such $t$, we have to impose only the requirements that:
\begin{itemize}
\item[1)] the $i^{th}$ vertex is sent to a combination of $|w_i|$ distinct basis elements and
\item[2)] any two adjacent vertices of $\Gamma$ share, via the embedding, exactly one element $e_j$.
\end{itemize}
If one of these conditions is not satisfied, then we end up with (at least) one line $\langle e_j \rangle$ which is not hit by the image of $\Gamma$ and that will produce an element in the orthogonal with square $-1$. So the image of the first vertex with weight $w_1$ must be a sum of $|w_1|$ distinct elements $e_i$. The second vertex is sent to a combination of $|w_2|$ elements, among which exactly one has already appeared in the image of the first vertex, and so on. 

The fact that the there are no bad vertices guarantees that it is possible to go on with this recipe and send every vertex with weight $w_i$ into a combination of $|w_i|$ distinct basis elements, where the repetitions between the images of different vertices occur exactly in correspondence of the edges. This provides a way to construct an irreducible embedding $\Gamma\hookrightarrow \langle -1\rangle^t$, which is unique up to isomorphism. Hence we find
\[t=1+\sum_{i=1}^k (|w(v_i)|-1).\]
Therefore, since the dimension of the maximal irreducible embedding of $\Gamma$ is as above, we have $r=k+b_2(X)\leq t$. We conclude:
\[b_2(X)\leq  t-k=1+\sum_{i=1}^k (|w(v_i)|-1)-k=1+\sum_{i=1}^k (|w(v_i)|-2).\]
\end{prf}

\begin{cor}\label{uniquenessgeneral}
In the setting above, the intersection form of the the manifolds with boundary $Y$ and maximal $b_2$ is uniquely determined, up to isomorphism.
\end{cor}

\begin{prf} Assume that $X_1$, $X_2$ are negative definite with no $(-1)$-class, $\partial X_1=\partial X_2=Y$ and with $b_2(X_i)$ maximal. Then, by uniqueness of the maximal irreducible embedding $\Gamma\hookrightarrow \langle -1\rangle^t$, we have that 
\[Q_{X_1}\simeq Q_{X_2}\simeq \Gamma^{\perp}\subseteq \langle -1\rangle^t.\]
\end{prf}

Now we specialize to the case of lens spaces. Start with $L(p,q)$ and take the expansion $p/q=[a_1,a_2,\ldots, a_l]$, where all the $a_i$'s are $\geq 2$. Call $\Lambda=\Lambda(p,q)$ the associated negative definite lattice with $l$ vertices (where $l=\length(p/q)$):
\begin{center}
$\Lambda=$
\begin{tikzpicture}
        \node[shape=circle,fill=black,inner sep=1.5pt,label=$-a_1$] (1)                  {};
        \node[shape=circle,fill=black,inner sep=1.5pt,label=$-a_2$] (2) [right=of 1] {}
        edge [-]               (1);
        \node[shape=circle,fill=black,inner sep=1.5pt,label=$-a_3$] (3) [right=of 2] {}
        edge [-]               (2);
        \node[shape=circle,fill=black,inner sep=1.5pt,label=$-a_{l-1}$] (4) [right=of 3] {} ;
                 \node at ($(3)!.5!(4)$) {\ldots};

        \node[shape=circle,fill=black,inner sep=1.5pt,label=$-a_l$] (5) [right=of 4] {}
        edge [-]               (4);
\end{tikzpicture}
\end{center}
We apply Riemenschneider's dots method \cite{riemen} to build a negative definite 4-manifold with boundary $L(p,p-q)$ whose intersection lattice will be called $\Lambda^{\nu}$. This is obtained by reading column-wise the entries of Table \ref{table}.
\begin{center}
\begin{table}[h!]
        \centering
        \begin{tabular}{ cccccccccccc}
            \multicolumn{4}{l}{$\overbrace{\rule{2.2cm}{0pt}}^{a_1-1}$}  &  & & & & & & & \\
           $\bullet$  & $\bullet$ &  $\cdots$ &  $\bullet$  &   &  &   &   &   &  &   &    \\
             \hline 
           & & &$\bullet$  & $\bullet$ &  $\cdots$ &  $\bullet$  &   &  &   &   &   \\
             \hline 
           & & &  & &   &    & $\ddots$   &  &   &   &   \\
             \hline 
           & & &  & &   &    &   &   $\bullet$  & $\bullet$ &  $\cdots$ &  $\bullet$\\
           & & &  \multicolumn{4}{l}{$\underbrace{\rule{2.2cm}{0pt}}_{a_2-1}$} & &  \multicolumn{4}{l}{$\underbrace{\rule{2.2cm}{0pt}}_{a_l-1}$} 
 \end{tabular}
  \caption{Riemenschneider's dots method.}
  \label{table}
  \end{table}
\end{center}

\noindent If we call $l^{\nu}$ the number of columns and set 
\[c_j=1+\#\{\mbox{dots in the $j^{\mbox{\tiny{th}}}$ column}\},\]
then we obtain the continued fraction expansion of $p/(p-q)$ as
\[\frac{p}{p-q}=[c_1,c_2,\ldots,c_{l^{\nu}}].\]
Before proving Theorem \ref{maxchi}, we need a lemma.

\begin{lem}\label{lem1}
\[\length(p/q)+\length(p/(p-q))=1+\sum_{i=1}^l(a_i-1).\]
\end{lem}

\begin{prf}
We know that $\length(p/q)=l$, so we compute $\length(p/(p-q))$. This is just the number $l^{\nu}$ of columns:
\begin{align*}
l^{\nu} = & (a_1-1)+(a_2-2)+\ldots+(a_l-2)\\
=\, &\sum_{i=1}^l(a_i-2)+1\\
=\, &\sum_{i=1}^l(a_i-1)-l+1.
\end{align*}
Therefore $\length(p/q)+\length(p/(p-q))=1+\sum_{i=1}^l(a_i-1)$. 
\end{prf}

\noindent By switching the roles of $p/q$ and $p/(p-q)$, it is clear from Lemma \ref{lem1} that
\begin{align*}
\length(p/q)+\length(p/(p-q))=& \rank(\Lambda)+\rank(\Lambda^{\nu})\\
=\,& l+l^{\nu}\\
=\,& 1+\sum_{i=1}^l(a_i-1)\\
=\, & 1+\sum_{i=1}^{l^{\nu}}(c_i-1).
\end{align*}

\begin{nrem}\label{b1zero}
In the book \cite[Section 12.3]{ozbagci} the authors made the following observation, which we will often use in this work. If $X$ is a Stein filling of $Y$, then the morphism $\pi_1(Y)\to \pi_1(X)$, induced by the inclusion, is surjective since $X$ can be built on $Y \times [0, 1]$ by attaching 2-, 3- and 4-handles only. In particular, $b_1(X)\leq b_1(Y)$ and if $Y$ is a lens space, then $b_1(X)=0$.
\end{nrem}

Theorem \ref{maxchi} follows now from Theorem \ref{maxtree}:

\begin{prf}[of Theorem \ref{maxchi}]
Let $Y=L(p,q)$, so that $\overline{Y}$ is the 3-manifold associated to $\Lambda^{\nu}$, with $\Lambda^{\nu}$ playing the role of $\Gamma$. The setting for lens spaces is coherent with the hypotheses of Theorem \ref{maxtree}:
\begin{itemize}
\item lens spaces arise as plumbings on trees with no bad vertices;
\item lens spaces are rational homology spheres;
\item contact structures on lens spaces are planar (\cite[Theorem 3.3]{schonenberger}), and therefore $b_2(X)=b_2^-(X)$ for any minimal filling $X$ (\cite{etnyre});
\item minimal fillings of planar contact structures have no $(-1)$-class, as proved in \cite[Corollary 1.8]{ghigginigolla};
\end{itemize}
Therefore, since any minimal filling $X$ of $(Y,\xi)$ has $b_1=0$, as seen in Remark \ref{b1zero}, we have:
\begin{align*}
\chi(X)= & 1+b_2(X)\\
\leq\, & 1+\left(1+\sum_{i=1}^{l^{\nu}} (c_i-2)\right)\\
=\, & 1+\left(1+\sum_{i=1}^{l^{\nu}} (c_i-1)-l^{\nu}\right)\\
=\, & 1+(l+l^{\nu}-l^{\nu})\\
=\, & 1+l.
\end{align*}
\end{prf}

\begin{prf}[of Theorem \ref{uniqueintform}]
The fact that the intersection form is uniquely determined is just a special case of Corollary \ref{uniquenessgeneral}. For the fundamental group, let $X$ be a filling with $b_2(X)=l=\length(p/q)$. We know that
\[Q_X\simeq Q_{X_{\Lambda(p,q)}},\]
and we look at the long exact sequence of the pair $(X,\partial X)$, with $\partial X=L(p,q)$:
\[\xymatrix{ H_2(L(p,q))\ar[r] & H_2(X)\ar[r]^-{Q_X} &H_2(X,L(p,q))\ar[r]& H_1(L(p,q))\ar[r]& H_1(X)\ar[r] &H_1(X,L(p,q))}\]
Notice that:
\begin{itemize}
\item[1)] $H_2(L(p,q))\simeq H^1(L(p,q))=0$;
\item[2)] $H_2(X)\simeq \Z^l$;
\item[3)] $H_2(X,L(p,q))\simeq H^2(X)\simeq \Z^l\oplus H_1(X)$;
\item[4)] $H_1(L(p,q))\simeq \Z/p\Z$;
\item[5)] $H_1(X,L(p,q))\simeq H^3(X)=0$;
\item[6)] $\det(Q_X)=p$.
\end{itemize}
Therefore, by substituting everything, it follows that $H_1(X)=0.$ But since, by Remark \ref{b1zero}, $\pi_1(X)$ is abelian, we have that $\pi_1(X)=0$, as wanted.
\end{prf}

\noindent We can now give a proof of Corollary \ref{homeo}.

\begin{prf}[of Corollary \ref{homeo}] To prove this corollary we need three facts:

\begin{itemize}
\item[1)] $X$ and $X_{\Lambda(p,q)}$ are both simply connected by Theorem \ref{uniqueintform};
\item[2)] $X$ and $X_{\Lambda(p,q)}$ have isomorphic intersection forms by Theorem \ref{uniqueintform};
\item[3)] the fundamental group of their boundary is $\pi_1(Y)\simeq \Z/p\Z$, with $p\in\{2,4,s^n,2s^n\}$, for some odd prime $s$ and positive integer $n$.
\end{itemize}

\noindent Then \cite[Proposition 0.6]{boyer} applies and tells that $X$ and $X_{\Lambda(p,q)}$ are homeomorphic.
\end{prf}

\section{Lower bound for the Euler characteristic} \label{lowerboundsection}


The goal of this section is to prove Theorem \ref{minchi}, i.e. that, among the virtually overtwisted structures on the lens spaces of the form $L(m^2,mk-1)$ with $(m,k)=1$, none of these can be filled by a Stein rational homology ball.  

The first thing to notice is that, thanks to Honda's classification result \cite{honda}, each tight contact structure on a lens space has a Legendrian surgery presentation which comes from placing the corresponding chain of unknots into Legendrian position with respect to the standard contact structure of $S^3$. So, by varying the rotation numbers of the various components of the link, we can describe all the tight contact structures that a lens space supports, up to isotopy.

Let $(Y,\xi)$ be a contact 3-manifold with $c_1(\xi)$ a torsion class. Then, in \cite{gompfsymp}, Gompf defined the invariant
\[d_3(Y,\xi)=\frac{1}{4}(c_1(X,J)^2-3\sigma(X)-2\chi(X))\in \Q,\]
where $(X,J)$ is any almost complex 4-manifold with boundary $\partial X=Y$ such that $\xi$ is homotopic to $TY\cap JTY$ (compare with Lemma 6.2.6 of \cite{ozbagci}).

\begin{lem} If $(Y,\xi)$ bounds a Stein rational homology ball, then $d_3(Y,\xi)=-\frac{1}{2}$.
\end{lem}

\begin{prf} The quantity $d_3=\frac{1}{4}(c_1^2-3\sigma-2\chi)$ does not depend on the chosen filling, and if $(Y,\xi)=\partial (X,J)$ with $H_2(X;\Q)=H_1(X;\Q)=0$, then
\[d_3(Y,\xi)=\frac{1}{4}(c_1(X,J)^2-3\sigma(X)-2\chi(X))=\frac{1}{4}(0-0-2)=-\frac{1}{2}.\]
\end{prf}

\noindent In the case of lens spaces, the computation of the $d_3$ invariant is as follows:
\[d_3=\frac{1}{4}(c_1^2-3	\sigma-2(1-\sigma))=\frac{1}{4}(c_1^2-\sigma-2),\]
because all the Stein fillings have $b_1=b_3=0$ (see Remark \ref{b1zero}) and $b_2=b_2^-$ by \cite{etnyre}.
This means that, if $(L(p,q),\xi)$ bounds a Stein rational ball, then for \emph{any} other filling $(X,J)$ we have:
\[-\frac{1}{2}=\frac{1}{4}(c_1(J)^2-\sigma(X)-2)\]
and hence
\begin{equation}\label{c2=sigma}
c_1(J)^2=\sigma(X).
\end{equation}

We want to compute $c_1(J)^2$ for the filling $(X,J)$ of $(L(p,q),\xi)$ realized as the plumbing described by the linear graph of the expansion $p/q$.
To do this, we need to specify the vector $r$ of rotation numbers for the components of the linear plumbing. If
\[\frac{p}{q}=[v_1,v_2,\ldots v_n]=v_1-\frac{1}{v_2-\frac{1}{\ddots -\frac{1}{v_n}}},\]
with all $v_i\geq 2$, then the quantity $c_1(J)^2$ is given by
\begin{equation}\label{c2=rot}
r^T(Q)^{-1}r,
\end{equation}
where $Q$ is the matrix
\[Q=
\begin{bmatrix}
    -v_1       & 1  \\
 1       & \ddots & \ddots \\
 & \ddots & \ddots & \ddots \\
  & & \ddots &\ddots &1\\
  & & & 1 & -v_n\\
\end{bmatrix}
,\]
which represents the intersection form of $X$ in the basis corresponding to the linear graph, where each vertex is a generator.
According to \cite{honda}, there are two universally tight contact structure on $L(p,q)$ up to isotopy (and just one on $L(p,p-1)$). Honda also characterizes the rotation number of each component of the link given by the chain of Legendrian unknots, whose associated Legendrian surgered manifold is $(L(p,q),	\xi_{st})$.

Let $y=(-v_1+2,-v_2+2,\ldots,-v_n+2)$ be the vector of these rotation numbers, i.e. the vector corresponding to one of the two universally tight (standard) structures on $L(p,q)$, the other one being $-y$. By construction, the rotation vectors representing the virtually overtwisted structures have components $x_i$ satisfying 
\[|x_i|\leq |y_i|,\]
with at least one index $\overline{\imath}$ for which $|x_{\overline{\imath}}|<|y_{\overline{\imath}}|$.

\noindent Consider the function $f:\R^n\to \R$ given by $z\mapsto \|z\|_{Q^{-1}}=z^T(Q)^{-1}z$ and notice that, by Equalities \eqref{c2=rot} and \eqref{c2=sigma},
\[f(y)=\sigma(P)=-n.\]
Theorem \ref{minchi} follows directly from Proposition \ref{concave} below, but first we need:

\begin{lem} \label{matrix}
All the entries of the matrix $Q^{-1}$ are strictly negative (in short: $Q^{-1}\ll 0$).
\end{lem}

\begin{prf}
The condition $Q^{-1}\ll 0$ is true if we show that $Q^{-1}x\ll 0$ holds whenever $x$ is a non-zero vector with non-negative components, i.e. $0\neq x \underline{\gg} 0$ (this is just a consequence of the fact that the columns of $Q^{-1}$ are the images of the vectors of the canonical basis).


So we need to check that: $0\neq x\;\underline{\gg}\; 0$ implies $Q^{-1}x\ll 0$. Rephrased in a different way (using the fact that $Q$ is a bijection), we will show that 
\[0\neq Qx \;\underline{\gg}\; 0 \Rightarrow x\ll 0.\]
The condition $Qx \;\underline{\gg}\; 0$ gives us a system
\[ \begin{cases}
-v_1x_1+x_2\geq 0\\
x_1-v_2x_2+x_3\geq 0\\
\ldots\\
x_{n-1}-v_nx_n\geq 0
\end{cases}
\]
where all the $v_i$'s are $\geq 2$. Let $k$ be an index with 
\[x_k=\max_{i}\{x_i\}_i.\]
We want to show that $x_k<0$. Suppose that $1<k<n$. Then 
\begin{align*}
x_{k-1}-v_kx_k+x_{k+1}&\geq 0 \\
x_{k-1}+x_{k+1}&\geq v_kx_k\\
\end{align*}
and therefore
\[ 2x_k\overset{(a)}\geq x_{k-1}+x_{k+1} \geq v_kx_k\overset{(b)}\geq  2x_k.\]
The inequality $(a)$ follows by the definition of $x_k$, while $(b)$ is true if $x_k\geq 0$ (if it is $<0$ then we would be already done). This implies that $x_{k-1}=x_{k+1}=x_k$ and so we can assume, by iterating this argument, that $k=1$ (the case $k=n$ is the same).
We have:
\begin{align*}
-v_1x_1+x_2&\geq 0 \\
x_2&\geq v_1x_1\\
\end{align*}
and again, as before 
\[x_1\geq x_2\geq v_1x_1\geq 2x_1.\]
Therefore $x_1\geq 2x_1$, so $x_1\leq 0$. To exclude $x_1=0$ just notice that if this were the case, then from $-v_1x_1+x_2\geq 0$ it would follow that $x_2=0$ (being $x_1=0$ the maximum among the $x_i$'s), and consequently all the remaining $x_3=\ldots=x_n=0$, contradicting the assumption $Qx\neq 0$.
\end{prf}

\begin{prop}\label{concave}
For any rotation vector $x$ corresponding to a virtually overtwisted structure (i.e. with components $|x_i|\leq |y_i|$, with at least one strict inequality) we have
\[f(x)>f(y).\]
\end{prop}

\begin{prf}
Inside $\R^n$ we look at the region $D=\{(x_1,\ldots, x_n),\,|\,x_i|\leq |y_i|\,\, \forall i\}$. The goal is to show that the minimum of $\left.f\right|_D:D\to\R$ is realized on the vectors which correspond to the universally tight structures $y$ and $-y$, lying on $\partial D$.

Since $Q$ (and hence $Q^{-1}$) is negative definite, $f$ is concave. Being $f$ a negative definite norm, we know that it is has a unique maximum, which is the origin. Moreover, the minimum of $\left.f\right|_D$ is reached on the boundary $\partial D$. The fact that it is realized on $y$ and $-y$ follows from Lemma \ref{matrix}.
\end{prf}

This implies that the contact structures encoded by the vector $x$ cannot bound any Stein rational ball:

\begin{prf}[of Theorem \ref{minchi}]
Let $(X,J)$ be the Stein filling of $(L(p,q),\xi_{vo})$ described by the Legendrian realization of the chain of unknots associated with the vector of rotation numbers $x$. If $p$ and $q$ are of the form $p=m^2$ and $q=mk-1$ as in the hypothesis, then we know that the universally tight contact structure (corresponding to the rotation vector $y$) admits a Stein rational ball filling it, so $f(y)=\sigma(X)$. By Proposition \ref{concave} we know that $f(x)>f(y)$, hence
\[c_1(X,J)^2=f(x)>f(y)=\sigma(X).\]
Since Equality \eqref{c2=sigma} is not satisfied, $(L(p,q),\xi_{vo})$ does not bound a Stein rational homology ball.
\end{prf}

\section{Coverings of tight structures on lens spaces and applications}\label{coveringsection}

In general, it can be hard to tell if the pullback of a tight contact structure on a 3-manifold along a given covering map is again tight. The situation is much easier if we restrict to lens spaces because of two reasons:
\begin{itemize}
\item[1)] tight structures are classified;
\item[2)] the fundamental group is finite cyclic, hence it is straightforward to determine their coverings.
\end{itemize}
If we start with a virtually overtwisted structure $\xi_{vo}$ on $L(p,q)$ we get an overtwisted structure $\pi^*\xi_{vo}$ on $S^3$, where $\pi:S^3\to L(p,q)$ is the universal cover. If $p$ is a prime number, then this is the only cover that $L(p,q)$ has, otherwise there is a bigger lattice of covering spaces depending on the divisors of $p$. 

Studying the behavior of coverings of a contact 3-manifold gives information about the fundamental group of its fillings:

\begin{thm}\label{pi1}
Let $Y=L(p,q)$ with $p$ prime, so that its fundamental group $\pi_1(Y)\simeq \Z/p\Z$ is simple. Let $\xi$ be a virtually overtwisted contact structure on $Y$, and $(X,J)$ a Stein filling of $(Y,\xi)$. Then $X$ is simply connected.
\end{thm}

\begin{prf}
Let $i:Y\hookrightarrow X$ be the inclusion of the boundary $Y=\partial X$. Being $(X,J)$ a Stein filling of $(Y,\xi)$, the induced morphism $i_*:\pi_1(Y)\to \pi_1(X)$ is surjective. Moreover, by simplicity of $\pi_1(Y)$, we have that $\ker i_*$ can either be:
 
\begin{itemize} 
\item $\ker i_*=1:$ 

In this case, take a finite cover $p:(\widehat{Y},\xi_{ot})\to (Y,\xi)$ for which $\xi_{ot}$ is overtwisted. Call $n$ the degree of such cover. Define the group
\[G=i_*p_*\pi_1(\widehat{Y})\leq \pi_1(X),\]
consider the covering space of $X$ associated to $G$ which is connected and call it $\widehat{X_G}$. Since $\deg(\widehat{X_G}\to X)=n$, we have that $\widehat{X_G}$ is compact.
We are in the case where $i_*$ is an isomorphism and so $\partial \widehat{X_G}$ contains a diffeomorphic copy of $\widehat{Y}$. But by lifting the Stein structure from $X$ to $\widehat{X_G}$ we get a Stein structure on $\widehat{X_G}$ which fills the \emph{connected} contact boundary: note that any Stein semi-filling of a lens space is actually a filling, i.e. its boundary is connected (this comes from the more general result \cite[Theorem 1.4]{os}). Therefore we obtained a Stein filling of $\partial \widehat{X_G}=(\widehat{Y},\xi_{ot})$. This is not possible since the overtwisted contact structures are not fillable (as proved in \cite{eliashberg} and \cite{gromov}).

\item $\ker i_*=\pi_1(Y):$ 

This tells us that $i_*$ is identically zero, and so that, by surjectivity, $\pi_1(X)=1$ as wanted.
\end{itemize} 
\end{prf}

Now we want to study more carefully the behavior of the virtually overtwisted contact structures under covering maps, in order to derive some consequences on the possible fundamental groups of the fillings. The driving condition is the following observation:
let $p:\widehat{Y}\to Y$ be a covering map between compact and connected 3-manifolds, and let $i:Y\hookrightarrow X$ be the inclusion of the boundary $Y=\partial X$. Then, by covering theory:
\begin{align*}
\exists \mbox{ covering }\widehat{X}\to X \mbox{ that r}& \mbox{estricts to a covering } \partial\widehat{X}\to \widehat{Y} \\ 
&\Updownarrow \\
\ker i_* &\leq p_*\pi_1(\widehat{Y}).
\end{align*}
The way we want to apply this is to deduce that $\ker i_*$ should be big enough not to be contained in those subgroups of $\pi_1(Y)$ for which we can associate an \emph{overtwisted} cover. For example, if $X$ is a Stein filling of $Y$ and we are able to construct overtwisted coverings of $Y$ associated to every maximal subgroups of $\pi_1(Y)$, then the kernel of $i_*$ is forced to be the whole $\pi_1(Y)$, being this one the only subgroup of $\pi_1(Y)$ not contained in any maximal subgroups. 
By surjectivity of $i_*$ we would then conclude that $X$ is simply connected.

This looks to be promising because in the case of lens spaces it is easy to determine all the maximal subgroups of the fundamental group. It is nevertheless not so immediate to understand the behavior of the contact structure under the pullback map of a covering, but in certain cases we can use a necessary condition of compatibility of Euler classes to get some results. To better explain this, let us consider the following:

\begin{exmp} Let $(L(34,7),\xi_{vo})$ be obtained by contact $(-1)$-surgery on the Legendrian link of Figure \ref{hopf(5,7)}. If we orient the two components in the counter-clockwise direction we get rotation numbers respectively $+3$ and $-5$.

\begin{figure}[ht!]
\centering
\includegraphics[scale=0.5]{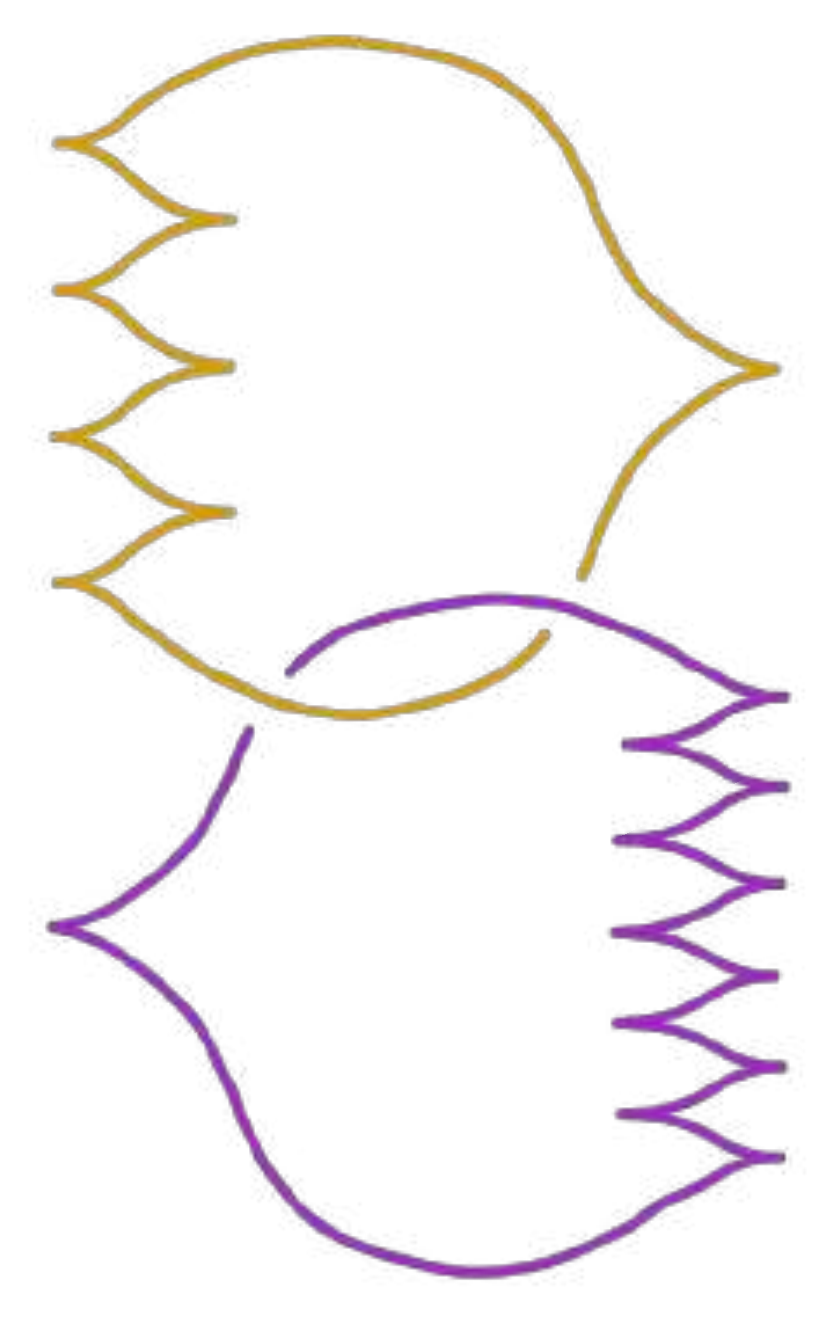}
\caption{Contact $(-1)$-surgery producing $L(34,7)$.}
\label{hopf(5,7)}
\end{figure}

After factoring $34=17\cdot 2$, we see that there are just two coverings:
\[L(17,7)\to L(34,7),\qquad L(2,7)\simeq L(2,1)\to L(34,7).\]
We will show that the given contact structure $\xi_{vo}$ on $L(34,7)$ lifts in both cases to an overtwisted structure. This tells us that, given any Stein filling $X$ of $L(34,7)$, the kernel on the inclusion map at the level of fundamental groups cannot be contained in $\Z/17\Z$ nor $\Z/2\Z$ and therefore is the whole $\Z/34\Z$, so $X$ is necessarily simply connected.
\begin{itemize}
\item The lift of $\xi_{vo}$ to $L(2,1)$ is overtwisted, because the only tight structure on $L(2,1)$ is universally tight and this one pulls backs to the tight structure on $S^3$, but since $\xi_{vo}$ is virtually overtwisted the lift to $S^3$ must be overtwisted.
\item To exclude that $\xi_{vo}$ pulls back to a tight structure on $L(17,7)$ we analyze the possible tight structures supported there. The fraction expansion of $17/7$ is
\[\frac{17}{7}=[3,2,4]\]
and so we see that there are 6 tight structures on $L(17,7)$ up to isotopy (and 3 up to contactomorphism, which are exhibited in Figure \ref{L(17,7)}).

For these structures we compute the Poincaré dual of the Euler class, viewed as an element of $\Z/17\Z\simeq H_1(L(17,7);\Z)$. The previous isomorphism is realized by choosing as a generator the meridional curve $\mu_1$ of the yellow curve with Thurston-Bennequin number $-2$.

Let $\xi$ be any of the three tight contact structures on $L(17,7)$ of Figure \ref{L(17,7)}. 
 
 \begin{figure}[h!]
\centering
\begin{subfigure}[t]{.3\textwidth}
  \centering
  \includegraphics[scale=0.4]{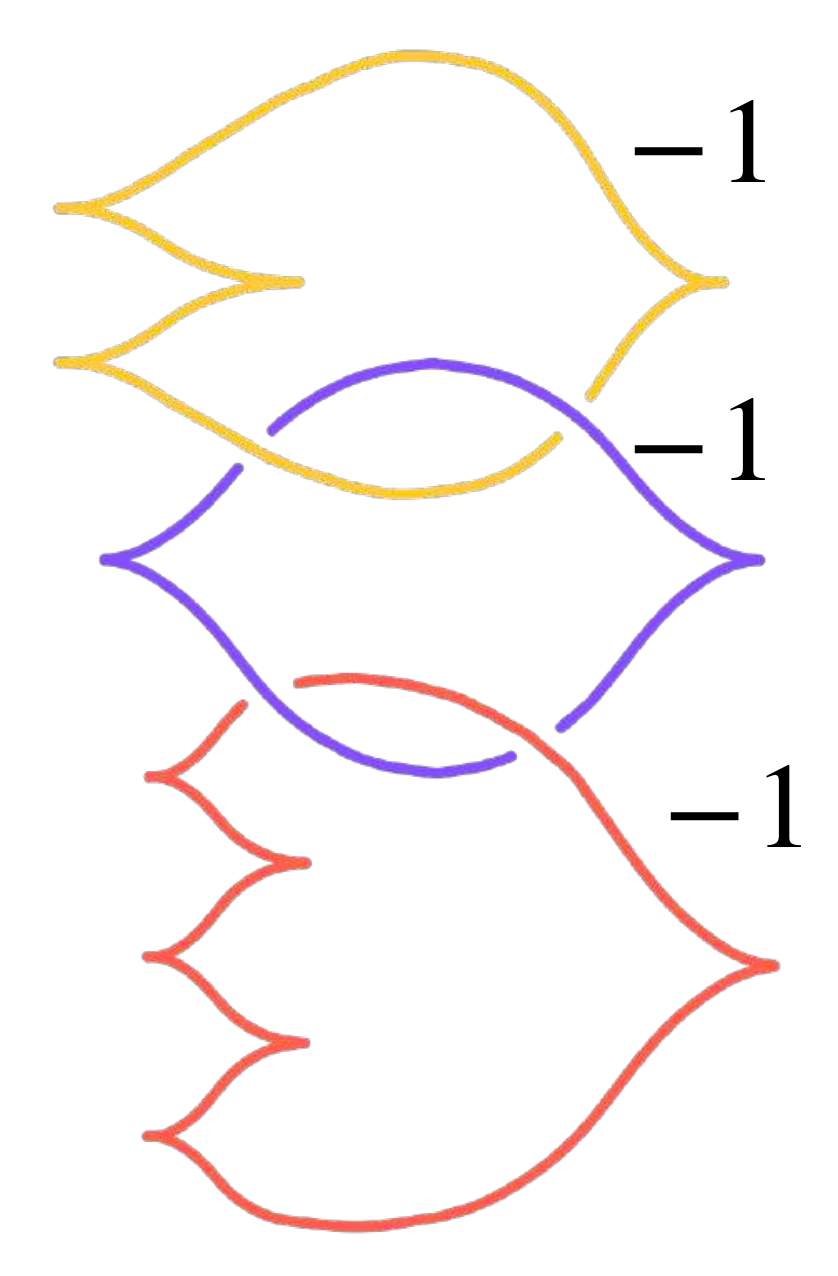}
  \caption{$\xi_1$}
  \label{xi1}
\end{subfigure}%
\begin{subfigure}[t]{.3\textwidth}
  \centering
 \includegraphics[scale=0.4]{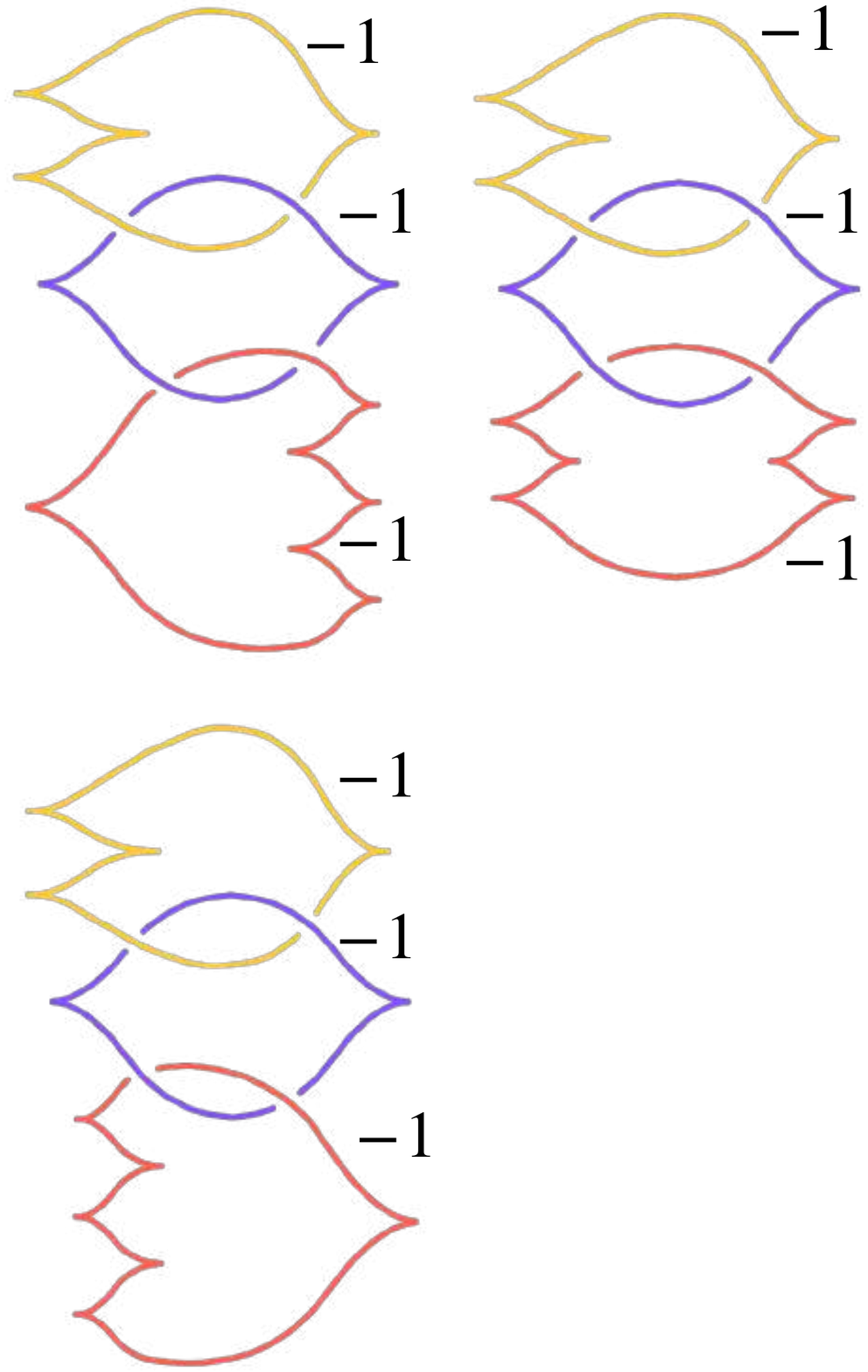}
  \caption{$\xi_2$}
  \label{xi2}
\end{subfigure}
\begin{subfigure}[t]{.3\textwidth}
  \centering
 \includegraphics[scale=0.4]{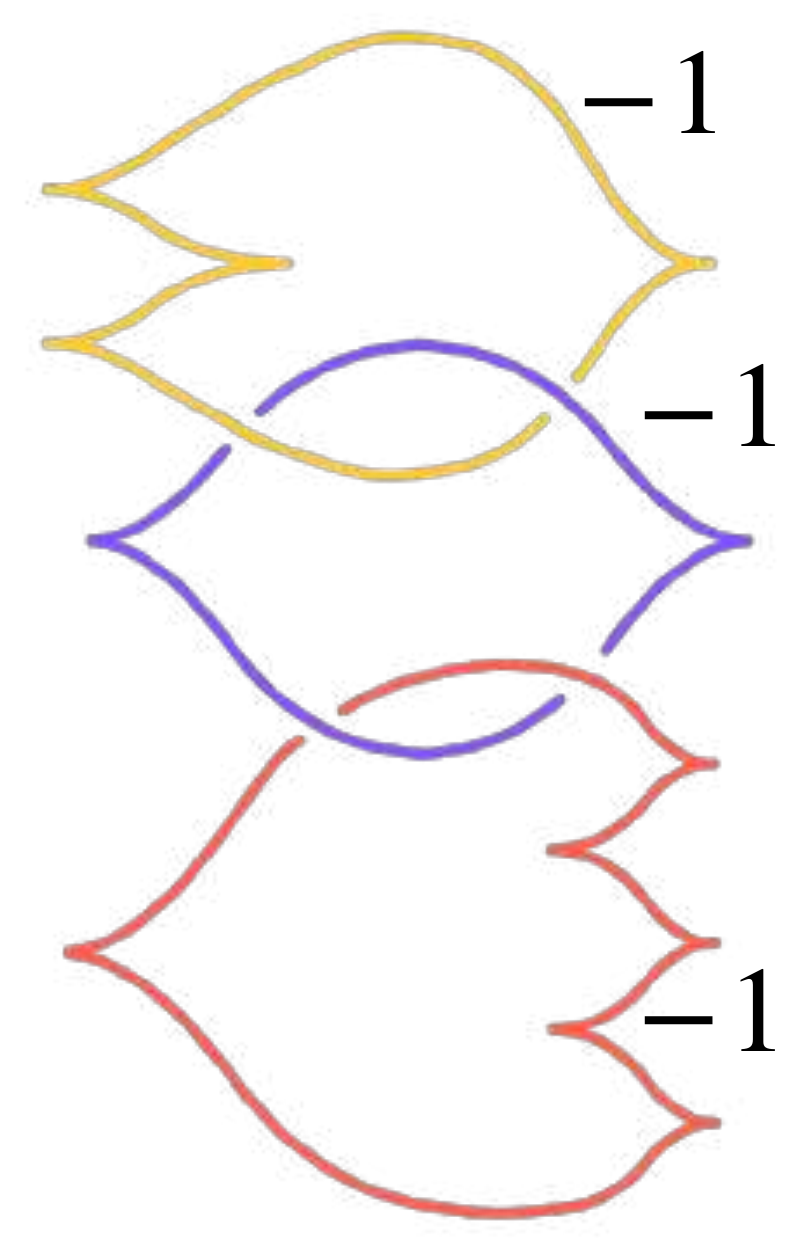}
  \caption{$\xi_3$}
  \label{xi3}
\end{subfigure}
\caption{Tight structures on $L(17,7)$.}
\label{L(17,7)}
\end{figure}

The class $\PD(e(\xi))$ is the image via the boundary map 
\[\partial: H_2(W,\partial W)\to H_1(\partial W)\simeq H_1(L(17,7))\]
of the Poincaré dual of the relative first Chern class of the Stein structure on $W$, where $W$ is the Stein domain described by the corresponding diagram of Figure \ref{L(17,7)}.
The Poincaré dual of the relative first Chern class is (see \cite[Proposition 8.2.4]{ozbagci})
\[\rot(K_1)[D_1,\partial D_1] +\rot(K_2)[D_2,\partial D_2] +\rot(K_3)[D_3,\partial D_3],\]
where the $K_i$'s are the three components of the link and the $[D_i,\partial D_i]$'s are the relative homology classes of the meridian disks of the 4-dimensional 2-handles attached to form the Stein filling $W$. Calling $\mu_i=\partial [D_i,\partial D_i]=[\partial D_i]$, for $i\in \{1,2,3\}$, the meridians of the attaching circles of these handles, we have

\[\PD(e(\xi))=\rot(K_1)\mu_1+\rot(K_2)\mu_2+\rot(K_3)\mu_3.\]

Let $Q$ be the matrix describing the intersection form of $W$, which is the same as the linking matrix
\[  \begin{bmatrix}
   -3 & 1 & 0 \\
 	1 & -2 & 1 \\
   0 & 1 & -4
   \end{bmatrix}.\]
From the exact sequence
\[\xymatrix{ H_2(W) \ar[r]^-Q & H_2(W,\partial W) \ar[r]^-{\partial} & H_1(\partial W)\simeq H_1(L(17,7))},\]
we get three linear relations
\[ \begin{cases}
-3	\mu_1+\mu_2=0\\
\mu_1-2\mu_2+\mu_3=0\\
\mu_2-4\mu_3=0.
\end{cases} \]
which tell us that $\mu_2=3\mu_1$ and $\mu_3=2\mu_2-\mu_1=5\mu_1$.
By putting everything together we get:

\begin{align*}
\PD(e(\xi))&=\partial(\PD(c_1(W,J)))\\
&=\partial (\rot(K_1)[D_1,\partial D_1] +\rot(K_2)[D_2,\partial D_2] +\rot(K_3)[D_3,\partial D_3])\\
&= \rot(K_1)\mu_1+\rot(K_2)\mu_2+\rot(K_3)\mu_3 \\
&=(\rot(K_1)+3\rot(K_2)+5\rot(K_3))\mu_1.
\end{align*}

If we substitute the values of the rotation numbers for the three different contact structures of Figure \ref{L(17,7)} we find:
\[ \PD(e(\xi_1))=11\mu_1, \qquad \PD(e(\xi_2))=\mu_1, \qquad \PD(e(\xi_3))=8\mu_1.\]

The contact structure described in Figure \ref{hopf(5,7)} we started from has 
\[ \PD(e(\xi))=12\mu, \]
with $\mu$ being the meridian of the yellow curve of Figure \ref{hopf(5,7)}.

Notice that $\mu$ is the image of the curve $\mu_1$ under the covering map $p:L(17,7)\to L(34,7)$. This is clear if we take the meridian curves of the single-component unknots with rational framing $-17/7$ and $-34/7$: in this case, the meridian of the curve upstairs is sent to the meridian downstairs, and when we expand from rational to integer surgery representation, we just glue in a series of thickened annuli to the neighborhood of the first component (before the final solid torus is attached), so that previous meridional curves still correspond via the covering map. This is well described in \cite[Section 2.3]{saveliev}. This explains why $\mu_1$ is sent to $\mu$ by the covering map.

At the level of the homology group $H_1$ the covering map is a multiplication by 2 (the degree of the covering) and by naturality we need to find 
\[p_*(\PD(e(p^*(\xi))))=2\PD(e(\xi))\in H_1(L(34,7)).\]
But 
\[2\cdot 11\neq \pm 2\cdot 12 \in \Z/34\Z, \qquad 2\cdot 1\neq\pm  2 \cdot 12 \in \Z/34\Z, \qquad 2\cdot 8\neq\pm  2 \cdot 12 \in \Z/34\Z,\]
therefore we have that none of the three structures of Figure \ref{L(17,7)} is the pullback of our starting structure of Figure \ref{hopf(5,7)}. But those were the only (up to contactomorphism) tight structures on $L(17,7)$, so we conclude that the pullback is necessarily overtwisted, as wanted.

Note that we could have excluded a priori the contact structure $\xi_1$ of Figure \ref{xi1}, this being universally tight. 
\end{itemize}

Similar computations can be done if we start with a Legendrian representation of the Hopf link of Figure \ref{hopf(5,7)1} with rotation numbers $\pm(-3,1)$, $\pm(-3,3)$, $\pm(-3,5)$, $\pm(-1,1)$, $\pm(-1,3)$, $\pm(-1,5)$. We made use of the software \emph{Mathematica} to carry out the computations and check that there is no tight structure on the double cover $L(17,7)$ with compatible Euler class.

\begin{figure}[h!]
\centering
\includegraphics[scale=0.6]{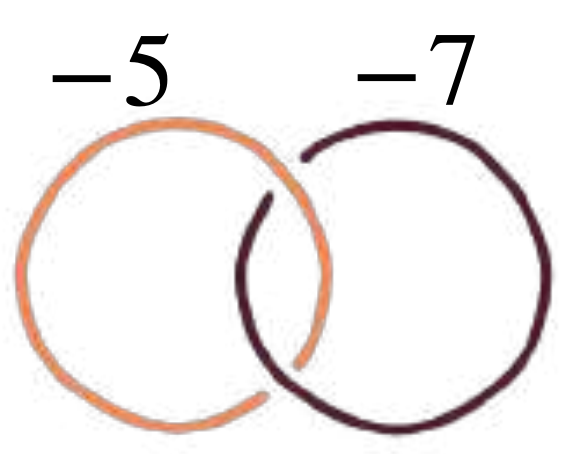}
\caption{Hopf link for $L(34,7)$.}
\label{hopf(5,7)1}
\end{figure}

The fact that the Stein fillings of these virtually overtwisted structure on $L(34,7)$ are simply connected can be deduced, as we just did, simply by looking at the two different coverings. This is something we already knew from the classification of fillings of those lens spaces obtained by contact surgery on the Hopf link, since the fraction expansion of $34/7$ has length 2, see Theorem \ref{theorem1}.
\end{exmp}

\begin{exmp}

Sometimes, an even quicker argument can be used to understand the behavior of a contact structure along certain covering maps. Let's take as an example $L(52,11)$, whose associated fraction expansion has length 3:
\[-\frac{52}{11}=[-5,-4,-3].\]
The two maximal subgroups of $\Z/52\Z$ are $\Z/4\Z$ and $\Z/26\Z$, and again, by running the computation of the Euler classes as above, we can determine which virtually overtwisted contact structure on the base cannot lift to a tight structure. 
But if we look at the covering of degree 13, we find $L(4,11)\simeq L(4,3)$ as total space, and since 
\[-\frac{4}{3}=[-2,-2,-2],\] 
we see that the only tight structure it supports is universally tight. 
Similarly, if we consider the covering $L(13,11)\to L(52,11)$ which has degree 4, we notice that
\[-\frac{13}{11}=[-2,-2,-2,-2,-2,-3]\] 
and hence also $L(13,11)$ supports only universally tight structures, among the tight ones. 
In the covering lattice of $L(52,11)$ it remains to study just the case of $L(26,11)$, for which the behavior can be more subtle (see next section, Theorem \ref{thmL5211}).

\[
\xymatrix{ & (L(13,11),{\color{red}\xi_{ot}})\ar[rr]^-{2:1}& & (L(26,11),{\color{blue}\xi_{?}})\ar[rd]_-{2:1} \\
(S^3,{\color{red}\xi_{ot}})\ar[ru]^-{13:1}\ar[rr]_-{2:1} & & (L(2,11),{\color{red}\xi_{ot}}) \ar[ru]_-{13:1} \ar[rd]_-{2:1}& & (L(52,11),{\color{blue}\xi_{vot}})\\
& & & (L(4,11),{\color{red}\xi_{ot}})\ar[ru]_-{13:1}
}
\]

\end{exmp}

\subsection*{A closer look to the coverings between lens spaces}

The test we made with the Poincaré duals gives only a necessary condition that does not guarantee that the pullback of a given tight contact structure is a tight contact structure simply because characteristic classes match.
So what can be said when there is compatibility between the Euler class of the contact structures of the base and of the covering?
We will try to present the idea of this subsection by starting from an example.

\begin{figure}[h!]
\centering
\begin{subfigure}[t]{.5\textwidth}
  \centering
  \includegraphics[scale=0.5]{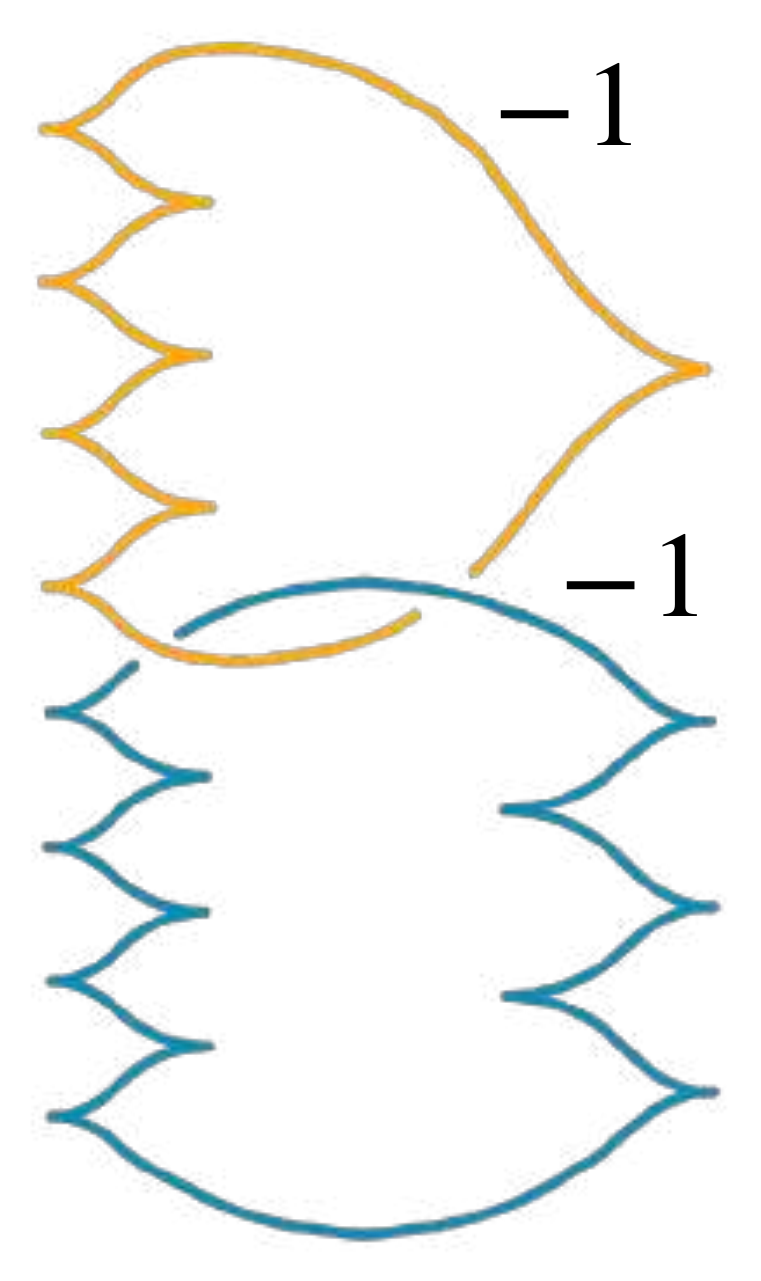}
  \caption{Contact structure on $L(34,7)$.}
  \label{L347}
\end{subfigure}%
\begin{subfigure}[t]{.5\textwidth}
  \centering
 \includegraphics[scale=0.5]{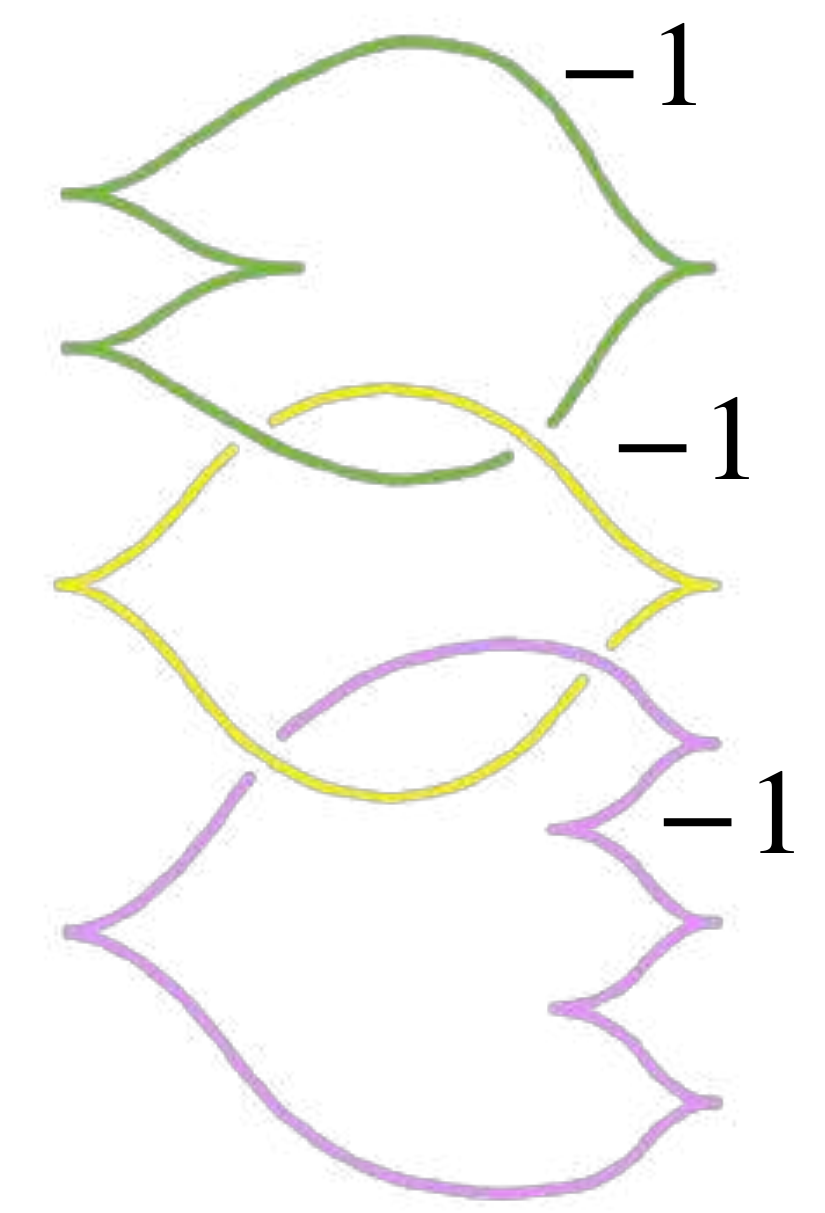}
  \caption{Contact structure on $L(17,7)$.}
  \label{L177}
\end{subfigure}
\caption{ }
\label{covering}
\end{figure}

Again, we choose to describe the double cover of $L(34,7)$. This time we fix the virtually overtwisted structure $\xi$ on $L(34,7)$ where the components of the link have rotation numbers $+3$ and $+1$ respectively, see Figure \ref{L347}. The computation shows that the Poincaré dual of the Euler class of $\xi$ is $+8\in \Z/34\Z$ (via the same identification of $H_1(L(34,7))\simeq \Z/34\Z$ as before). On the double cover $L(17,7)$ we take the tight structure $\hat{\xi}$ corresponding to the rotation vector $(1, 0, -2)$, as showed in Figure \ref{L177}. 

\noindent By running the computation, we find that $\PD(e(\hat{\xi}))=+8\in \Z/17\Z$, so that the covering map 
\[p:L(17,7)\to L(34,7)\]
takes $\PD(e(\hat{\xi}))$ to $16=2\PD(e(\xi))$, as it should certainly happen if $\hat{\xi}$ were isotopic to $p^*\xi$. But we will show that this is not the case, and argue that $p^*\xi$ is instead \emph{overtwisted}.

To do this, we need to use the description of tight structures on lens spaces of \cite{honda}, which we recall after the following definition:
\begin{defn}\label{convexdivid}
Given a contact 3-manifold $(Y,\xi)$, a \emph{contact vector field} $v$ on $Y$ is a vector field whose flow preserves the contact planes. A smooth surface $\S\subseteq (Y,\xi)$ is \emph{convex} if there exists a contact vector field $v$ on $Y$ transverse to $\S$.
The \emph{dividing set} of $v$ on $\S$ is defined as 
\[\Gamma=\{x\in\S\;\mid\;v(x)\in\xi_x\}.\]
\end{defn}

Giroux proved in \cite{girouxconv} that the dividing set is a 1-dimensional submanifold, whose isotopy type is independent of the choice of the contact vector field. We now focus on the case when $\S=T^2$. If the contact structure is tight in a neighborhood the torus, then the diving set for a convex torus consists of an even number of parallel circles. By identifying $T^2$ with $\R^2/\Z^2$, we can talk about the \emph{slope} of these circles as a pair of numbers, which depends on the choice of the identification: when $T^2=\partial D^2\times S^1$, we use the meridian curve as one direction.

\paragraph{Honda's algorithm.}

\begin{figure}[h!]
\centering
\includegraphics[scale=0.5]{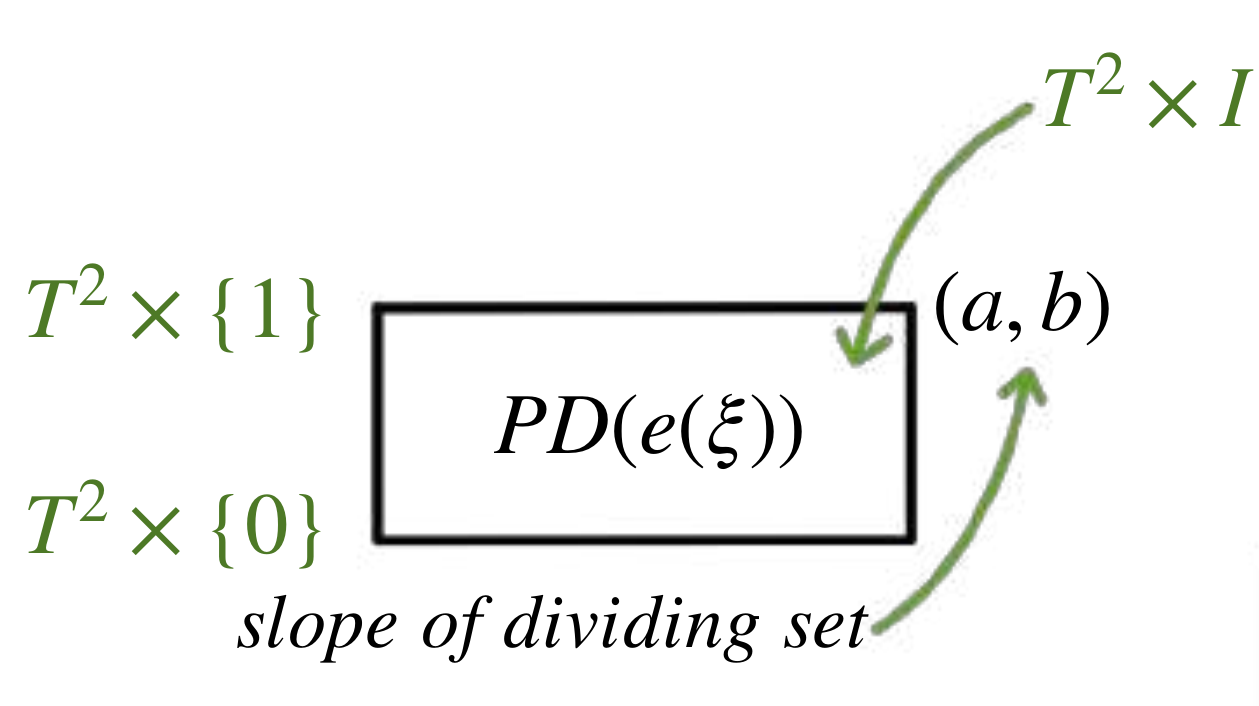}
\caption{Anatomy of a basic slice.}
\label{basicslice}
\end{figure}

In \cite[Section 4.3]{honda} it is explained how to cut a lens space, endowed with a tight contact structure, into two \emph{standard solid tori} and other pieces called \emph{basic slices}. With \emph{standard solid torus} we mean a small tubular neighborhood of a Legendrian knot, with standard coordinates on its boundary, see \cite[Section 2]{etnyrenotes}. On the other hand, a \emph{basic slice} is an oriented thickened torus $T^2\times I$ with a tight contact structure on it, such that
\begin{itemize}
\item the two boundary components are convex;
\item the minimal integral representatives of $\Z^2$ corresponding to the slopes at the extremes form a $\Z$–basis of $\Z^2$;
\item every convex torus parallel to the boundary has slope between the slopes of the extremes.
\end{itemize}
Each basic slice supports a unique tight contact structure, up to contactomorphism, but up to isotopy there are two classes: the isotopy class is determined by the sign of the (Poincaré dual of the) Euler class of the contact structure restricted to that basic slice. We always assume that the boundary tori are oriented according to the initial orientation on $T^2\times I$. A schematic picture of a basic slice is represented in Figure \ref{basicslice}.

\noindent The contact structure on the lens space is then encoded in the sequence of slopes on each basic slice and in the corresponding signs. We recall how the algorithm of Honda works for the lens space $L(p,q)$:

\noindent we start from the expansion $-p/q=[-a_1,\ldots,-a_n]$, with $a_i\geq 2$ for every $i$. Then we compute
\begin{align*}
-p_1/q_1:= & [-a_1,\ldots,-a_n+1] \\
-p_2/q_2:= & [-a_1,\ldots,-a_n+2] \\
-p_3/q_3:= & [-a_1,\ldots,-a_n+3] \\
\vdots
\end{align*}
until we get, after $k=a_n-1$ steps, to a rational number such that the length of its expansion is a number $m$, smaller than $n$, say
\[-p_k/q_k= [-b_1,\ldots,-b_m]. \]
This first set of numbers $\{-p_1/q_1,\ldots,-p_k/q_k\}$ will constitute the first block. Then we continue
\begin{align*}
-p_{k+1}/q_{k+1}:= & [-b_1,\ldots,-b_m+1] \\
-p_{k+2}/q_{k+2}:= & [-b_1,\ldots,-b_m+2] \\
\vdots
\end{align*}
until we get, after $h=b_m-1$ steps, to a rational number $-p_h/q_h$ such that the length of its expansion is less than $m$. This set of numbers $\{-p_k/q_k,\ldots,-p_h/q_h\}$ will constitute the second block. We go on this way until we reach the rational number $-1/1$. 

In total, we will produce an ordered set of blocks of ordered rational numbers which increase from $-p_1/q_1$ to $-1/1$, such that the numbers in each block have an associated continued fraction expansion of the same length. The boundary numbers, i.e. those which determine a change of length, appear twice: once at the bottom of a block, and then immediately after at the top of the following block. For example, the number $-14/3$ can appear twice, once as $[-5,-4,-1]$ and once as $[-5,-3]$. The expansion $[-5,-4,-1]$ determines the end of the block with length 3, while $[-5,-3]$ determines the start of the block of length 2. We record these rational numbers $-p_i/q_i$ as pairs of coprime integers $(-q_i,p_i)$. These numbers correspond to the slope of the dividing sets of the contact structure under analysis, when restricted to the corresponding basic slice.

Then we remove a standard torus from the lens space $(L(p,q),\xi)$ and we picture what is left in the following way: we draw the basic slices starting from the slope $-p_1/q_1$ until $-1/1$, divided into the blocks as described above. At the end of this thickened torus we draw the other basic torus. 

As we explained, every basic slice comes equipped with boundary slopes described by two rational numbers, which are represented by pairs $(-q,p)$ and $(-q',p')$. Honda proved that taking the difference of these values gives the Poincaré dual of the Euler class restricted there, up to sign, as an element of $H_1(T^2)\simeq \Z\oplus\Z$, written in the basis $(\partial D^2,S^1)$ specified by the lower solid torus (see \cite[Section 4.7.1]{honda}). As mentioned above, the isotopy class of the unique (up to contactomorphism) contact structure on each basic slice is specified by the sign of the restriction of the Poincaré dual of the Euler class. Within a single block of basic slices, the only thing that matters is how many positive and negative signs we have, but not where these are placed: this is a consequence of a property of shuffling, which says that rearranging the signs within a block gives an isotopic contact structure, see \cite[Section 4.4.5]{honda}. This is coherent with the fact that, when drawing a Legendrian unknot with its stabilizations, we don't need to remember if we first stabilized positively or negatively, but just the final result. 

\textbf{To sum up}: if we start from a chain of $n$ Legendrian unknots, we get $n$ blocks (one for every component) of $a_i-2$ basic slices each (where $-p/q=[-a_1,\ldots ,-a_n]$). Every positive/negative stabilization that we see in the Legendrian link corresponds to a plus/minus sign in the corresponding block. Notice that when a coefficient in the expansion is $-2$, then its corresponding block will be empty, reflecting the fact that there is no choice of placing stabilizations in a Legendrian knot with Thurston-Bennequin number $-1$.

\vspace{0.5cm}

\begin{figure}[h!]
\centering
\begin{subfigure}[t]{.5\textwidth}
  \centering
  \includegraphics[scale=0.5]{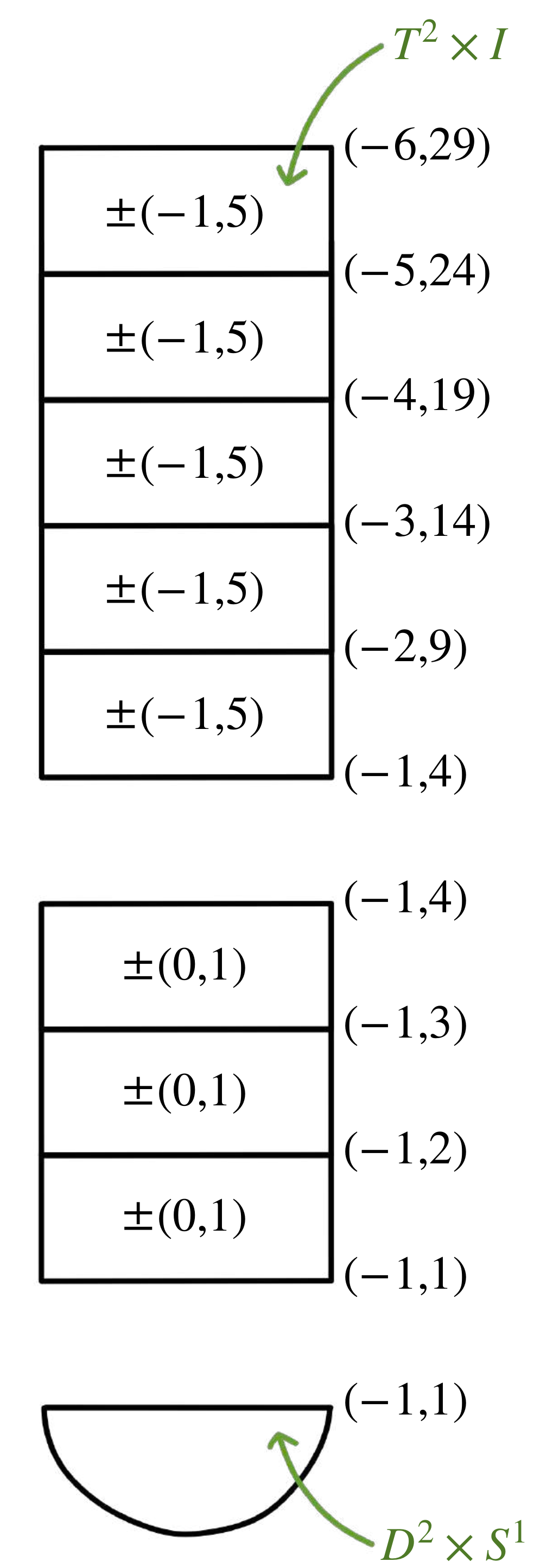}
  \caption{Subdivision into basic slices.}
  \label{blocks}
\end{subfigure}%
\begin{subfigure}[t]{.5\textwidth}
  \centering
 \includegraphics[scale=0.5]{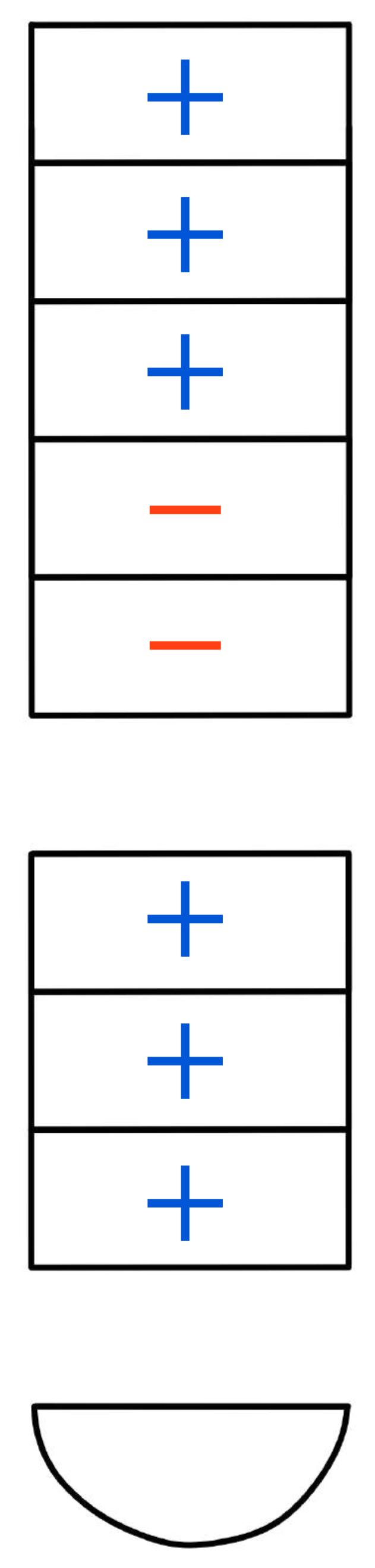}
  \caption{Signs of the basic slices.}
  \label{signs}
\end{subfigure}
\caption{Description of the contact structure $\xi$ on $L(34,7)$.}
\label{blockssigns}
\end{figure}

For example, see Figure \ref{blocks}: the algorithm applied to $L(34,7)$ gives two blocks of 5 and 3 basic slices respectively, where the slopes of the dividing sets on the boundary are indicated there.

From this picture it is also easy to calculate the Poincaré dual of the Euler class of the structure we choose according to the signs of each basic slice (indicated with colors blue and red in Figure \ref{signs}). By capping off with the solid torus below, we make the first $S^1$-factor of $T^2\times I$ nullhomologous, so we can just focus on the second entry in homology. The Poincaré dual of the Euler class of the structure is finally understood in the first homology group of the lens space once we glue the other solid torus (above). The structure in Figure \ref{blockssigns} has $\PD(e)$ given by 
\begin{equation}\label{PDdown}
5+5+5-5-5+1+1+1=8\in \Z/34\Z,
\end{equation}
and it is exactly the one resulting from contact $(-1)$-surgery on the Legendrian Hopf link of Figure \ref{L347}, where the component with Thurston-Bennequin number $-4$ has rotation $+3$ (corresponding to the three pluses in the lower block), and the other one has rotation $+1$ (corresponding to the upper block with three pluses and two minuses).

\begin{figure}[h!]
\centering
\begin{subfigure}[t]{.5\textwidth}
  \centering
  \includegraphics[scale=0.5]{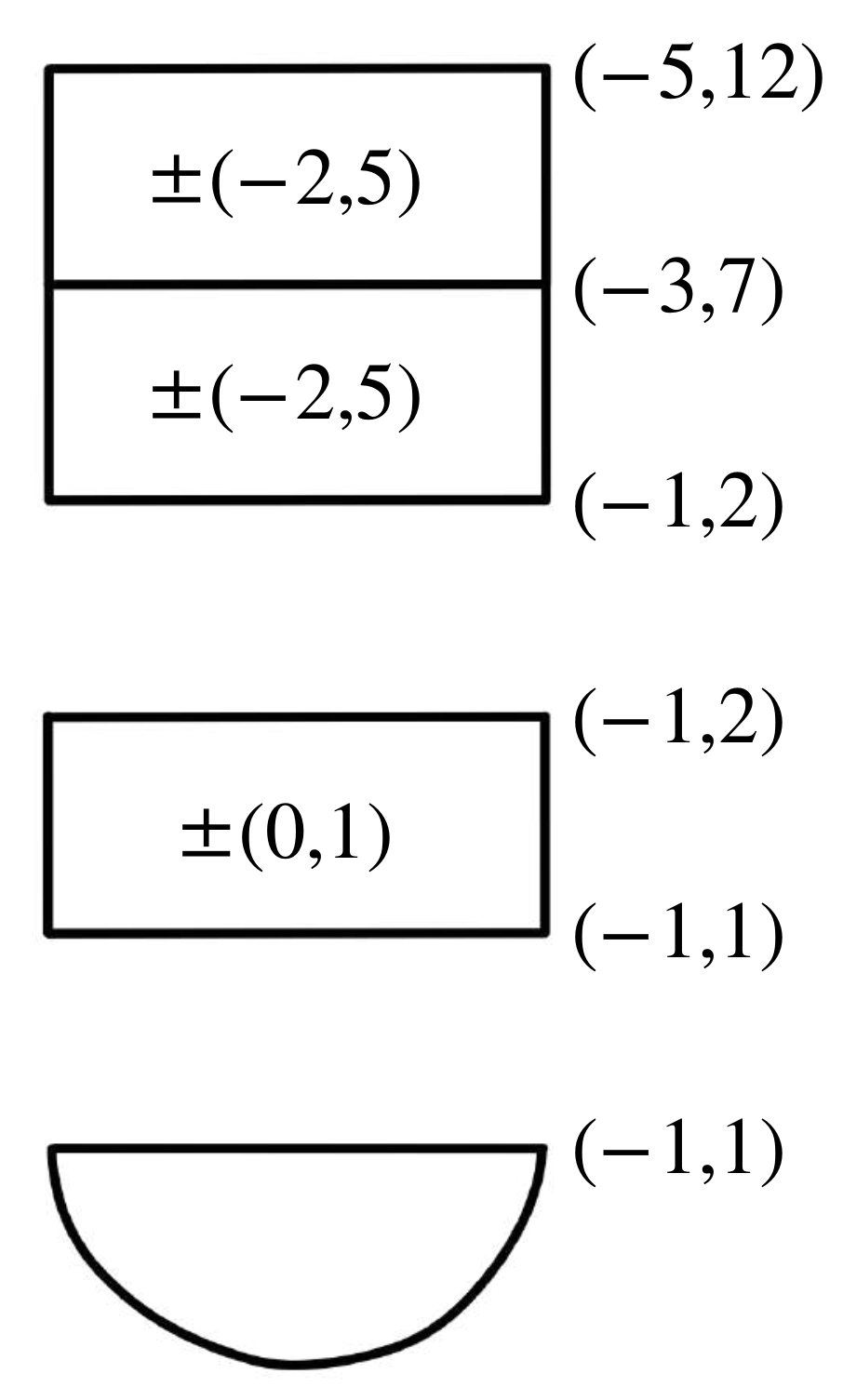}
  \caption{Subdivision into basic slices.}
  \label{blocks2}
\end{subfigure}%
\begin{subfigure}[t]{.5\textwidth}
  \centering
 \includegraphics[scale=0.5]{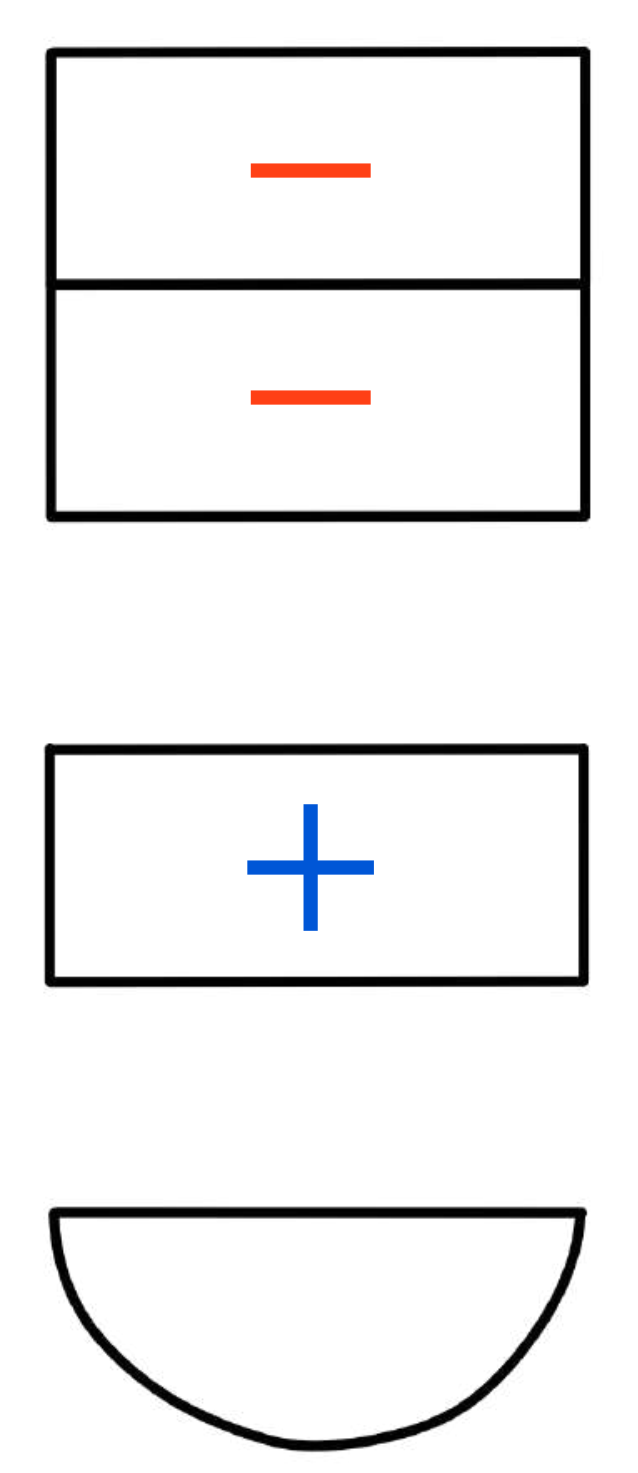}
  \caption{Signs of the basic slices.}
  \label{signs2}
\end{subfigure}
\caption{Description of a contact structure on $L(17,7)$.}
\label{blockssigns2}
\end{figure}

Now we look at the double covering map, which, on every basic slice, looks like
\[(z,w)\to (z,w^2),\]
where $z$ is the coordinate corresponding to $\partial D^2$ (which will be capped off when the lower solid torus is glued), and $w$ is the coordinate of the other $S^1$-factor.
We split $L(17,7)$ with a tight structure into two solid tori: the first one is pictured in Figure \ref{blocks2}, and subdivided into a block of two basic slices, plus a single basic slice, plus a standard solid torus; the other solid torus is a standard torus which will be glued on top of the uppermost basic slice and which is not pictured. 

All the tight contact structures on $L(17,7)$ are encoded in the decomposition of the represented solid torus into these pieces: by choosing the sign of the basic slices we produce all the different (up to isotopy) 6 tight structures that $L(17,7)$ supports. As a double check, one can think at the different Legendrian representatives of the 3-components link made by a chain of unknots with Thurston-Bennequin numbers $-2$, $-1$ and $-3$ (notice indeed that $-17/7=[-3,-2,-4]$). The candidate tight contact structure on $L(17,7)$ which should be the pullback of the one on $L(34,7)$ described by Figure \ref{blocks} has the single basic slice with positive sign, and the other two in the block with negative signs. This corresponds to the choice of the rotation numbers for the components of the link to be $+1$, $0$ and $-2$: the link on which contact $(-1)$-surgery should give the pullback structure on $L(17,7)$ along the covering map is pictured in Figure \ref{L177}. The reason why this is the correct candidate is because, among the 3 different (up to contactomorphism) contact structures on $L(17,7)$, this is the only case where we have compatibility of Euler classes: the computation (which can be performed in two different ways) shows that the Poincaré dual of the Euler class upstairs is $-9\equiv 8 \pmod{17}$, which gets sent to $16=2\cdot 8 \in \Z/34\Z$, which, as we already computed in Equation \eqref{PDdown}, is the double of the Poincaré dual of the Euler class downstairs.

\begin{figure}[h!]
\centering
\begin{subfigure}[t]{.5\textwidth}
  \centering
  \includegraphics[scale=0.5]{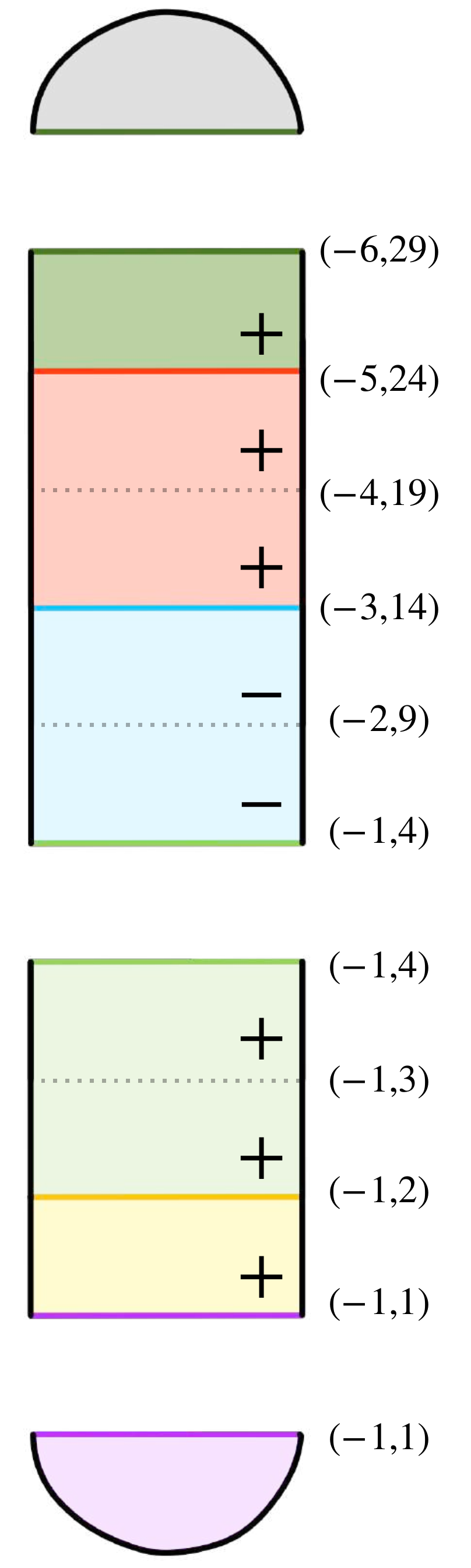}
  \caption{$L(34,7)$}
  \label{colors347}
\end{subfigure}%
\begin{subfigure}[t]{.5\textwidth}
  \centering
 \includegraphics[scale=0.5]{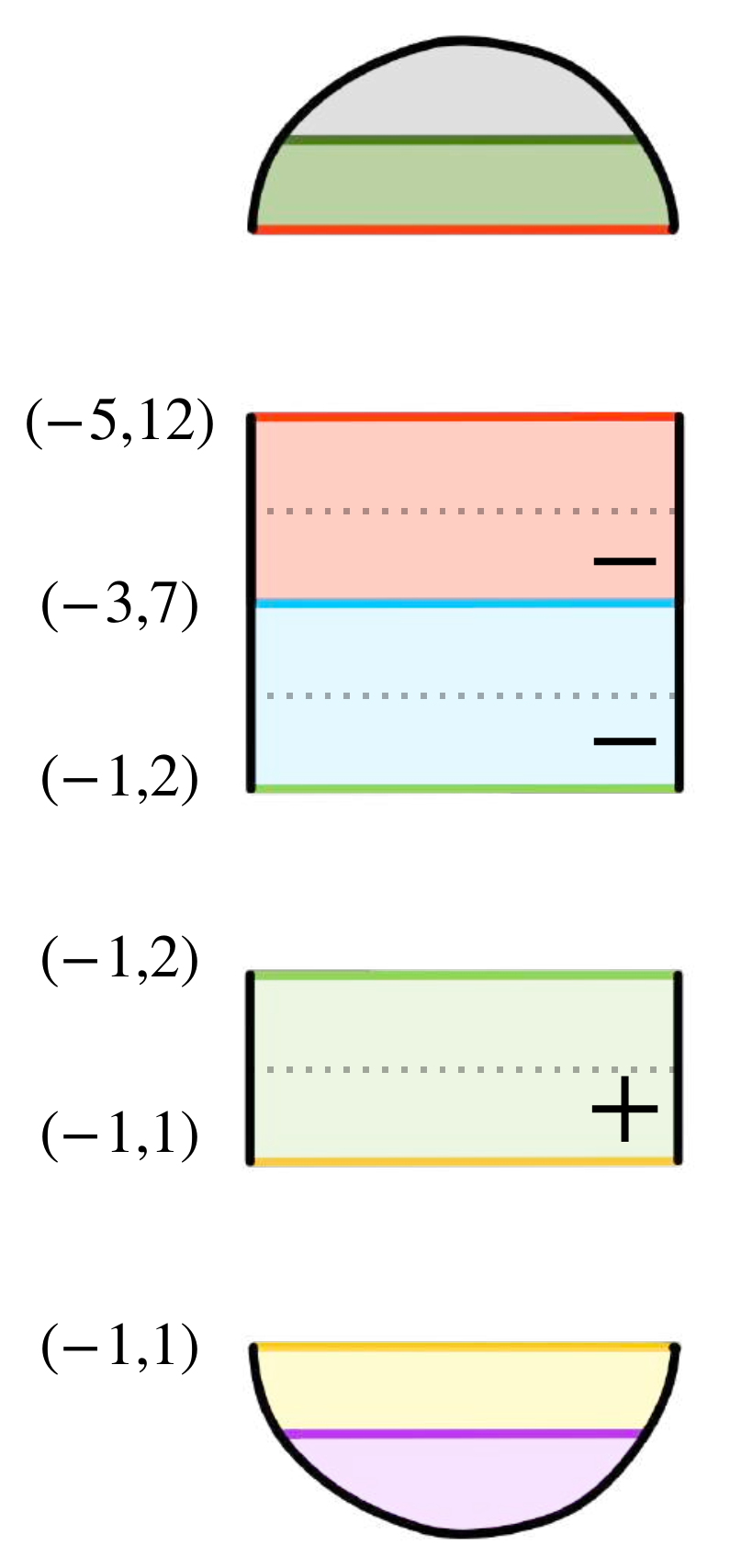}
  \caption{$L(17,7)$}
  \label{colors177}
\end{subfigure}
\caption{Behavior of slices under the covering map.}
\label{colors}
\end{figure}

But now we argue that there cannot be compatibility in the signs of the basic slices of $L(34,7)$ and $L(17,7)$. Indeed, once a sign for a basic slice downstairs is chosen, then its lift should have the same sign, see \cite[Section 1.1.4]{honda2}.
By lifting the dividing sets of the various convex tori we see where the different basic slices go: Figure \ref{colors} is describing this by means of colors. Computations show that the lowest basic slice of $L(34,7)$ is pulled back inside the standard torus, and the same is true for the uppermost slice. Therefore the behavior of the contact structure upstairs is regulated by what happens to the central slices, i.e. from the yellow line $(-1,2)$ to the red line $(-5,24)$.

But here we finally see the contradiction. While:
\begin{itemize}
\item[1)] the positive slices from yellow $(-1,2)$ to green $(-1,4)$ lift to a positive slice in $L(17,7)$ and
\item[2)] the negative slices from green $(-1,4)$ to blue $(-3,14)$ lift to a negative slice in $L(17,7)$, we have that
\item[3)]  the \emph{positive} slices from blue $(-3,14)$ to red $(-5,24)$ lift to a \emph{negative} slice in $L(17,7)$
\end{itemize}
and this is not possible. No matter how we decide to shuffle the basic slices in each single block (see \cite[Section 4.4.5]{honda}), we always end up with a contradicting situation (as proved in Theorem \ref{thmcoverovertwisted}).

This tells us that, even if there is a tight virtually overtwisted structure on $L(17,7)$ whose Euler class is compatible with the  structure $\xi$ we chose on $L(34,7)$, the pullback of $\xi$ along the double covering map is \emph{overtwisted}, as claimed.

\begin{thm}\label{thmcoverovertwisted} Any virtually overtwisted structure on $L(34,7)$ lifts to an overtwisted one along the double cover 
\[L(17,7)\to L(34,7).\]
\end{thm}

\begin{prf}
We argue here using the behavior of the basic slices described in Figure \ref{colors}.
Look at the three basic slices in $L(17,7)$, Figure \ref{colors177}, regardless of the signs. Call $\hat{\xi}$ the pullback of a given $\xi$ on $L(34,7)$ and compare the Poincaré dual of their Euler classes. Assuming that the structures are both tight, we see that the choice of the sign of the red basic slice in $L(17,7)$ contributes to a $\pm 5$ for $\PD(\hat{\xi})$ and, pushed down, to a $\pm 10$ for $\PD(\xi)$. The same is true for the light blue slice, while the green slice gives a $\pm 1$ for $\PD(\hat{\xi})$ and a $\pm 2$ for $\PD(\xi)$. Moreover, inside $L(34,7)$ we have two extra slices (dark green and yellow in Figure \ref{colors347}), whose signs can be chosen independently. Requiring compatibility of Euler classes means to impose
\[\PD(\xi)\equiv \PD(\hat{\xi})\;\pmod{17}.\]
Therefore, according to what we have just said:
\[{\color{red}\pm 10}\, {\color{cyan}\pm 10} \, {\color{green}\pm 2}\,  {\color{olive}\pm 5}\, {\color{orange}\pm 1}\equiv {\color{red}\pm 5}\,	{\color{cyan}\pm 5} \, {\color{green}\pm 1} \;\pmod{17}\]
which is the same as
\[{\color{red}\pm 5}\, {\color{cyan}\pm 5} \, {\color{green}\pm 1}\,  {\color{olive}\pm 5}\, {\color{orange}\pm 1}\equiv 0\;\pmod{17}.\]
Clearly, this can be done only in two ways, namely by choosing all pluses or all minuses. And these correspond exactly to the two universally tight structures, for which we already knew that there is compatibility. Therefore, among the virtually overtwisted structures there cannot be a coherent choice of signs resulting in compatible Euler classes.
\end{prf}

\begin{thm}\label{thmL5211} Any virtually overtwisted structure on $L(52,11)$ lifts to an overtwisted one along all of its non-trivial covers.
\end{thm}

\begin{prf} At the end of previous section we argued that in the covering lattice of $L(52,11)$ the only case which was more subtle to describe was the double cover
\[L(26,11)\to L(52,11),\]
because otherwise we already knew that virtually overtwisted structures on the base would lift to overtwisted structures. We analyze this remaining case as we did before, by looking for compatibility between the signs of the basic slices and the count of the possible Euler classes.
Figure \ref{L5211L2611} shows where the basic slices go, from $L(26,11)$ to $L(52,11)$. 
\begin{figure}[h!]
\centering
\begin{subfigure}[t]{.5\textwidth}
  \centering
  \includegraphics[scale=0.5]{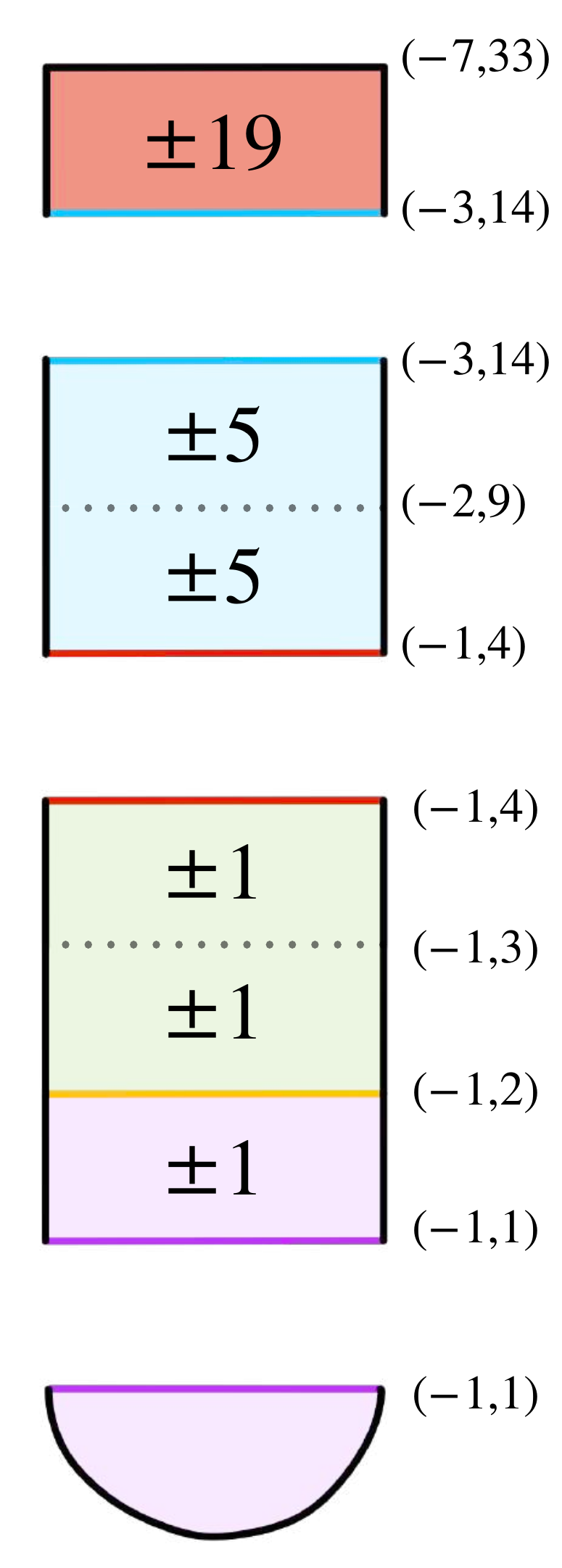}
  \caption{$L(52,11)$}
  \label{blocksL5211}
\end{subfigure}%
\begin{subfigure}[t]{.5\textwidth}
  \centering
 \includegraphics[scale=0.5]{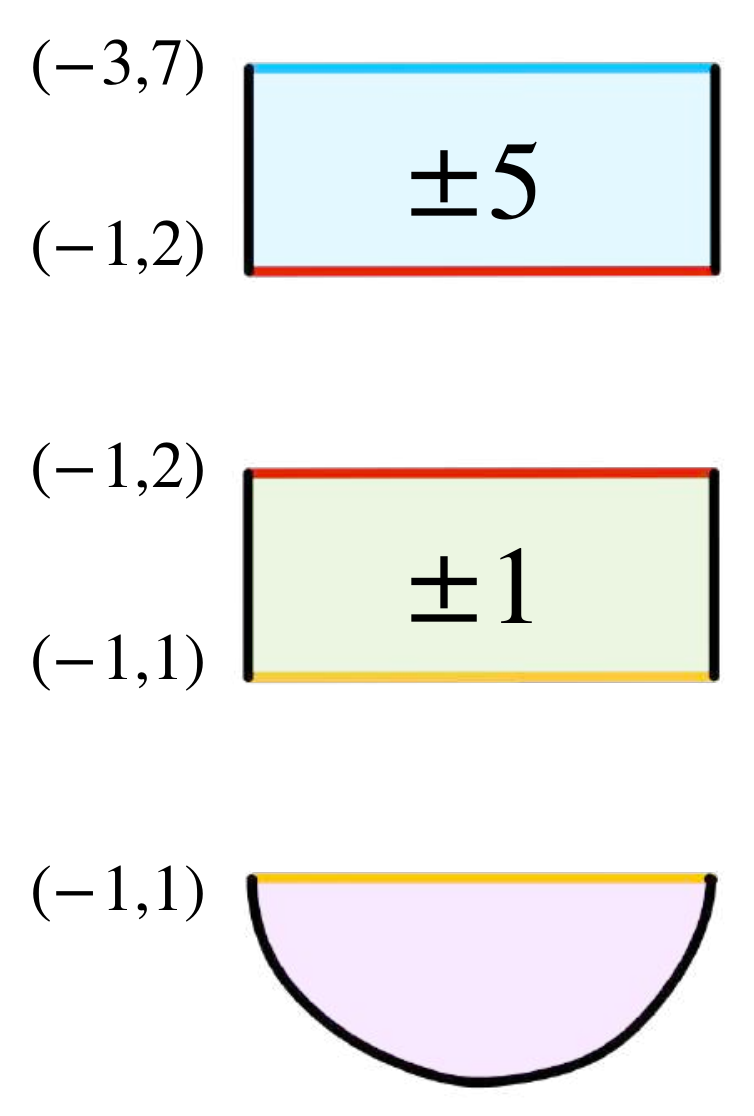}
  \caption{$L(26,11)$}
  \label{blocksL2611}
\end{subfigure}
\caption{Behavior of slices under the covering map.}
\label{L5211L2611}
\end{figure}
The count of the Poincaré duals of the two Euler classes gives
\[{\color{cyan}\pm 10}\, {\color{green}\pm 2} \, {\color{red}\pm 19}\,  {\color{violet}\pm 1} \equiv {\color{cyan}\pm 5}\,	{\color{green}\pm 1} \;\pmod{26}\]
which is the same as
\[{\color{cyan}\pm 5}\, {\color{green}\pm 1} \, {\color{red}\pm 19}\,  {\color{violet}\pm 1} \equiv 0 \;\pmod{26}.\]
Again, we see that this can be done only in two ways, namely by choosing all pluses or all minuses, which correspond exactly to the two universally tight structures. Therefore, among the virtually overtwisted structures there cannot be a coherent choice of signs resulting in compatible Euler classes.
\end{prf}

We can finally give a proof of Theorem \ref{thmcovering}, which stated that every virtually overtwisted contact structure on $L(p,q)$ lifts along a degree $d$ covering to a structure which is overtwisted, if $q<p<dq$. 

\begin{prf}[of Theorem \ref{thmcovering}]
If the pullback of the contact structure were tight, it should fit with the description of tight structures according to the basic slices subdivision. We claim that the lower solid torus $\mathcal{H}_1$ until the level $-p/q$ gets all pulled back into the standard solid torus whose dividing set has slope $-1/1$. This comes from the fact that the curve with slope $(-q,p)$ pulls back to the one with slope $(-dq,p)$, according to the behavior
\[\xymatrix{
-\frac{p}{dq}\ar@{|->}[r]^-{\cdot d}&-\frac{p}{q}. }\]
But by assumption $-\frac{1}{1}<-\frac{p}{dq}$, and since the slopes of the dividing sets are increasing when read from top to bottom (compare with Figure \ref{blocks}), the claim follows.

\noindent Since we are considering a virtually overtwisted structure on $L(p,q)$, the pullback of $\mathcal{H}_1$ \emph{cannot be tight}, otherwise it would be universally tight, being it a subset of a solid torus in standard coordinates (which does not support virtually overtwisted structures). Therefore, we must have here an overtwisted disk, as wanted.
\end{prf}

\begin{cor} Let $p_1,\,p_2$ be prime numbers, not necessarily distinct, and let $q$ be an integer such that $p_i<q<p_1p_2$ for $i=1,\,2$. Then each non-trivial covering of $(L(p_1p_2,q),\xi_{vo})$ is overtwisted, for any virtually overtwisted structure $\xi_{vo}$.
\end{cor}

\begin{prf}
It is a direct consequence of previous theorem, since 
\[p_1p_2<dq,\]
where $d$ is either $p_1$ or $p_2$ (which are the only possible degrees for a non-trivial covering).
\end{prf}

\begin{nrem} The hypothesis of Theorem \ref{thmcovering} can be relaxed by just requiring that $p'<dq'$, where $p'$ and $q'$ are determined as follows: let 
\[-\frac{p}{q}=[-a_1,\ldots ,-a_n]\]
be the continuous fraction expansion, with $a_i\geq 2$ for each $i=1,\ldots, n$. Then define $p'$ and $q'$ as
\[[-a_1,\ldots ,-a_n+1]=-\frac{p'}{q'}.\]
In this way we have
\[-\frac{p}{q}<-\frac{p'}{q'}\]
so that the requirement $-1<-p'/q'$ is less restrictive. The reason why Theorem \ref{thmcovering} stays true with this weaker assumption is that the description of a contact structure via basic slices shows as the smallest slope (hence on top of the uppermost block) precisely the slope $-p'/q'$ (see \cite[Section 4.6]{honda}). To ask that, from this level down, the solid torus is pulled back inside the standard torus in the covering guarantees the existence of an overtwisted disk in the covering space, as argued in the proof of Theorem \ref{thmcovering}. 

There is another description of the two numbers $p'$ and $q'$ which is intrinsic in the sense that does not involve the computation of the continued fraction expansion: given $p$ and $q$, let $q^*$ be the multiplicative inverse of $q$, modulo $p$, i.e. $0<q^*<p$ and
\[q^*q\equiv 1\pmod{p}.\]
If we put $p'=p+q^*$, then $q'$ is the multiplicative inverse of $q^*$, modulo $p'$, i.e.
\[q'q^*\equiv 1\pmod{p'}.\]
\end{nrem}

\subsection*{Comparing $\pi_1$ and $\chi$ of a filling}

The goal of this section is to see some applications to concrete examples of Theorem \ref{pi1b2}, which is proved below.

\begin{prf}[of Theorem \ref{pi1b2}] 
Take the universal covering $\widetilde{X}\to X$, which has degree $d$, whose boundary is the (connected) covering $L(p',q')\to L(p,q)$ of  degree $d$.
The Euler characteristics of the fillings satisfy
\[\chi(\widetilde{X})=d\chi(X)\]
and hence, by Theorem \ref{maxchi},
\[\chi(X)=\frac{\chi(\widetilde{X})}{d}\leq \frac{1+l'}{d}.\]
\end{prf}

\begin{cor} Let $X$ be a Stein filling of a lens space $L(p,q)$ with a virtually overtwisted structure, and let $d$ be a divisor of $p$. If
\[2d>1+\length((p/d)/q),\]
then the fundamental group of $X$ cannot be $\Z/d\Z$.
\end{cor}

\begin{prf}
It follows by contradiction from Theorem \ref{pi1b2} if we look at the associated $d$-covering $\widehat{X}\to X$ and remember that $2\leq \chi(X)$, as proved in Section \ref{lowerboundsection} and \cite[Proposition A.1]{gollastar}. The number $\length((p/d)/q)$ has to computed after reducing $q$ modulo $p/d$.
\end{prf}

\vspace{0.7cm}

Sometimes, depending on the arithmetic of the rational numbers, it happens that the behavior of the basic slices of a covering is never compatible with the choice of signs determining the Euler classes, and this guarantees the covering itself to be overtwisted, which in turn implies that all the fillings are simply connected. But there are cases where a non-trivial cover of a tight virtually overtwisted structure stays as such, and so we need other arguments to calculate the fundamental group of a filling.

\begin{figure}[ht!]
\begin{subfigure}[t]{.5\textwidth}
  \centering
  \includegraphics[scale=0.5]{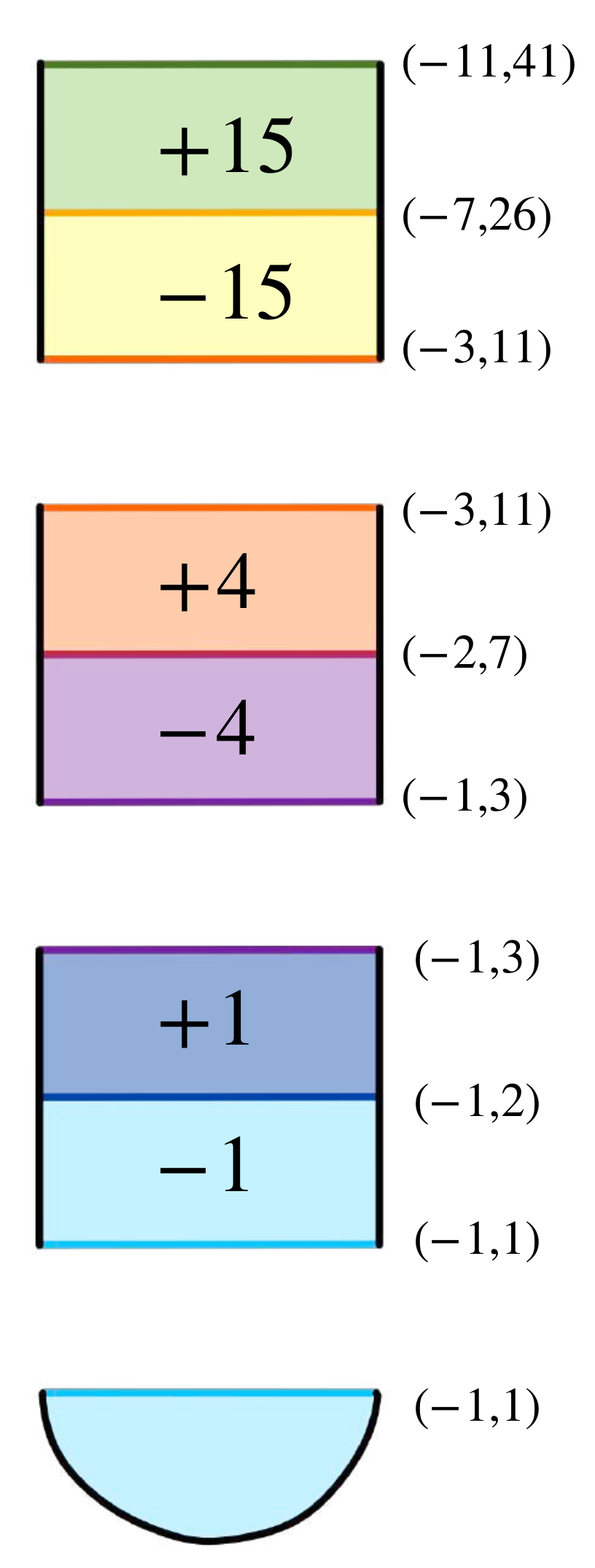}
  \caption{$L(56,15)$}
  \label{blocksL5615}
\end{subfigure}
\begin{subfigure}[t]{.5\textwidth}
  \centering
 \includegraphics[scale=0.5]{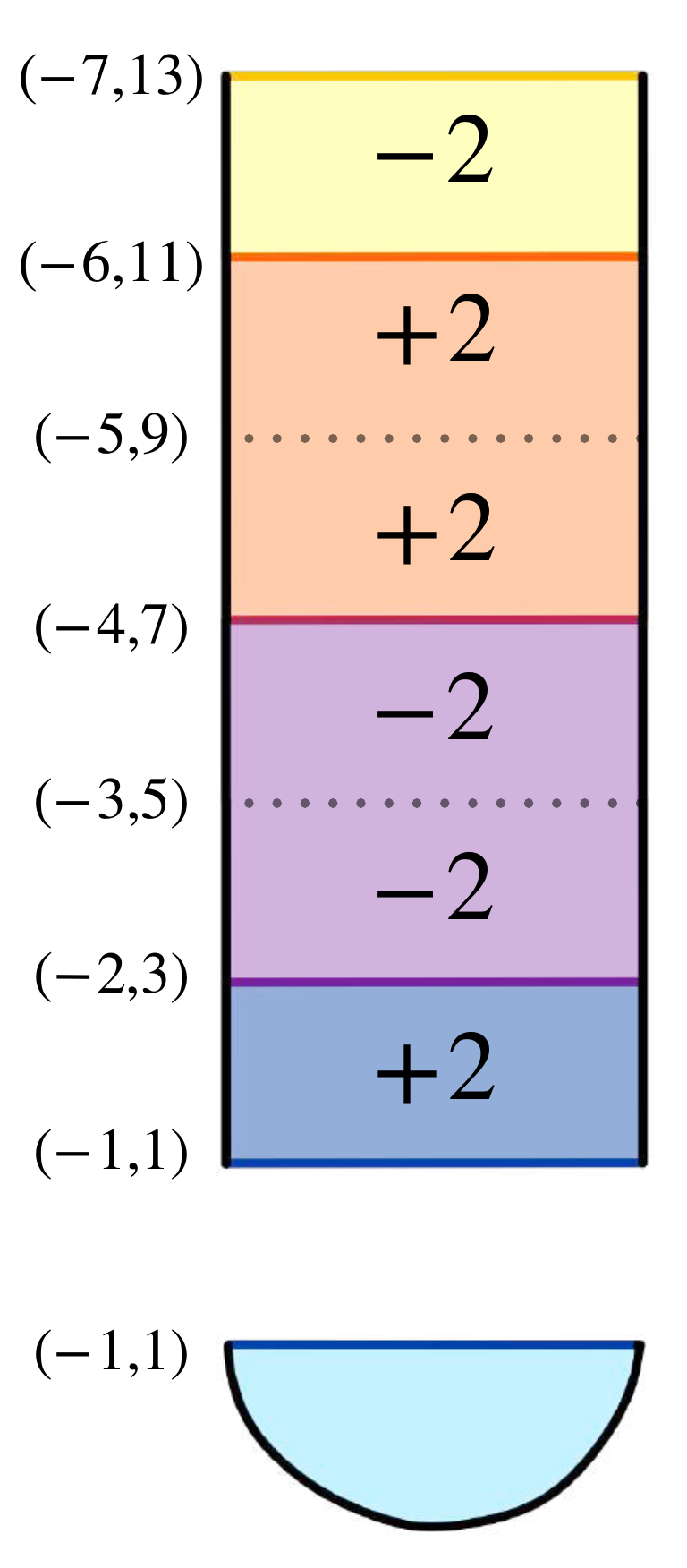}
  \caption{$L(28,15)$}
  \label{blocksL2815}
\end{subfigure}
\caption{Compatible choice of signs for a covering map.}
\label{L5615L2815}
\end{figure}

\begin{figure}[ht!]
\centering
\begin{subfigure}[t]{.45\textwidth}
  \centering
 \includegraphics[scale=0.45]{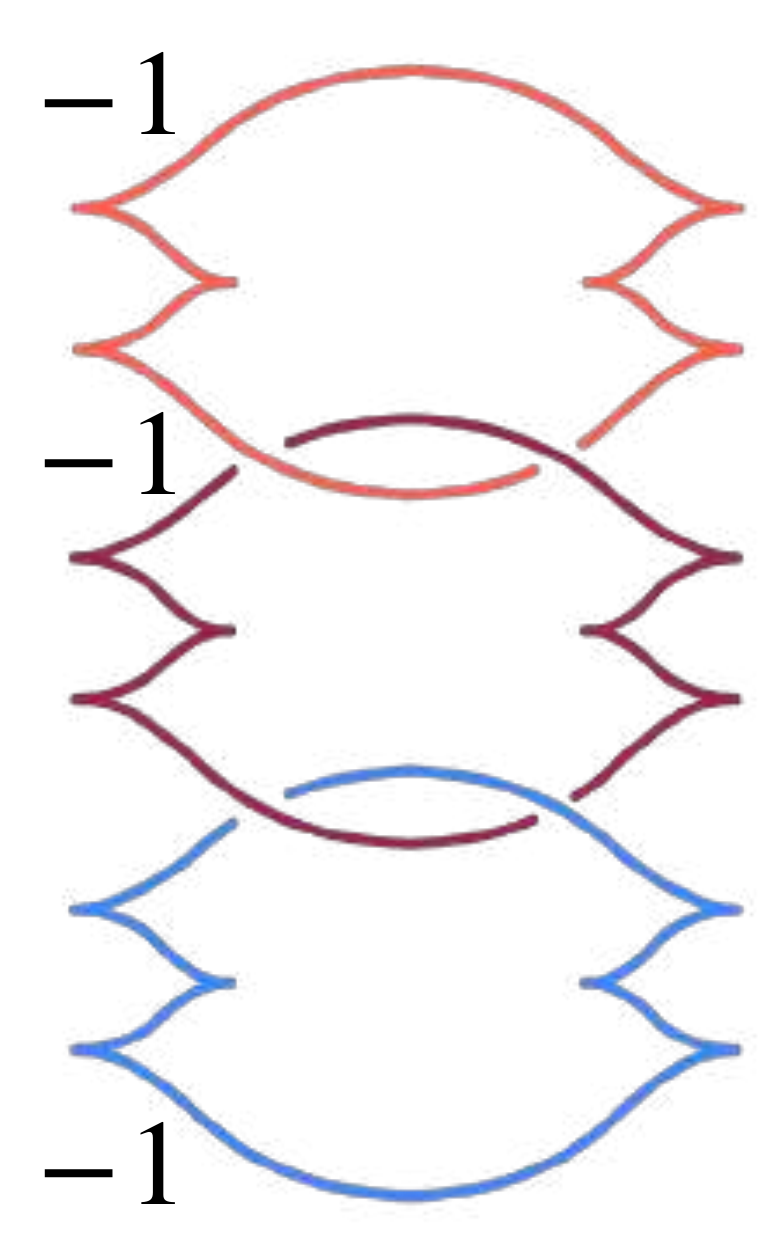}
\caption{$L(56,15)$}
  \label{L5615}
\end{subfigure}
\begin{subfigure}[t]{.45\textwidth}
  \centering
 \includegraphics[scale=0.5]{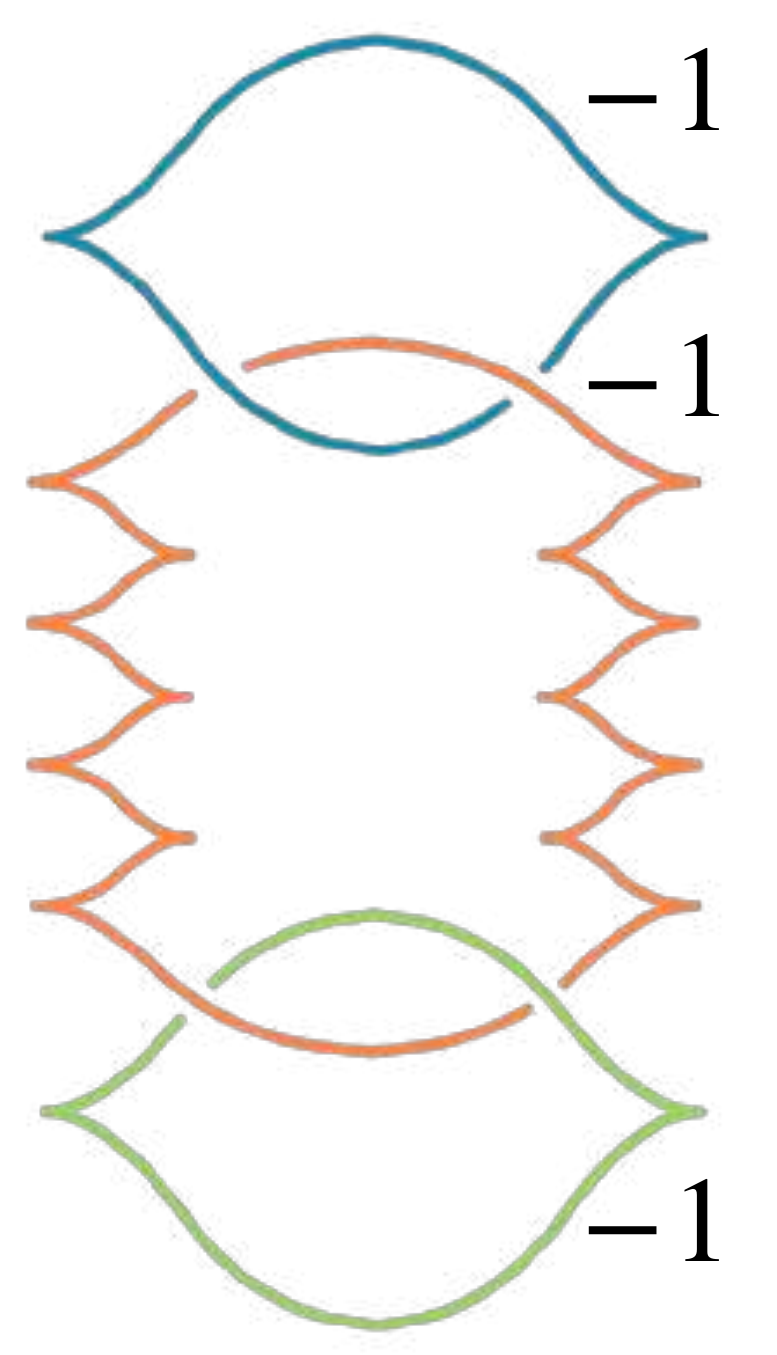}
  \caption{$L(28,15)$}
  \label{L2815}
\end{subfigure}
\caption{Contact surgery producing lens spaces.}
\label{L5615L2815links}
\end{figure}

A compatible case is illustrated for example by Figure \ref{L5615L2815},
which represents the double cover
\[L(28,15)\to L(56,15),\]
where the contact structures on the two lens spaces are specified by Figure \ref{L5615L2815links}.

If we look at the lattice of coverings of $L(56,15)$ we see that this contact structure $\xi$ (Figure \ref{L5615}) lifts to an overtwisted one along some (at least one) covering maps:

\[
\xymatrix{  & (L(7,1),{\color{blue}\xi_{?}})\ar[r]^-{2:1} & (L(14,1),{\color{blue}\xi_{?}}) \ar[r]^-{2:1}& (L(28,15),{\color{blue}\xi_{vot}})\ar[rd]^-{2:1} \\
(S^3,{\color{red}\xi_{ot}}) \ar[ru]^-{7:1}\ar[rd]_-{2:1}& & & &  (L(56,15),{\color{blue}\xi_{vot}}) \\
 & (L(2,1),{\color{red}\xi_{ot}})\ar[r]_-{2:1})\ar[ruu]_-{7:1} & (L(4,3),{\color{red}\xi_{ot}}) \ar[r]_-{2:1}\ar[ruu]_-{7:1}& (L(8,7),{\color{red}\xi_{ot}})\ar[ru]_-{7:1} 
 }\]
 
Therefore we cannot apply directly the criterion of previous section to conclude that the Stein fillings of $(L(56,15),\xi)$ are simply connected. By the fact that lifting $\xi$ to $L(8,7)$ results in an overtwisted structure, we get that the kernel of $i_*$ cannot be contained in $\Z/8\Z$, where
\[i:L(56,15)\hookrightarrow X\]
is the inclusion of the boundary of any Stein filling $X$. We have that
\[\pi_1(X)=\frac{\Z/56\Z}{\ker i_*},\]
so the possibilities are:
\begin{itemize}
\item $\ker i_*=\Z/7\Z$, which gives $\pi_1(X)=\Z/8\Z$,
\item $\ker i_*=\Z/14\Z$, which gives $\pi_1(X)=\Z/4\Z$,
\item $\ker i_*=\Z/28\Z$, which gives $\pi_1(X)=\Z/2\Z$,
\item $\ker i_*=\Z/56\Z$, which gives $\pi_1(X)=1$. The following proposition proves that this is the only possibility.  
\end{itemize}

\begin{prop} \label{pi=1}
Let $X$ be a Stein filling of $(L(56,15),\xi)$, with $\xi$ described by the diagram of Figure \ref{L5615}. Then $\pi_1(X)=1$.
\end{prop}

\begin{prf} 
Consider the Stein filling $X_{\Lambda}$ of $(L(56,15),\xi)$ described by the diagram of Figure \ref{L5615}. We want to compute the $d_3$ invariant of the contact structure on the boundary:
\[d_3(\xi)=\frac{1}{4}(c_1(X_{\Lambda})^2 -3\sigma(X_{\Lambda})-2\chi(X_{\Lambda})).\]
The first Chern class $c_1(X_{\Lambda})$ is zero, because it evaluates as $\rot_i=0$ on the three generators of $H_2(X_{\Lambda})$. Moreover, $\sigma(X_{\Lambda})=-3$ and $\chi(X_{\Lambda})=4$. Therefore
\[d_3(\xi)=\frac{1}{4}.\]
Notice that $c_1(\xi)=0$ because it is the restriction of $c_1(X_{\Lambda})$, which is 0 itself. Being any contact structure on a lens space planar \cite{schonenberger}, we can apply \cite[Corollary 1.5]{c1planar} and conclude that any Stein filling of $(L(56,15),\xi)$ has vanishing $c_1$. 

We want to compute $d_3(\xi)$ using the Stein filling $X$. For what we have just said $c_1(X)^2=0$ and we also have $\sigma(X)=1-\chi(X)$. So:
\[\frac{1}{4}=d_3(\xi)=\frac{1}{4}(c_1(X)^2 -3\sigma(X)-2\chi(X))=\frac{1}{4}(-3+\chi(X)).\]
This tells us that 
\[\fbox{$\chi(X)=4.$}\]

Now we analyze the possibilities for its fundamental group case by case.
\begin{itemize} 
\item[i)] Suppose that $\pi_1(X)=\Z/8\Z.$ Then we pass to the universal covering $\widetilde{X}\to X$, of degree 8, whose boundary is the (connected) degree-8 covering $L(7,1)\to L(56,15)$. By Theorem \ref{pi1b2} we have
\[\chi(X)\leq \frac{1+\length(7/1)}{8}=\frac{1+1}{8}=\frac{1}{4},\]
which is impossible. So $\pi_1(X)\neq \Z/8\Z$.

\item[ii)] If $\pi_1(X)=\Z/4\Z$, we pass to the universal covering and since $\length(14/1)=1$, we get $\chi(X)\leq 1/2$. This is not possible, hence $\pi_1(X)\neq \Z/4\Z$. 

\item[iii)] Again, we take the universal covering $\widetilde{X}\to X$, of degree 2, whose boundary is the (connected) degree-2 covering $L(28,15)\to L(56,15)$. By Theorem \ref{pi1b2}, we have
\[\chi(X)\leq \frac{1+\length(28/15)}{2}=\frac{1+3}{2}=2\]
and hence $\chi(X)\leq 2$, which is not possible. Hence $\pi_1(X)\neq \Z/2\Z$.

\item[iv)] We conclude that any Stein filling of $(L(56,15),\xi)$ is simply connected. (Note that by \cite[Theorem 1.3]{menke} we already know that in fact there is a unique filling obtained by attaching three 2-handles to $B^4$ along the link of Figure \ref{L5615}). 
\end{itemize}
\end{prf}
 
The result proved in Theorem \ref{pi1b2} is that somehow for a Stein filling $X$ of $(L(p,q),\xi_{vo})$ "the bigger $\pi_1(X)$ is, the smaller its Euler characteristic is forced to be". Of course this is in general spoiled by the quantity $l'$, appearing in the statement, which depends on the numbers $p/d$ and $q$ (one should first reduce $q$ modulo $p/d$, in case it were bigger).

On the other hand, if $p$ is small, then by Theorem \ref{thmcovering} we have a bigger chance of finding coverings of $L(p,q)$ which are overtwisted, and hence apply our criterion to bound the cardinality of $\pi_1(X)$. Despite this, there are examples (found by Marco Golla) of fillings with non-trivial fundamental group: let 
\[-\frac{p}{q}=[-4,-2n,-4], \;n>1, \]
and consider the Legendrian representative of the 3-components link associated to this continued fraction expansion where the first and third components have rotation number +2, while the middle one has rotation number 0. The fillings of this virtually overtwisted structure can be completely described by using, for example, the works of \cite{menke} and \cite{mcduff}. In particular, there is a filling which is obtained from a boundary connected sum of two rational homology balls with $\pi_1=\Z/2\Z$ (corresponding to the two $-4$) by attaching a single Weinstein 2-handle (corresponding to the central $-2n$): this handle attachment does not kill the whole $\Z/2\Z\ast \Z/2\Z$, resulting in a non simply-connected filling.

\chapter{Hypersurface singularities and contact structures}\label{boundarymilnor}

It is known that the lens space $L(2n,1)$ supports a virtually overtwisted contact structure arising as the boundary of the Milnor fiber of a complex hypersurface singularity \cite{pichon}. In this chapter we study the problem of realizing other $(L(p,q),\xi)$ in such a way, obtaining a series of necessary conditions for this to happen. The driving question of this chapter is therefore:

\begin{center}
Are the virtually overtwisted structures on lens spaces realizable as the boundary of the Milnor fiber of some complex hypersurface singularity?
\end{center}

Following \cite{nemethi}, we review some terminology of singularity theory. Let $f:(\C^3,0)\to (\C,0)$ be the germ of a complex analytic function with a singularity at the origin, and let 
\[K=f^{-1}(0)\cap S_{\varepsilon}\]
be the \emph{link of the singularity}, where $S_{\varepsilon}$ is the sphere of radius $\varepsilon$ centered at the origin.
Milnor proved in \cite{milnor} that there exists $\varepsilon_0>0$ such that $\forall\, 0<\varepsilon<\varepsilon_0$ the map
\[f/|f|:S_{\varepsilon}\smallsetminus K\to S^1=\{z\in\C:\,|z|=1\}\]
is a smooth fibration. For any such $\varepsilon$ there exists $\delta_{\varepsilon}$ so that $\forall\, 0<\delta<\delta_{\varepsilon}$ the restriction
\[
f:B_{\varepsilon}\cap f^{-1}(\partial D_{\delta}) \to \partial D_{\delta}
\]
is a smooth fibration, whose diffeomorphism type does not depend on $\varepsilon$ and $\delta$ (where $B_{\varepsilon}$ denotes the open ball centered in 0 with radius $\varepsilon$, and $D_{\delta}$ denotes the closed ball centered in 0 with radius $\delta$). This is what we refer to as the \emph{Milnor fibration} of $f$. The fiber 
\[F=F_{\varepsilon,\delta}=B_{\varepsilon}\cap f^{-1}(\delta)\]
is the \emph{Milnor fiber} of $f$ (we omit $\varepsilon$ and $\delta$ from the notation of $F$). If $f$ is the germ of an isolated singularity, then we have a diffeomorphism $\partial F\simeq K$, but in the case of non-isolated singularity $K$ is not smooth (while $F$ and $\partial F$ are always smooth manifolds). The boundary of the Milnor fiber comes with an extra structure, as explained below. 

Contact topology shows up in singularity theory in the following way: the Milnor fiber $F$ of a singularity comes with a Stein structure $J$ which makes it a Stein filling of its boundary $\partial F$ equipped with the contact structure 
\[\xi=T\partial F\cap JT\partial F,\]
which is always tight, see \cite{eliashberg}. Therefore, from a complex germ $f:(\C^3,0)\to (\C,0)$ we obtain a contact 3-manifold $(\partial F,\xi)$ with a Stein filling $(F,J)$. We have a dichotomy: 

\begin{itemize}
\item the singularity is isolated. 
In this case the structure $\xi$ is universally tight, see \cite{lekili}. Another work on this topic is \cite{dynamics}, where the authors show that the link of the hypersurface singularity 
\[z^p+xy=0\]
is $L(p,p-1)$ with its unique tight contact structure (universally tight). We will prove (see Corollary \ref{L(p,p-1)}) that this is the only lens space arising as the link of an isolated hypersurface singularity.

\item The singularity is not isolated. By contrast with previous point, this is the only case where a virtually overtwisted contact structure could arise. A good source of examples is given by the \emph{Hirzebruch singularity}
\[z^2+xy^n=0\]
with $n>1$, for which the boundary of the associated Milnor fiber is $L(2n,1)$. This type of singularity is studied in \cite[Section 6]{pichon}.
\end{itemize}

We study those lens spaces $L(p,q)$ with a tight contact structure $\xi$ arising as the boundary of the Milnor fiber of a hypersurface singularity $f:(\C^3,0)\to (\C,0)$. Theorem \ref{mainthm} gives a partial answer to a question raised by \cite[open problems 24.4.2]{nemethi}.

\begin{thm}\label{mainthm}
Let $\xi_{vo}$ be a virtually overtwisted structure on $L(p,q)$. If we are in one of the cases below, then $(L(p,q),\xi_{vo})$ is not the boundary of the Milnor fiber of any complex hypersurface singularity:
\begin{itemize}
	\item[a)] $p/q=[a_1,a_2,\ldots,a_n]$ and $a_i$ is odd for some $i$;
	\item[b)] $p/q=[2x_1,2x_2]$;
	\item[c)] $p/q=[2x_1,2x_2,\ldots,2x_n]$, with $x_i>1$ for every $i$ ($n\geq 3$) and either:
\begin{itemize}
		\item[i)] $q^2\not\equiv 1 \pmod{p}$ or
		\item[ii)] $q^2\equiv 1 \pmod{p}$ and $n$ is even.
\end{itemize}
\end{itemize}
\end{thm}

The first step in proving this theorem is to characterize those contact structures $\xi$ which can appear in the context of hypersurface singularity: the fact that $c_1(\xi)$ vanishes imposes certain conditions on the coefficients of the continued fraction expansion of $p/q$, which allow us to prove part $(a)$ of Theorem \ref{mainthm}. Section \ref{secvan} deals with these numerical restrictions, adapted to the language of contact geometry. To prove parts $(b)$ and $(c)$ we need to look closely at the topology of the Milnor fibration and analyze its monodromy. In order to derive our statements we study the integral orthogonal group of the intersection form of the Milnor fiber, imposing a further restriction coming from a theorem of A'Campo \cite{acampo}. This is explained in Section \ref{mainproof}. The chapter ends with an open question which focuses on the limits of Theorem \ref{mainthm}: where are these techniques failing?

\section{Vanishing of the rotation numbers}\label{secvan}

The goal of this section is to prove Theorem \ref{thmc1=0}, which is a special case of \cite[Corollary 1.5]{c1planar}. We prove it using elementary techniques that do not involve the Ozsváth-Szabó contact invariant (defined in \cite{ozsvath2005heegaard}). Theorem \ref{thmc1=0} will be the starting point in the proof of Theorem \ref{mainthm}. 

\begin{thm}\label{thmc1=0}
Let $L$ be a linear chain of Legendrian unknots in the standard contact $S^3$ and let $(L(p,q),\xi)$ be the contact 3-manifold obtained by Legendrian surgery on $L$. If $c_1(\xi)=0$, then $\rot (L_i)=0$ for every component $L_i$ of the link $L$.
\end{thm}

To prove Theorem \ref{thmc1=0} we need some notation that makes the computation easier, and a few more results. Let 
\[\frac{p}{q}=[a_1,a_2,\ldots, a_n],\]
and let $Q$ be the intersection form of the plumbed 4-manifold associated to the linear graph
\begin{center}
\begin{tikzpicture}
        \node[shape=circle,fill=black,inner sep=1.5pt,label=$-a_1$] (1)                  {};
        \node[shape=circle,fill=black,inner sep=1.5pt,label=$-a_2$] (2) [right=of 1] {}
        edge [-]               (1);
        \node[shape=circle,fill=black,inner sep=1.5pt,label=$-a_3$] (3) [right=of 2] {}
        edge [-]               (2);
        \node[shape=circle,fill=black,inner sep=1.5pt,label=$-a_{n-1}$] (4) [right=of 3] {} ;
                 \node at ($(3)!.5!(4)$) {\ldots};

        \node[shape=circle,fill=black,inner sep=1.5pt,label=$-a_n$] (5) [right=of 4] {}
        edge [-]               (4);
\end{tikzpicture}
\end{center}
written in the ordered basis given by the meridian of each attaching circle. We know that
\[H^2(L(p,q))\ni c_1(\xi)=0 \Leftrightarrow \PD(c_1(\xi))=0\in H_1(L(p,q)),\]
and we express
\[\PD(c_1(\xi))=\sum_{i=1}^nr_i\mu_i,\]
where $r_i$ and $\mu_i$ are respectively the rotation number $\rot(L_i)$ and the meridian of the $i^{th}$ component of the Legendrian link describing $(L(p,q),\xi)$ (compare with \cite[Proposition 8.2.4]{ozbagci}). By looking at the linear plumbing graph, we can find relations among the $\mu_i$'s and get a set of equations coming from $Q$, which hold in $H_1(L(p,q))$:
\[
\begin{cases}
-a_1\mu_1+\mu_2=0	\\
\mu_1-a_2\mu_2+\mu_3=0 \\
\vdots \\
\mu_j-a_{j+1}\mu_{j+1}+\mu_{j+2}=0 \\
\vdots \\
\mu_{n-1}-a_n\mu_n=0.
\end{cases}
\]
We choose $\mu_1\in H_1(L(p,q))$ as a generator and we lift it from $\Z/p\Z\simeq H_1(L(p,q))$ to $\Z$. Hence, from now on, we identify $\mu_1$ with $1\in \Z$ and all the other $\mu_i$'s with integer numbers according to the recursive expressions:
\[
\begin{cases}
\mu_1=1\\
\mu_2=a_1\\
\mu_i=a_{i-1}\mu_{i-1}-\mu_{i-2}.
\end{cases}
\]
Thanks to Lemma \ref{lemma0}, we do not need to make computations in $\Z/p\Z$, hence we can use previous relations, which hold over $\Z$. Then $c_1(\xi)$ is 0 exactly when 
\begin{equation}\label{chernrotations}
\sum_{i=1}^n r_i\mu_i\equiv 0\pmod{p}.
\end{equation}
Define
\[
\begin{cases}
\D[-1]=0;\\
\D[0]=1;\\
\D[i]=-a_i\D[i-1]-\D[i-2]
\end{cases}
\]
and note that
\begin{itemize}
\item $\det(Q)=\D[n]=\pm p$;
\item $\sign(\D[i])=(-1)^i$, hence $\D[i]=(-1)^i|\D[i]|$;
\item $\mu_i> \mu_{i-1}$;
\item $|r_i|\leq a_i-2$ and $r_i\equiv a_i \pmod{2}$;
\item $\D[i]=(-1)^i\mu_{i+1}$ (proved by induction), therefore $|\D[i]|>|\D[i-1]|$.
\end{itemize}

\begin{lem}\label{lemma0}
Equation \eqref{chernrotations} is satisfied in $\Z/p\Z$ if and only if it is satisfied in $\Z$.
\end{lem}

\begin{prf}
We prove by induction on $n$ that 
\begin{equation}\label{inequality}
p>\abs[\Big]{\sum_{i=1}^n r_i\mu_i}.
\end{equation}
This will tell that Equation \eqref{chernrotations} can only be satisfied with 
\[\abs[\Big]{\sum_{i=1}^n r_i\mu_i}=0.\]
If $n=1$ then the two sides of Inequality \eqref{inequality} are respectively $|a_1|$ and $|r_1|$, so it is true. 
The first interesting case is then $n=2$:
\[a_2a_1-1>r_1+r_2a_1?\] 
We have:
\begin{align*}
r_1+r_2a_2 < &(a_1-2)+(a_2-2)a_1\\
= & a_1a_2-a_1-2\\
< & a_1a_2-1.\qquad \checkmark
\end{align*}
The general case now:
\begin{align*}
\abs[\Big]{\sum_{i=1}^nr_i\mu_i}= & \abs[\Big]{\sum_{i=1}^nr_i(-1)^{i-1}\D[i-1]}\\
\leq &\left(\abs[\Big]{\sum_{i=1}^{n-1}r_i(-1)^{i-1}\D[i-1]}\right)+|r_n\D[n-1]| \\
<&|\D[n-1]|+|r_n\D[n-1]|&& \mbox{(induction)} \\
\leq &|\D[n-1]|+|(a_n-2)\D[n-1]|.
\end{align*}
There are two possibilities for the right-hand side, according to the parity of $n$.
\begin{itemize}
\item[1)] $n$ is even (hence $\D[n]>0$):
\begin{align*}
|\D[n-1]|+|(a_n-2)\D[n-1] |=& -\D[n-1]-a_n\D[n-1]+2\D[n-1]\\
=&-a_n\D[n-1]+\D[n-1]\\
<& -a_n\D[n-1]-\D[n-2] \mbox{ (because $n$ is even)}\\
=&\D[n]\\
=&{|\det(Q)|}=p.
\end{align*}
\item[2)] $n$ is odd (hence $\D[n]<0$):
\begin{align*}
|\D[n-1]|+|(a_n-2)\D[n-1] |=& a_n\D[n-1]-\D[n-1]\\
<& a_n\D[n-1]+\D[n-2] \mbox{ (because $n$ is odd)}\\
=&-\D[n]\\
=&{|\det(Q)|}=p.
\end{align*}
\end{itemize}
\end{prf}

\begin{lem}\label{lemma1} 
\[a_n\mu_n-a_{n-1}\mu_{n-1}>0.\]
\end{lem}

\begin{prf} We prove it by induction on $n$. If $n=2$ then the formula is just
\[a_2\mu_2-a_1\mu_1=a_2a_1-a_1>0.\qquad\checkmark\]
In general, assuming the result true for $n-1$, we have
\begin{align*}
a_n\mu_n-a_{n-1}\mu_{n-1}=& a_n(a_{n-1}\mu_{n-1}-\mu_{n-2})-a_{n-1}\mu_{n-1}\\
=& a_{n-1}\mu_{n-1}(a_n-1)-a_n\mu_{n-2}\\
>& a_{n-2}\mu_{n-2}(a_n-1)-a_n\mu_{n-2}&& \mbox{(induction)}\\
=& (a_{n-2}a_n-a_{n-2}-a_n)\mu_{n-2}\\
\geq & 0.
\end{align*}
\end{prf}

\begin{lem}\label{lemma2} If $r_n\neq 0$, then we have
\[|r_n\mu_n|-|r_{n-1}\mu_{n-1}|>0.\]
\end{lem}

\begin{prf}
\begin{align*}
|r_n\mu_n|-|r_{n-1}\mu_{n-1}|>& |\mu_n|-(a_{n-1}-2)\mu_{n-1}\\
=& a_{n-1}\mu_{n-1}-\mu_{n-2}-a_{n-1}\mu_{n-1}+2\mu_{n-1}\\
=& 2\mu_{n-1}-\mu_{n-2}\\
> & 0.
\end{align*}
\end{prf}

\begin{lem}\label{lemma3}
If $r_n\neq 0$, then we have
\[|r_n\mu_n|-\abs[\Big]{\sum_{i=1}^{n-1}r_i\mu_i}>0.\]
\end{lem}

\begin{prf} We prove it by induction on $n$. If $n=2$, the formula is the same as the one of Lemma \ref{lemma2}, and so we know it holds. Now we do the general case. Assume the inequality holds for $n-1$ and let $j\leq n-2$ be the biggest integer such that $r_j\neq 0$ (note that if $r_i=0,\;\forall i\leq n-2$, then by Lemma \ref{lemma2} we would be done after the second line in the following computation). 
We have 
\begin{align*}
|r_n\mu_n|-\abs[\Big]{\sum_{i=1}^{n-1}r_i\mu_i} >& |\mu_n|-\abs[\Big]{\sum_{i=1}^{j-1}r_i\mu_i}-|r_{n-1}\mu_{n-1}| \\
>& \mu_n-|r_j\mu_j|-|r_{n-1}\mu_{n-1}|\hspace{1.7cm}  \mbox{(induction)} \\
\geq & \mu_n-(a_j-2)\mu_j-(a_{n-1}-2)\mu_{n-1}\\
= & a_{n-1}\mu_{n-1}-\mu_{n-2}-(a_j-2)\mu_j-(a_{n-1}-2)\mu_{n-1} \\
= & 2\mu_{n-1}-\mu_{n-2}-(a_j-2)\mu_j	\\
> & \mu_{n-1}-(a_j-2)\mu_j\\
= & \mu_{n-1}+\mu_j+\mu_j-a_j\mu_j \\
>&\mu_{j+1}+\mu_{j-1}+\mu_{j-1}-a_j\mu_j \qquad\qquad (n-2\geq j\Rightarrow n-1\geq j+1) \\
=&\mu_{j-1} \hspace{4.8cm} (\mu_{j+1}+\mu_{j-1}=a_{j}\mu_{j})\\
>&0.
\end{align*}
\end{prf}

\begin{lem} \label{lemma4}
\[\sum_{i=1}^{n}r_i\mu_i=0\Longrightarrow r_i=0,\;\forall i.\]
\end{lem}

\begin{prf}
\[\sum_{i=1}^{n}r_i\mu_i=0\Longrightarrow r_n\mu_n=-\sum_{i=1}^{n-1}r_i\mu_i\Longrightarrow |r_n\mu_n|=\abs[\Big]{\sum_{i=1}^{n-1}r_i\mu_i}.
\]
But if $r_n\neq 0$, then we should have a strict inequality by Lemma \ref{lemma3}, hence $r_n=0$. By applying this repeatedly we get to
\[r_n=r_{n-1}=\ldots=r_1=0.\]
\end{prf}

\noindent We can finally give the following:

\begin{prf}[of Theorem \ref{thmc1=0}]
By combining Lemmas \ref{lemma0} and \ref{lemma4}, we have that
\[c_1(\xi)=0\Longleftrightarrow \sum_{i=1}^n r_i\mu_i\equiv 0\pmod{p}\; \Longleftrightarrow \;\sum_{i=1}^n r_i\mu_i=0\Longleftrightarrow r_i=0\,\forall i,\]
where $r_i=\rot (L_i)$.
\end{prf}

\section{Proof of Theorem \ref{mainthm} }\label{mainproof}

\noindent In the universally tight case for lens spaces we know (see \cite{nemethippp} and \cite{bhupal}) that all the fillings come from algebraic geometry, and Choi and Park show in the article \cite{park} that the theory of surface singularities describes all the fillings of small Seifert 3-manifold equipped with the  canonical contact structure.

\vspace{0.5cm}

\noindent \textbf{Key fact:} suppose there is a polynomial function $f:(\C^3,0)\to (\C,0)$ such that 
\[\partial(F,J)=(L(p,q),\xi),\]
where $(F,J)$ is the Milnor fiber of $f$ and $\xi$ is the contact structure on the boundary induced by complex tangencies, as explained above. The Legendrian link representation of $(L(p,q),\xi)$ describes more than a contact 3-manifold: the components $L_i$ of the link can be thought of as the attaching circles of the 2-handles of the Stein domain $(F,J)$, as Corollary \ref{uniqueness} explains. The first Chern class $c_1(F,J)\in H^2(F;\Z)$ evaluates on each 2-handle as the correspondent rotation number, and $c_1(\xi)$ is the restriction of $c_1(F,J)$. But the tangent bundle of $F$ is stably trivial (since the normal bundle of a complex hypersurface of $\C^3$ is trivial), hence $c_1(F,J)=0$, and also $c_1(\xi)=0$ on the boundary $L(p,q)$. 

\begin{prf}[of Theorem \ref{mainthm}a]
If $(L(p,q),\xi_{vo})$ is the boundary of the Milnor fiber of a complex hypersurface singularity, then $c_1(\xi_{vo})=0$ and, by Theorem \ref{thmc1=0}, all $\rot(L_i)$ are zero. \\
\noindent Let $p/q=[a_1,a_2,\ldots,a_n]$ and remember that 
\[a_i\equiv \rot(L_i)\pmod{2}.\] 
Since all the rotation numbers are zero, the conclusion follows.
\end{prf}

Theorem \ref{thmc1=0} implies also the following corollary, which will be used later.

\begin{cor}\label{uniqueness}
If $(L(p,q),\xi)$ is the boundary of the Milnor fiber of a hypersurface singularity, then it has a unique Stein filling, which is the Milnor fiber itself.
\end{cor}

\begin{prf}
The fact that $\rot(L_i)=0$ for every $i$ implies, by \cite[Theorem 1.3]{menke}, that from the chain of Legendrian unknots producing $(L(p,q),\xi)$ we can forget about those components with $\tb(L_i)\neq -1$ and look for Stein fillings of the 3-manifold $Y$ which is left. Then, all the Stein fillings of $(L(p,q),\xi)$ will be uniquely obtained by attaching the 2-handles (corresponding to the forgotten components) to the Stein fillings of $Y$. From the link diagram of $(L(p,q),\xi)$ we see that $Y$ is a connected sum of $(L(n_j,n_j-1),\xi_{st})$, where each of the prime factor corresponds to a string of $-2$ in the expansion of $p/q$. By \cite[Theorem 16.9]{steinandback}, $Y$ admits a unique Stein filling, because each factor $(L(n_j,n_j-1),\xi_{st})$ does (by the work of Lisca \cite{lisca}). This concludes the proof.
\end{prf}

\vspace{0.5cm}

\noindent Another consequence of Theorem \ref{thmc1=0} is:

\begin{cor}\label{L(p,p-1)}
Let $(L(p,q),\xi)$ be a lens space with a contact structure arising as the link of an isolated hypersurface singularity. Then $q=p-1$.
\end{cor}

\begin{prf}
From Theorem \ref{thmc1=0} we have that all the rotation numbers are zero and \cite[Theorem 2.1]{lekili} says that $\xi$ is universally tight. By the work of \cite{honda} we know that a universally tight structure on a lens space is the result of contact $(-1)$-surgery on a link where all the stabilizations appear on the same side, i.e. when the (absolute values of the) rotation numbers are maximal. Therefore, every Legendrian knot must have Thurston-Bennequin number equal to $-1$:
\[-\frac{p}{q}=[-2,-2,\ldots,-2]\qquad\Rightarrow\qquad q=p-1.\]
\end{prf}

\subsection*{Proof of Theorem \ref{mainthm}b}

Corollary \ref{uniqueness} says that if $(L(p,q),\xi_{vo})$ arises as $\partial (F,J)$, then the Stein filling $F$ is uniquely determined: topologically it is given by the plumbing of spheres according to the expansion of $p/q$. The monodromy $\varphi$ of the Milnor fibration induces, in cohomology, a homomorphism
\[\varphi^*:H^*(F;\Z)\to H^*(F;\Z)\]
such that the alternating sum of the traces is zero, by \cite[Theorem 1]{acampo}:
\[\trace(\varphi^*_0)-\trace(\varphi^*_1)+\trace(\varphi^*_2)-\trace(\varphi^*_3)+\trace(\varphi^*_4)=0.\]
In our case, the Stein fillings of those lens spaces with $c_1(\xi)=0$ are simply connected, hence $\varphi^*_1=\varphi^*_3=0$. Moreover, $\varphi^*_0:\Z\to\Z$ is the identity and $\varphi^*_4=0$. Hence, previous equation simply reads as:
\begin{equation}\label{trace}
1+\trace(\varphi^*_2)=0.
\end{equation}
The intersection form of $F$ must be preserved by the homomorphism $\varphi^*_2$ and we are therefore led to study its isometry group. If in this group there is no element whose trace is $-1$, then Equality \eqref{trace} cannot be satisfied, and we conclude that the polynomial function $f:(\C^3,0)\to (\C,0)$ with $(\partial F,\xi)=(L(p,q),\xi_{vo})$ does not exist.

\begin{prf}[of Theorem \ref{mainthm}b]
Our goal is to prove that there is no $f:(\C^3,0)\to (\C,0)$ with non-isolated singularity at the origin, whose Milnor fiber has boundary $L(p,q)$. Note that 
\[x_1x_2>1,\]
otherwise the induced contact structure is universally tight. 

Assume by contradiction that such $f$ exists. Then, by Corollary \ref{uniqueness}, we know that the Milnor fiber $F$ has negative-definite intersection form isomorphic to
\[-M=
\begin{bmatrix}
   -2x_1 & 1 \\
   1 & -2x_2 \\
\end{bmatrix}.
\]
By Equality \eqref{trace}, we must have $\trace(\varphi^*_2)=-1$. The morphism $\varphi^*_2:H^2(F;\Z)\to H^2(F;\Z)$ is induced by a diffeomorphism which preserves the intersection form and therefore it is represented by an integral matrix $A$ with
\[\begin{cases}
{|\det(A)|}=1 \\
\trace(A)=-1
\end{cases}
\]
and such that 
\[A(-M)A^T=-M.\]
We show now that such matrix cannot exist. We change sign to work with a positive definite matrix:
\[AMA^T=M\Longrightarrow
 \begin{bmatrix}
   a_1 & a_2 \\
   a_3 & a_4 \\
\end{bmatrix}
\begin{bmatrix}
   2x_1 & -1 \\
   -1 & 2x_2 \\
\end{bmatrix}
\begin{bmatrix}
   a_1 & a_3 \\
   a_2 & a_4 \\
\end{bmatrix} = 
\begin{bmatrix}
   2x_1 & -1 \\
   -1 & 2x_2 \\
\end{bmatrix}. \]
We get equations:
\[\begin{cases}
2x_1a_1^2-2a_1a_2+2x_2a_2^2=2x_1 \\
2x_1a_3^2-2a_3a_4+2x_2a_4^2=2x_2
\end{cases}
\]
that can be rewritten as

\begin{subequations}
\begin{empheq}[left=\empheqlbrace]{align}
  & (2x_1-1)a_1^2+(a_1-a_2)^2+(2x_2-1)a_2^2=2x_1 \label{star1}
  \\
  & (2x_1-1)a_3^2+(a_3-a_4)^2+(2x_2-1)a_4^2=2x_2. \label{star2} 
\end{empheq}
\end{subequations}

\noindent From Equations \eqref{star1} and \eqref{star2} it follows that $a_1^2\leq 1$ and $a_4^2\leq 1$. But since $\trace(A)=a_1+a_4=-1$, we have that either
\[ \begin{cases}
a_1=0 \\
a_4=-1
\end{cases}
\qquad \mbox{or}\qquad 
\begin{cases}
a_1=-1 \\
a_4=0.
\end{cases}\]
We do the first case $(a_1=0,\,a_4=-1)$, the other one is the same. From $AMA^T=M$ we also get 
\[
 \begin{bmatrix}
   0 & a_2
\end{bmatrix}
\begin{bmatrix}
   2x_1 & -1 \\
   -1 & 2x_2 \\
\end{bmatrix}
\begin{bmatrix}
   a_3 \\
   -1  \\
\end{bmatrix} = -1, \]
which gives the equation $a_2(a_3+2x_2)=1$. Hence $a_2=a_3+2x_2=1$ (or both $-1$, but the conclusion is the same).
From Equation \eqref{star1} we have $x_2a_2^2=x_1$, which gives $x_2=x_1$. Then
\[\pm 1=\det(A)=\det 
\begin{bmatrix}
   0 & \pm 1 \\
   a_3 & -1  \\
\end{bmatrix} 
\Longrightarrow a_3=\pm 1.\]
We are in the case where $a_3+2x_2=1$, so either $x_1=x_2=0$ or $x_1=x_2=1$, which are both contradicting the condition $x_1x_2>1$.
\end{prf}

\subsection*{Proof of Theorem \ref{mainthm}c}

Remember that in order to have $c_1(\xi)=0$ we need all the rotation numbers to be zero and, in particular, all the coefficients in the expansion to be even. Let $-M$ be the negative definite intersection lattice associated to the linear plumbing of spheres
\begin{center}
\begin{tikzpicture}
        \node[shape=circle,fill=black,inner sep=1.5pt,label=$-2x_1$] (1)                  {};
        \node[shape=circle,fill=black,inner sep=1.5pt,label=$-2x_2$] (2) [right=of 1] {}
        edge [-]               (1);
       \node[shape=circle,fill=black,inner sep=1.5pt,label=$-2x_3$] (3) [right=of 2] {}
        edge [-]               (2);  
        \node[shape=circle,fill=black,inner sep=1.5pt,label=$-2x_{n-1}$] (4) [right=of 3] {} ;
                         \node at ($(3)!.5!(4)$) {\ldots};

        \node[shape=circle,fill=black,inner sep=1.5pt,label=$-2x_n$] (5) [right=of 4] {}
        edge [-]               (4);
\end{tikzpicture}
\end{center}

We are looking for a matrix $A$ representing the monodromy of a Milnor fibration on the second cohomology group, that respect the intersection form of the Milnor fiber (i.e. $A(-M)A^T=M$) and whose trace is $-1$.

What we need to understand is the integral orthogonal group $O_{\Z}(-M)$ of the negative definite lattice $(\Z^n,-M)$, which is isomorphic to $O_{\Z}(M)$. The latter is studied in the article \cite{gerstein}, where the following theorem is proved:

\begin{thm*}[\cite{gerstein}] 
Let $M$ be the integer matrix
\[
\begin{bmatrix}
    2x_1       & -1  \\
 -1       & \ddots & \ddots \\
 & \ddots & \ddots & \ddots \\
  & & \ddots &\ddots &-1\\
  & & & -1 & 2x_n\\
\end{bmatrix}
\]
with $n\geq 2$ and $2x_i\geq 3\;\forall i$. 
\begin{itemize}
\item[i)] If $x_i\neq x_{n+1-i}$ for some $i$, then $O_{\Z}(M)=\{\pm \id\}$.
\item[ii)] If $x_i=x_{n+1-i}$ for every $i=1,\ldots, n$, then $O_{\Z}(M)=\{\pm \id,\pm \rho\}$, where $\rho$ is the isometry that inverts the order of a basis.
\end{itemize}
\end{thm*}

Therefore we can derive:

\begin{prf}[of Theorem \ref{mainthm}c]
The condition $q^2\not\equiv 1 \pmod{p}$ can be rephrased in terms of the coefficients of the expansion by saying that $x_i\neq x_{n+1-i}$ for some $i$, see \cite[Appendix]{orlik}. On the other hand, $x_i=x_{n+1-i}$ for all $i$ if and only if $q^2\equiv 1 \pmod{p}$.

In the first case, the theorem of Gerstein quoted above tells us that if there is an integer matrix $A$ with $AMA^T=M$, then $A=\pm	\id$. Since this does not have trace $-1$, there cannot be a hypersurface singularity whose Milnor fiber has boundary $(L(p,q),\xi_{vo})$, otherwise we would have a contradiction with A'Campo's Equality \eqref{trace}.

In the second case, again A'Campo's formula cannot be satisfied because a hypersurface singularity would produce a Milnor fibration with monodromy $A\in\{\pm \id,\pm \rho\}$:
\[\rho=
\begin{bmatrix}
     & & & 1  \\
  &   & 1  \\
  & \reflectbox{$\ddots$} \\
   1\\
\end{bmatrix}
\]
and, if $n$ is even, then $\trace(\rho)=0\neq -1$.
\end{prf}

\begin{que} There are cases which are not covered by Theorem \ref{mainthm}: every time that in the continued fraction expansion of $-p/q$ there is a $-2$, the orthogonal group $O_{\Z}(M)$ is harder to understand. Nevertheless, in the easier case when $q^2\equiv 1\pmod{p}$ we have a complete description of $O_{\Z}(M)$, but inside this group there is a matrix with trace $-1$ if the length of the expansion is odd. A simple case is for example 
\[-\frac{p}{q}=-\frac{12}{7}=[-2,-4,-2].\] 
Our techniques indeed do not exclude that the isometry 
\[-\rho=
\begin{bmatrix}
    0 & 0 &-1  \\
  0 &-1  & 0  \\
   -1 &0  & 0  \\
\end{bmatrix}
\]
is the morphism induced by the monodromy of the Milnor fibration of a certain non-isolated hypersurface singularity producing $L(12,7)$ with the virtually overtwisted contact structure represented by a Legendrian link with the rotation numbers of the three components equal to zero. How can we deal with cases like this and exclude the existence of such a monodromy? Further works will hopefully clarify this problem and either find the polynomial function or rule out this possibility as well.
\end{que}

\chapter{Artin presentations and contact geometry}\label{artinandcontact}

We present in this final chapter an interesting approach to 3-manifolds, through the theory of Artin presentations. As we will see, to the datum of such a presentation it is possible to associate a (planar) open book decomposition of a 3-manifold and therefore a contact structure on it.  We highlight the main results of the theory and derive some consequences related to contact geometry, see Section \ref{artingeom}. What follows is the result of several discussions held in Budapest between the author and Fabio Gironella.

\vspace{0.7cm}

A group presentation is called \emph{Artin presentation} (of length $n$) if it is of the form
\[p_n=\langle x_1,\ldots, x_n\;|\;r_1,\ldots, r_n\rangle,\]
for some $n$, with the requirement that, in the free group on $n$ generators $\{x_1,\ldots,x_n\}$, the following equality holds:
\[\displaystyle\prod_{i=1}^{n} x_i = \displaystyle\prod_{i=1}^{n} r_i^{-1}x_ir_i.\]
Such group presentations came to the attention of topologists when it was proved that a group $G$ admits an Artin presentation if and only if $G$ is the fundamental group of a closed oriented 3-manifold, see \cite{winkelnkemper}. More specifically, the datum of an Artin presentation $p_n$ corresponds to a pair $(\Sigma_{n+1},\varphi)$, where $\Sigma_{n+1}$ is a planar surface with $n+1$ boundary components and $\varphi$ is an element in the mapping class group $\Gamma(\Sigma_{n+1})$, where homeomorphisms are the identity on the boundary and isotopies are relative to the boundary. Vice versa, given such a pair, it is possible to construct an Artin presentation. In what follows, we highlight the steps in the construction of this bijection.

\section{From \texorpdfstring{$(\Sigma,\varphi)$}{(S,f)} to \texorpdfstring{$p$}{p}}

\begin{figure}[h!]
\centering
\begin{subfigure}[t]{.5\textwidth}
  \centering
  \includegraphics[scale=0.5]{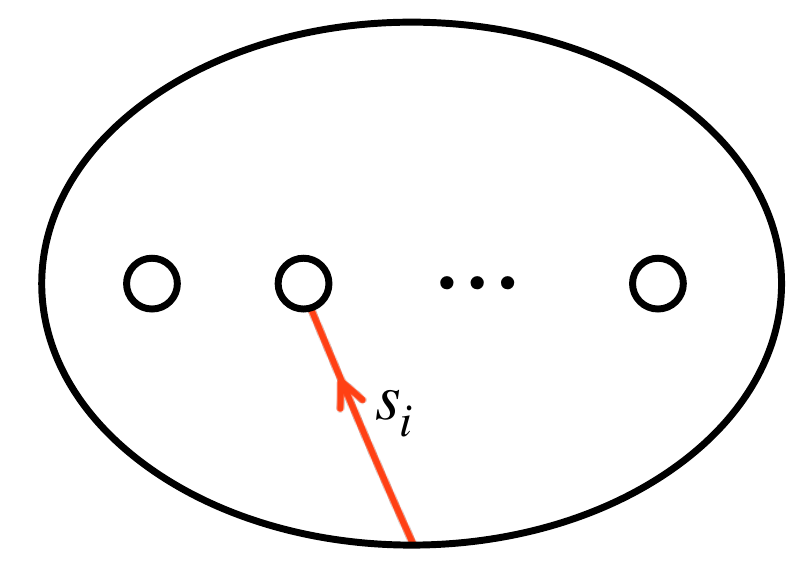}
  \caption{The arc $s_i$.}
  \label{arcSi}
\end{subfigure}%
\begin{subfigure}[t]{.5\textwidth}
  \centering
 \includegraphics[scale=0.5]{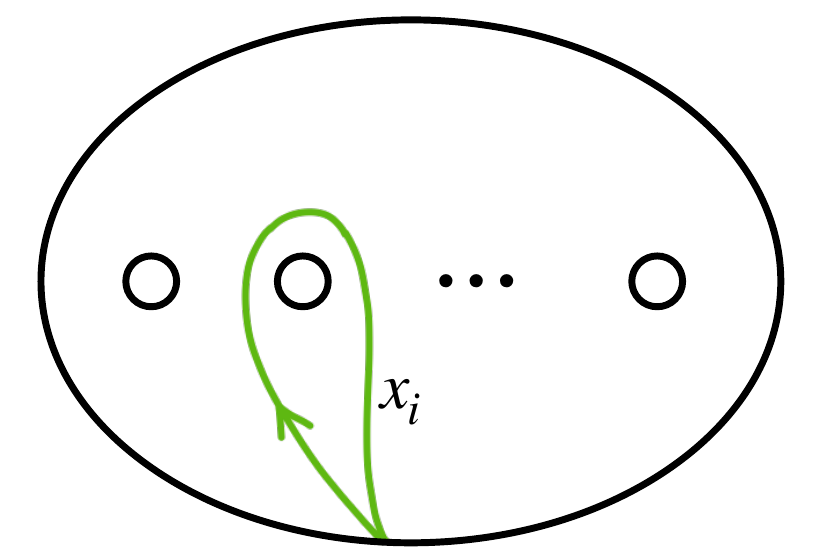}
  \caption{The loop $x_i$.}
  \label{loopXi}
\end{subfigure}
\caption{ }
\end{figure}

Let $s_i$ be the oriented arc starting on the outer component with endpoint on the $i^{th}$ inner component, as in Figure \ref{arcSi}. Denote by $x_i$ the oriented loops as in Figure \ref{loopXi}, which generate the fundamental group $\pi_1(\Sigma)\simeq F_n$ for $i=1,\ldots,n$, and are based at $x_i(0)$. From now on, we will consider the fundamental group of $\Sigma$ based at $x=x_i(0)$.

Now define the $i^{th}$ relation $r_i(\varphi)$ by looking at the action of the diffeomorphism $\varphi$ on $\pi_1(\S,x)$ (notice that this is possible because $\varphi$ is the identity near the boundary of $\Sigma$):
\[r_i=r_i(\varphi)=s_i\ast \varphi(s_i^{-1}),\]
written as a word in the alphabet $\{x_1^{\pm 1},\ldots, x_n^{\pm 1}\}$, where $\ast$ denotes the concatenation of paths (read from left to right). This way we get a collection of words $r_1,\ldots,r_n$ which produces the Artin presentation $\langle x_1,\ldots, x_n\;|\;r_1,\ldots, r_n\rangle$, compare with \cite{winkelnkemper}. We denote the Artin presentation obtained this way as $p=p(\varphi)$.

The set of Artin presentations $\mathcal{P}_n$ of length $n$ can be endowed with the following operation, see \cite[page 227]{winkelnkemper}. Take two Artin presentations $p,p'\in\mathcal{P}_n$ and inside the relations $r_j'$ of $p'$ substitute each $x_i'$ with $r_i^{-1}x_ir_i$, for all $i,j=1,\ldots,n$. Denote by $R_j$ the new set of relations produced this way (which are now written in the alphabet $\{x_1^{\pm 1},\ldots, x_n^{\pm 1}\}$). Finally, set the product $p\cdot p'$ to be:
\[p\cdot p'=\langle x_1,\ldots, x_n\;|\;r_1 R_1,\ldots, r_n R_n\rangle.\]

\begin{prop} The above operation is compatible with the group operation of the mapping class group:
\[p(\varphi\circ \psi)=p(\varphi)\cdot p(\psi).\]
\end{prop}

\begin{prf} The operation in the mapping class group (i.e. the composition of diffeomorphisms) is to be read from right to left, as usual. Before doing the computation, notice that from the definition of the $i^{th}$ relation $r_i(\psi)$ associated to $\psi$, we have:
\begin{equation}\label{condition1}
\psi(s_i^{-1})=s_i^{-1}\ast r_i(\psi).
\end{equation}
Moreover, since $x_j$ can be written as $s_j\ast \partial_j\ast s_j^{-1}$, with $\partial_j$ being the $j^{th}$ inner boundary component run around in the clockwise direction, we get that:
\begin{equation}\label{condition2}
\varphi(x_j)=r_j^{-1}(\varphi)\cdot x_j\cdot r_j(\varphi)
\end{equation}
Then we compute:
\begin{align*}
r_i(\varphi\circ \psi)= & s_i\ast (\varphi\circ \psi)(s_i^{-1}) \\
 = &  s_i \ast \varphi(s_i^{-1})\ast \varphi (r_i(\psi)), && \mbox{using Equality \eqref{condition1}} \\
 = & r_i(\varphi)\cdot \varphi(r_i(\psi)), && \mbox{using Equality \eqref{condition2}}\\
 = & r_i(\varphi)\cdot R_i\\
 = & (r(\varphi)\cdot r(\psi))_i.
\end{align*}
\end{prf}
This proposition will be used in the coming section, where we construct a pair $(\Sigma,\varphi)$ from a given Artin presentation $p$, with the property that $p(\varphi)$ is again $p$.

\section{From \texorpdfstring{$p$}{p}  to \texorpdfstring{$(\Sigma,\varphi)$}{(S,f)} }

\begin{lem}\label{artindt} 
The Artin presentation associated to the (positive) Dehn twist around the $i^{th}$ boundary-parallel inner component is:
\[p(\t_{\p_i})=\langle x_1,\ldots, x_n\;|\; 1,\ldots,1,x_i,1,\ldots,1\rangle.\]
\end{lem}

\begin{prf} 
The calculation follows from Figure \ref{artindehn}: note that the arc $s_j$ is disjoint from the curve $\partial_i$ for $j\neq i$, hence $\tau_{\partial_i}(s_j)=s_j$ and $r_j=s_j\ast s_j^{-1}=1$. On the other hand, the relation $r_i$ is computed by concatenating the arc $s_i$ with the path starting at $s_i(1)$, going around the $i^{th}$ hole and going back to $s_i(0)$. This concatenation is precisely the generator $x_i$ in the fundamental group, hence $r_i=x_i$.
\begin{figure}[h!]  
  \centering
 \includegraphics[scale=0.5]{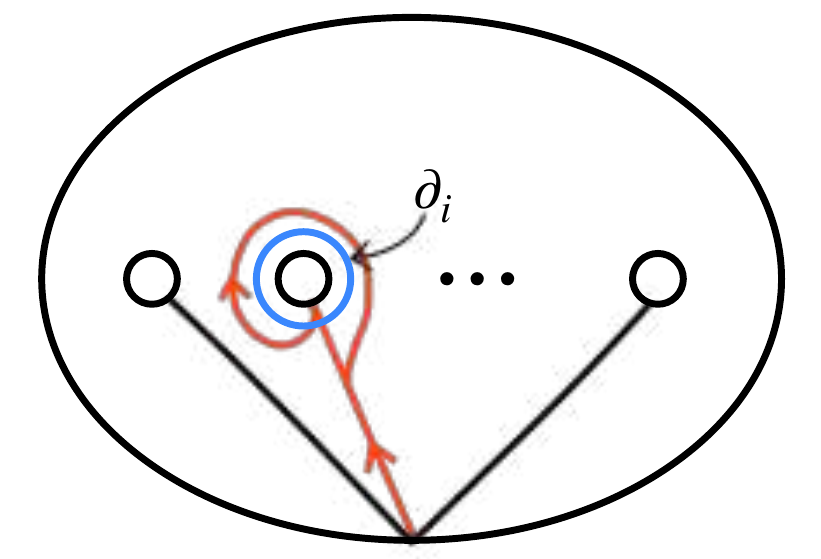}
\caption{}
\label{artindehn}
\end{figure}

\end{prf}

\begin{nrem} Lemma \ref{artindt} can be easily generalized for an arbitrary power $\t_{\p_i}^k$, producing an Artin presentation where the only non-trivial relation is the $i^{th}$ one, namely $r_i(\t_{\p_i}^k)=x_i^k$. Another consequence of Lemma \ref{artindt} is that composing a diffeomorphism $\varphi$ with a such a Dehn twist has the effect of modifying the Artin presentation, according to the operation described in previous section, as follows:
\[p(\t_{\p_i}\circ \varphi)=p(\t_{\p_i})\cdot p(\varphi)=\langle x_1,\ldots, x_n\;|\;r_1(\varphi),\ldots,x_ir_i(\varphi),\ldots, r_n(\varphi)\rangle.\]
This means that every relation except the $i^{th}$ one remains the same, while $r_i(\t_{\p_i}\circ \varphi)$ is simply the concatenation $x_ir_i(\varphi)$.
\end{nrem}

Starting from an Artin presentation $p$ of length $n$, we define a group automorphism 
\[\pi_1(p):F_n\simeq \pi_1(\Sigma_{n+1})\to F_n\simeq \pi_1(\Sigma_{n+1})\]
by sending the generator $x_i$ to $r_i^{-1}x_ir_i$. This uniquely defines the morphism $\pi_1(p)$ and, consequently, a continuous map (unique up to homotopy, which is \emph{not} relative to the boundary)
\[\varphi':\Sigma_{n+1}\to \Sigma_{n+1}\]
which induces the prescribed action on the fundamental group ($x_i\mapsto r_i^{-1}x_ir_i$). This is true because the surface $\Sigma_{n+1}$ is a classifying space for $F_n$. By \cite{artin}, we can actually choose $\varphi'$ to be a diffeomorphism of $\Sigma_{n+1}$, which is the identity on the boundary. We still have to modify the map $\varphi'$ close to the boundary, to find the mapping class of a diffeomorphism whose associated Artin presentation is precisely $p$:

\begin{prop}
There is a unique element $\varphi\in\Gamma(\Sigma_{n+1})$ such that $p(\varphi)=p$.
\end{prop}

\begin{prf} 
We modify $\varphi'$ preserving the fact that the generators $x_i$ are sent to $r_i^{-1}x_ir_i$ for $i=1,\ldots, n$. To this end, notice that the Dehn twists $\t_{\p_1},\ldots, \t_{\p_n}$ act as the identity on $\pi_1(\S,x)$.


Write the relations $r_i$ of $p$ and the relations associated to $\varphi'$ minimally as
\begin{align*}
r_i=& x_i^{a_i}\nu_i && \mbox{ with $\nu_i$ not starting with $x_i^{\pm 1}$,} \\
r_i(\varphi')=& x_i^{b_i}\mu_i && \mbox{ with $\mu_i$ not starting with $x_i^{\pm 1}$.} 
\end{align*}
We impose:
\[r_i^{-1}x_ir_i=r_i(\varphi')^{-1}x_ir_i(\varphi')\in \pi_1(\S,x).\]
Since we are in a free group, the equality $\nu_i^{-1}x_i\nu_i=\mu_i^{-1}x_i\mu_i$ implies $\nu_i=\mu_i$. Hence
\[r_i(\varphi') r_i^{-1}=x_i^{b_i} x_i^{-a_i}= x_i^{b_i-a_i} \Longrightarrow r_i(\varphi')= x_i^{b_i-a_i} r_i.\]
So we see that by composing $\varphi'$ with a suitable power of a boundary-parallel Dehn twist, using Lemma \ref{artindt}, we get:
\[r_i(\t_{\p_i}^{a_i-b_i}\circ \varphi')=r_i(\t_{\p_i}^{a_i-b_i}) r_i(\varphi')=x_i^{a_i-b_i}x_i^{b_i-a_i} r_i=r_i.\]
In this way, we see that the composition 
\[\varphi=\t_{\p_1}^{a_1-b_1}\circ\ldots\circ \t_{\p_n}^{a_n-b_n}\circ \varphi'\]
is the required element in the mapping class group $\Gamma(\Sigma_{n+1})$ with $p(\varphi)=p$. Uniqueness comes from the fact that homotopic diffeomorphisms of a surface are isotopic (\cite[Theorem 6.3]{epstein}), and from the following observation: 
notice that if we take a diffeomorphism $f\in\Gamma(\Sigma_{n+1})$ inducing the identity on $\pi_1(\S,x)$ and cap $n$ boundary components with punctured disks and one with a disk, then the class of $f$ is trivial in the pure mapping class group (see \cite[Theorem 4]{birmanmapping}), hence a product of boundary parallel Dehn twists around those $n$ boundary components (by Birman exact sequence \cite[Theorem 4.6]{farbmarg}). So, if we had picked a map $\varphi''$ different from $\varphi'$ in the first instance, then 
\[\varphi'\circ(\varphi'')^{-1}=\t_{\p_1}^{m_1}\cdots \t_{\p_n}^{m_n},\]
because $\varphi'$ and $\varphi''$ induce the same action on $\pi_1(\S,x)$. Hence correcting $\varphi''$ through the algorithm above would produce the same $\varphi$.
\end{prf}

\section{Matrix of relations and matrix of multiplicities}
The first step to start dealing with Artin presentations through the lens of contact geometry is to identify the \emph{matrix of relations} (Definition \ref{matrixrel}) with the \emph{matrix of multiplicities} (Definition \ref{matrixmult}). 

\begin{defn}\label{matrixrel}
Given an Artin presentation $p\in\mathcal{P}_n$ we define an $n\times n$ matrix $A(p)$ by setting $A(p)_{i,j}$ equal to the total exponent of $x_i$ in the relation $r_j$. The matrix $A(p)$ is called the \emph{matrix of relations} of the Artin presentation $p$.
\end{defn}

Recall that in Section \ref{classificationthm} we defined the capping maps

\begin{minipage}[c]{.40\textwidth}
\centering
\[\xymatrix{
\Gamma(\Sigma) \ar[r]^-{\mbox{cap}} \ar[rd]_-{m(-)} & \Gamma(\Sigma_{0,2}) \ar[d]^-{\simeq} \\
& \Z
}\]
\end{minipage}%
\hspace{10mm}%
\begin{minipage}[c]{.40\textwidth}
\centering
\[\xymatrix{
\Gamma(\Sigma) \ar[r]^-{\mbox{cap}} \ar@/_/[rdd] _-{m(-,-)} & \Gamma(\Sigma_{0,3}) \ar[d]^-{\simeq} \\
& \Z\oplus\Z\oplus\Z \ar[d]^-{\mbox{pr}_3}\\
& \Z
}\]
\end{minipage}

\noindent which are induced by capping off all but one (respectively two) interior component. In the first case we get the mapping class group of the annulus, which is isomorphic to $\Z$, generated by a positive Dehn twist along the core curve. In the second case we get the mapping class group of a pair of pants, isomorphic to a free abelian group of rank 3. By projecting onto the third summand (that takes the role of the outer boundary component) we get the joint multiplicity around the other two components. We denote by $m(-)$ the multiplicity of a single hole, and by $m(-,-)$ the joint multiplicity of a pair of holes.

\begin{defn}\label{matrixmult}
Given a diffeomorpshim class $\varphi\in \Gamma(\S_{n+1})$ of a planar surface, we define the \emph{matrix of multiplicities} $A(\varphi)$ as an $n\times n$ matrix by setting $A(\varphi)_{i,j}=m(i,j)$ for $i\neq j$, and $A(\varphi)_{i,i}=m(i)$.
\end{defn}

In Section \ref{secmatrixmult} (Step 2) we showed how to calculate the matrix of multiplicities starting from a factorization of $\varphi$ as a product of positive and negative Dehn twists.

\begin{thm}\label{mul=rel}
Given $\varphi\in \Gamma(\S_{n+1})$, consider the associated Artin presentation $p(\varphi)\in\mathcal{P}_n$ and the matrices $A(p(\varphi))$ and  $A(\varphi)$. Then:
\[A(p(\varphi))=A(\varphi).\]
\end{thm}

\begin{prf} By \cite{winkelnkemper} we have that, given two Artin presentations $p,p'\in\mathcal{P}_n$, the matrices of relations satisfy the following equality:
\begin{equation}\label{matrelcomp}
A(p\cdot p')= A(p)+A(p').
\end{equation}
Similarly, given two diffeomorphisms $\varphi,\varphi'\in\Gamma(\S)$ it follows from Definition \ref{matrixmult} that 
\begin{equation}\label{matmulcomp}
A(\varphi\circ \varphi')=A(\varphi)+A(\varphi').
\end{equation}
Therefore, thanks to the Equalities \eqref{matrelcomp} and \eqref{matmulcomp}, it is enough to prove the theorem in the case when $\varphi$ is a Dehn twist (around an arbitrary curve $\gamma$). Notice that $A(p(\varphi^{-1}))=-A(p(\varphi))$, from 
\[A(p(\varphi^{-1}))+A(p(\varphi))=A(p(\varphi^{-1})\cdot p(\varphi))=A(p(\varphi^{-1}\circ\varphi))=A(p(\id))= \text{\large 0}_{n\times n}. \]
Further, we reduce to the case when $\gamma$ is a curve as the ones drawn in Figure \ref{gammast}, that we refered to as \emph{standard}. The curves $\gamma_{st}$ are realized as the union of the green arcs together with a choice, at every inner boundary component, of an orange or purple arc (for a total of $2^n$ standard curves). This way we see that for any simple closed curve on the surface $\S_{n+1}$, there is a standard representative in the same homology class.
\begin{figure}[h!]  
  \centering
 \includegraphics[scale=0.6]{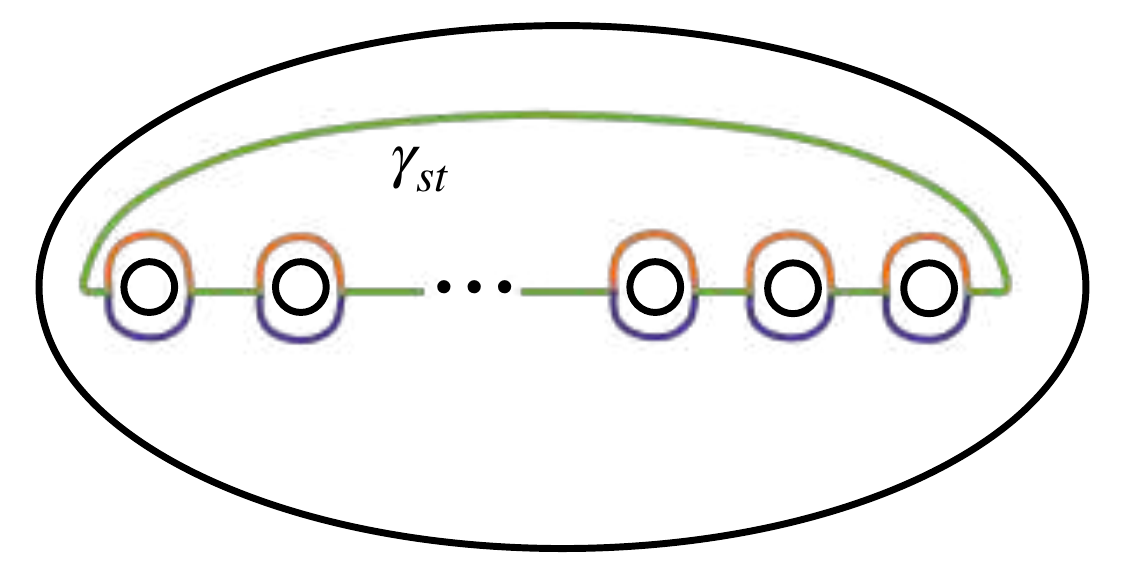}
\caption{Standard curves on $\S_{n+1}$.}
\label{gammast}
\end{figure}

\noindent Notice that $A(p(\t_{\gamma}))$ and $A(\t_{\gamma})$ are independent of the free homotopy type of $\gamma$, they just depends on the homology class of $\gamma$, as the following argument shows: we can change coordinates through a diffeomorphism $f$ and get
\begin{align*}
A(p(\t_{\gamma})) = & A(p(f\t_{\gamma_{st}} f^{-1})) \\
= & A(p(f)p(\t_{\gamma_{st}})p(f^{-1})) \\
= & A(p(f))+A(p(\t_{\gamma_{st}}))+A(p(f^{-1})) \\
= & A(p(\t_{\gamma_{st}})).\\
\\
A(\t_{\gamma}) = & A(f\t_{\gamma_{st}} f^{-1}) \\
= & A(f)+A(\t_{\gamma_{st}})+A(f^{-1}) \\
= & A(\t_{\gamma_{st}}).
\end{align*}

\noindent Hence, it remains to check that the theorem holds for $\varphi=\t_{\gamma_{st}}$.

We calculate $A(p(\t_{\gamma_{st}}))$ in one specific case (Figure \ref{sigma135above}). This has nothing special, but it makes the general argument clear.

\begin{figure}[h!]  
  \centering
 \includegraphics[scale=0.6]{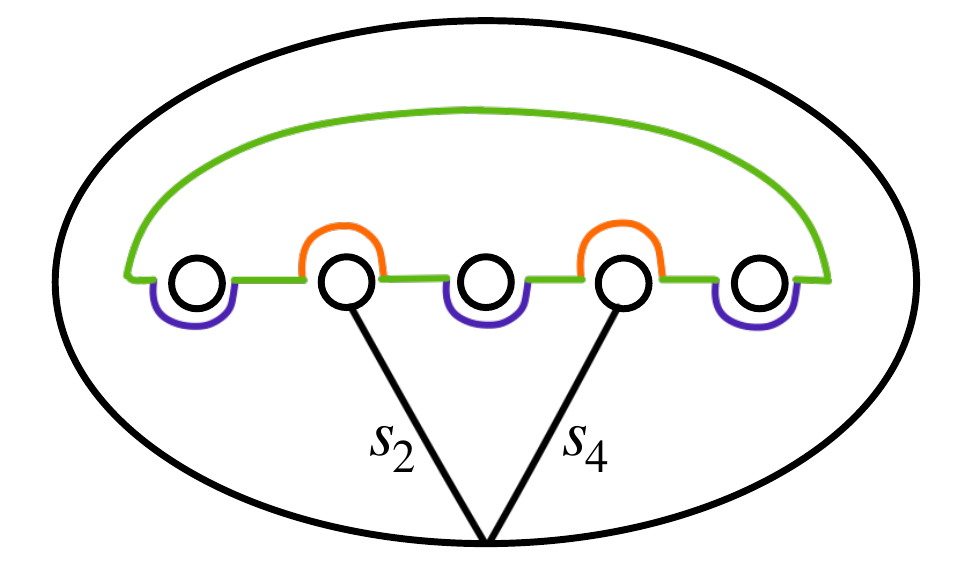}
\caption{}
\label{sigma135above}
\end{figure}

To this end, we start by noticing that $r_2(\t_{\gamma_{st}})=r_4(\t_{\gamma_{st}})=1$, being the curve disjoint from the arcs $s_2$ and $s_4$. Now, since we are only interested in the matrix of relations $A(p(\t_{\gamma_{st}}))$ we change the homotopy representative of $\gamma_{st}$ and keep its homology class fixed: this result in the new curve (see Figure \ref{sigma135below}) which we use to determine the three relations $r_1$, $r_3$ and $r_5$.

\begin{figure}[h!]  
  \centering
 \includegraphics[scale=0.6]{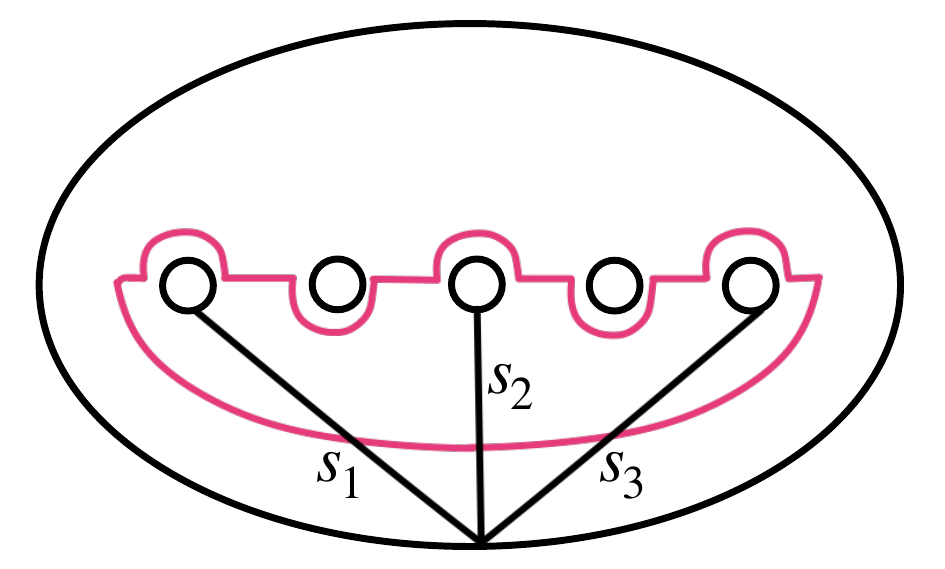}
\caption{}
\label{sigma135below}
\end{figure}

Then one computes $r_1=r_3=r_5=x_1x_3x_5$. These are \emph{not} the relations of $p(\t_{\gamma_{st}})$, but the exponential count of generators appearing there is the same as the entries of $A(p(\t_{\gamma_{st}}))$. From Figure \ref{sigma135above} it is easy to calculate the matrix of multiplicities $A(\t_{\gamma_{st}})$, and therefore we conclude:

\[A(p(\t_{\gamma_{st}}))=
 \begin{bmatrix}
  1 & 0 & 1 & 0 & 1 \\
  0 & 0 & 0 & 0 & 0 \\
    1 & 0 & 1 & 0 & 1 \\
  0 & 0 & 0 & 0 & 0 \\
    1 & 0 & 1 & 0 & 1 
       \end{bmatrix}
\overset{\checkmark}{=}A(\t_{\gamma_{st}}).\]
\end{prf}

\begin{cor}
The matrix of relations $A(p)$ is symmetric for any Artin presentation $p$.
\end{cor}

\begin{prf}
Take a diffeomorphism $\varphi\in\Gamma(\S)$ with $p(\varphi)=p$. By Theorem \ref{mul=rel}, we have 
\[A(p)=A(p(\varphi))=A(\varphi)\] 
and $A(\varphi)$ is symmetric.
\end{prf}

\section{Connection with contact geometry}\label{artingeom}
The bridge to contact geometry is Giroux correspondence \cite{giroux}, which produces a contact 3-manifold $(Y,\xi)$ out of an open book decomposition $(\Sigma,\varphi)$. Since the theory of Artin presentations deals with planar surfaces, we always get planar contact structures. One might be interested in understanding certain properties of $\xi$, such as fillability, through the Artin presentation itself.

A planar contact structure $\xi$ on a 3-manifold $Y$ is Stein fillable if and only the monodromy of a compatible (planar) open book decomposition admits a factorization into positive Dehn twists, as proved by \cite{wendl}. This can be translated in the language of Artin presentation: 

\vspace{0.5cm}

\minibox[frame]{The Artin presentation $p\in \mathcal{P}_n$ corresponds to a Stein fillable contact 3-manifold if and \\ only if $p$ is \emph{quasi-positive}, i.e. it can be written as a product of conjugated of $p(\t_{\gamma_{st}})$.}

\vspace{0.5cm}

\noindent Given the datum of an Artin presentation $p$, we can test quickly whether it can correspond to a Stein fillable contact 3-manifold in the following way. 

\begin{prop} 
Let $p$ be an Artin presentation that determines the contact 3-manifold $(Y,\xi)$ and assume that $(Y,\xi)$ is Stein fillable. Then
\[A(p)_{i,i}\geq A(p)_{i,j} \geq 0,\; \forall i,j.\]
\end{prop}

\begin{prf}
Remember that $p$ can be realized as $p(\varphi)$ for a unique element $\varphi\in\Gamma(\S)$. The contact 3-manifold $(Y,\xi)$ is Stein fillable if and only if $\varphi$ admits a factorization into \emph{positive} Dehn twists. If this happens, then the matrix of multiplicities $A(\varphi)$ would be a sum of $A(\t_{\gamma_{st}})$ for some standard curves. But such matrices have the property that their entries are all non-negative and that the diagonal elements are the biggest among the ones on that same row (and column), and this property is stable under the composition of homeomorphisms, which in turns translate to a sum of the corresponding matrices. And we know by Theorem \ref{mul=rel} that $A(p)=A(\varphi)$. This concludes the proof.
\end{prf}

\begin{exmp}
The result of \cite[Theorem 1.3]{conway} is an example of an infinite family of planar contact 3-manifolds which are hyperbolic, universally tight, with vanishing Heegaard Floer contact invariant but not fillable. To rule out the fillability property the authors use two theorems of \cite{niederkruger} and \cite{baldwin} which deal with holomorphic curves and Heegaard-Floer theory. Our tools represent a low-tech solution that gives the same conclusion: the family of monodromies that the authors deal with in \cite{conway} is 
\[ \varphi=\t_{\alpha}^{-n_1-1}\t_{\beta}^p\t_{B_1}^{n_1} \t_{B_2}^{n_2}\t_{B_3}^{n_3}\t_{B_4}^{n_4}\in \Gamma(\S_{4}),\]
with the supporting curves as in Figure \ref{conwayexmp}.

\begin{figure}[h!]  
  \centering
 \includegraphics[scale=0.6]{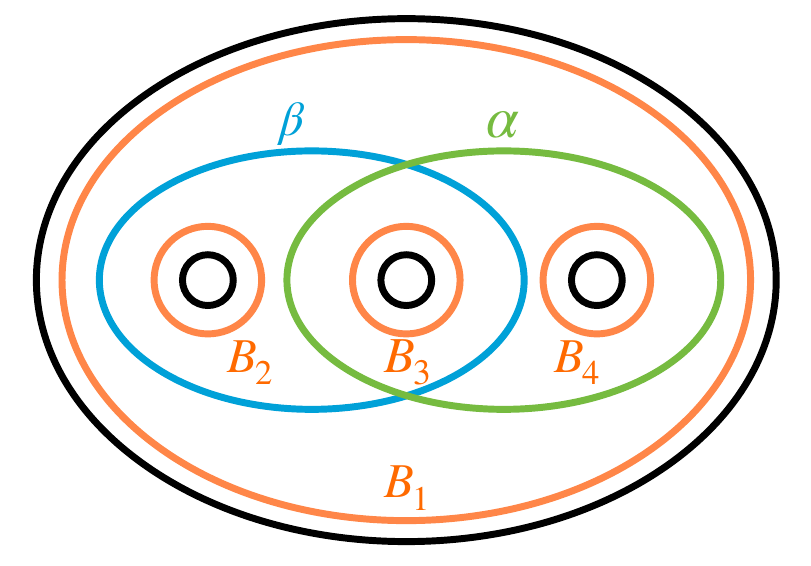}
\caption{Curves configuration.}
\label{conwayexmp}
\end{figure}

We compute the matrix of multiplicities $A(\varphi)$ by looking at Figure \ref{conwayexmp} and at the definition of $\varphi$. For example, $A(\varphi)_{1,1}$ is the sum of the multiplicities of the curves $B_2$, $\beta$ and $B_1$, which are respectively $n_2$, $p$ and $n_1$. We check that the matrix does not satisfy the necessary condition for being fillable (the element in position $(2,3)$ is negative):
\[ A(\varphi)=
\begin{bmatrix}
n_1+n_2+p & p+n_1 & n_1 \\
p+n_1 & n_3+p-1 & \red{-1} \\
n_1 & \red{-1} & n_4-1
 \end{bmatrix}.
\]
\end{exmp}

\begin{exmp}
In \cite{winkelnkemper}, Winkelnkemper proves that the only Artin presentations of length 2 are of the form
\[ p_{a,b,c}=\langle x_1,x_2\;|\;x_1^a(x_1x_2)^{a+c}, x_2^b(x_1x_2)^{b+c}\rangle,\]
with $a,b,c\in \Z$. This is consistent with the fact that the mapping class group $\Gamma(\S_3)$ is isomorphic to $\Z^3$, generated by the three Dehn twists around the boundary components, see Figure \ref{sigma3}. 

\begin{figure}[h!]  
  \centering
 \includegraphics[scale=0.6]{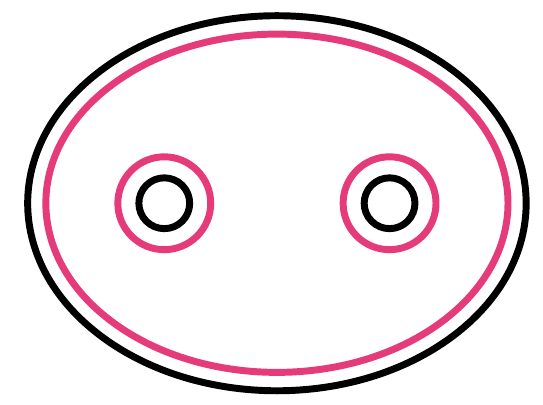}
\caption{Curves generating $\Gamma(\S_3)$.}
\label{sigma3}
\end{figure}

We see also that 
\[p_{a,b,c}=p(\t_{\p_0}^c\t_{\p_1}^a\t_{\p_2}^b)\]
and that 
\[A(p_{a,b,c})=A(\t_{\p_0}^c\t_{\p_1}^a\t_{\p_2}^b)=
\begin{bmatrix}
a+c & c \\
c & b+c
\end{bmatrix}.
\]
It is therefore clear that $p_{a,b,c}$ corresponds to a Stein fillable contact 3-manifold if and only if $a\geq 0$, $b\geq 0$ and $c\geq 0$. One could have noticed that $\t_{\p_0}^c\t_{\p_1}^a\t_{\p_2}^b$ is right-veering if and only if $a,b,c\geq 0$. Since right-veering-ness is a necessary condition for tightness (and hence fillability), this gives us the same conclusion. The advantage of working with Artin presentation is that we just have to deal with matrices instead of drawing curves and arcs. 
\end{exmp}

\newpage

\bibliographystyle{alpha} 

\bibliography{phd_thesis_CORRECTED.bib}
\thispagestyle{plain}

\end{document}